\documentclass[10pt]{amsart}
\usepackage{latexsym}
\usepackage{amssymb}
\usepackage{amsmath}

\numberwithin{equation}{section}

\newcommand{\ben}{\begin{enumerate}}
\newcommand{\een}{\end{enumerate}}

\newcommand{\bec}{\begin{center}}
\newcommand{\eec}{\end{center}}

\newcommand{\beq}{\begin{equation}}
\newcommand{\eeq}{\end{equation}}

\newcommand{\bdm}{\begin{displaymath}}
\newcommand{\edm}{\end{displaymath}}

\newcommand{\R}{\mathbb{R}}

\author{Nikolai  K. Nikolski}
\address{Institut de Math\'ematiques de Bordeaux, Universit\'e Bordeaux I, 
 33405 Talence, France}
\email{nikolski@math.u-bordeaux.fr}

\author{Igor E. Verbitsky}
\address{Department of Mathematics, University of Missouri, Columbia, Missouri 65211, USA}
\email{verbitskyi@missouri.edu}

\thanks{Supported in part
by NSF grant DMS-1161622.}

\begin{document}

\title[Fourier multipliers for weighted $ L^{2}$ spaces]
{Fourier multipliers for weighted $ L^{2}$ spaces\\ with L\'evy-Khinchin-Schoenberg weights}

  \begin{abstract}
 We present a class of weight functions $ w$ on the circle $ \mathbb{T}$,  
called L\'evy-Khinchin-Schoenberg ({\rm LKS}) weights, for which we are 
able to completely characterize (in terms of a capacitary inequality) all 
Fourier multipliers for the weighted space $ L^{2}(\mathbb{T},w)$. We show that the multiplier algebra is nontrivial if and only if $ 
1/w\in L^{1}(\mathbb{T})$, and in this case multipliers satisfy the Spectral 
Localization Property (no ``hidden spectrum"). On the other hand, the Muckenhoupt $ (A_{2})$ condition responsible for the basis property 
of exponentials $ (e^{ikx})$ is more or less independent of the Spectral 
Localization Property and {\rm LKS} 
requirements. Some more complicated compositions of {\rm LKS} weights are 
considered as well. 
\end{abstract} 

\maketitle

\vspace{.2in}

\section{Introduction: Fourier-Hadamard multipliers\\ 
and the Spectral Localization Property (SLP)}

 Given a (nonnegative) finite Borel measure 
$ \mu $ on the circle $ \mathbb{T}={\,} \{z\in \mathbb{C}:{\,} \vert z\vert 
=1\}$, we define Fourier-Hadamard multipliers for the space $ L^{p}(\mathbb{T},\mu 
)$, $ 1\leq {\,} p\leq {\,} \infty $, as\ \par  
  \centerline{ \it sequences of complex numbers $ (\lambda_{n})_{n\in \mathbb{Z}}$ such 
that  the map{\,}}
  \centerline{ {\it $T\colon e^{inx}\longmapsto {\,} \lambda_{n}e^{inx}{\,} 
(n\in \mathbb{Z})$ extends to a bounded linear 
operator on}{\,} $L^{p}(\mathbb{T},\mu )$\rm .}
\ \par 
{\,} {\,} \rm The mapping $ T$ defined by a multiplier $ (\lambda_{n})_{n\in 
\mathbb{Z}}$ is also called a multiplier and denoted by $ T={\,} T_{\lambda }$. 
The set of all $ L^{p}(\mathbb{T},\mu )$-multipliers endowed with the 
obvious operator (multiplier) norm is a unital Banach algebra of sequences 
on $ \mathbb{Z}$ denoted by $ {\rm Mult}(L^{p}(\mu ))$. If $ \mu $ is absolutely 
continuous with respect to Lebesgue measure $ m$, and $ \mu ={\,} wm$, 
we denote the corresponding class of multipliers by
 $ {\rm Mult}(L^{p}(w))$. It is clear that $ (\lambda_{n})_{n\in \mathbb{Z}}\in 
{\rm Mult}(L^{p}(\mu ))$ $ \Rightarrow $ $ \sup_{n\in \mathbb{Z}}\vert \lambda_{n}\vert 
\leq {\,} \Vert T_{\lambda }\Vert <{\,} \infty $, so that always 
$$
 {\rm Mult}(L^{p}(\mu))\subset {\,} l^{\infty }(\mathbb{Z}) .
$$

Despite the fact that multipliers play an important role in Fourier 
analysis, the only cases we know where the algebra $ {\rm Mult}(L^{p}(\mu ))$ 
has been characterized explicitly are $ \mu =  m$ and $ p= 1,{\,} 2,{\,} \infty 
$: 
$$ 
{\rm Mult}(L^{2}(m))={\,} l^{\infty }(\mathbb{Z}), \quad  {\rm Mult}(L^{1}(m))={\,} 
{\rm Mult}(L^{\infty }(m))= \mathcal{ FM}(\mathbb{T}) , 
$$
 $ \mathcal{ F}\mu ={\,} (\hat \mu (n))_{n\in \mathbb{Z}}$  being the Fourier 
transform on $ \mathbb{T}$, and $ \mathcal{ M}(\mathbb{T})$ the space of all complex Borel 
measures on $ \mathbb{T}$. 

Clearly, $\lambda= (\lambda_{n})_{n\in \mathbb{Z}} \in {\rm Mult}(L^{p}(\mu ))$ if and only if 
$\lambda =  \mathcal{ F} k$, where $k$  
is a pseudo-measure on $\mathbb{T}$ such that the corresponding convolution operator $ T_{\lambda }f= 
k \star f$
is bounded: there exists a positive constant $C$ for which  
$$
||k \star f||_{L^p(\mu)} \le C ||f||_{L^p(\mu)}, 
$$
for all trigonometric polynomials $f$. 

Many sufficient conditions are known for a 
sequence $ (\lambda_{n})_{n\in \mathbb{Z}}$ to be a multiplier (mostly 
for the ``flat" case $ \mu ={\,} m$, or for $ \mu ={\,} wm$, where 
$ w$ satisfies the Muckenhoupt condition $ (A_{p})$), starting with 
the famous theorems of J. Marcinkiewicz, S. Mikhlin, L. H\"ormander, S. Stechkin, 
followed by their improvements via the Littlewood-Paley theory, etc.; 
see \cite{Duo2001}, \cite{Gra2008}, \cite{Tor1986}.\ \par 
{\,} {\,} {\,} However, the known results leave open many 
natural questions, in particular, 
\ \par 
  \centerline{ {\it what does the spectrum of a multiplier $ T_{\lambda },$ 
$ \lambda ={\,} (\lambda_{n})_{n\in \mathbb{Z}},$   look like?} 
}
\ \par 
This question is important for solving convolution equations, spectral theory 
of discrete operators of Schr\"{o}dinger type, etc. 
Of course, obviously by definition the eigenvalues $ (\lambda_{n})$ 
are in the spectrum $\sigma(T_{\lambda })$ of $ T_{\lambda }$, and a natural conjecture can 
be, 
{\it whether we have  
$$\sigma (T_{\lambda })= {\rm clos}\Big \{\lambda_{j}:{\,} 
j\in \mathbb{Z}\Big \} , 
$$
for an arbitrary $ T_{\lambda }\in {\rm Mult}(L^{p}(\mu ))$?} If this 
is true for every $ T_{\lambda }\in {\rm Mult}(L^{p}(\mu ))$, we say, following 
\cite{Nik2009}, that the {\it Spectral Localization Property} (SLP) holds.   
Clearly, the SLP is equivalent to the following {\it inverse closedness 
property}
$$ 
\Big (T_{\lambda }\in {\rm Mult}(L^{p}(\mu )), \Big \vert \lambda 
_{j}\Big \vert \geq {\,} \delta >0, \forall j\Big ){\,} \Rightarrow {\,} 
T_{\lambda }^{-1{\,} } is{\,} \textit{bounded} \, \, (\textit{i.e., is in}{\,\, } {\rm Mult}(L^{p}(\mu 
) ).
$$
 
It is known that the algebras $ {\rm Mult}(L^{p}(m))$, 
$ p\not= 2$, {\it do not have the {\rm SLP}}. (For $ p=1$ and $ {\rm Mult}(L^{1}(m))=$ 
$ \mathcal{ FM}(\mathbb{T})$ this is the so-called {\it Wiener-Pitt-Shreider 
phenomenon}, see \cite{GRS1960}, and for $ 1<p<\infty $, $ p\not= 2$, 
 its generalization to $ L^{p}(m)$ spaces due to S. Igari and M. Zafran, 
see \cite{GMcG1979}.)  In this paper, we give nontrivial examples of 
algebras $ {\rm Mult}(L^{2}(w))$ ($ w$ is not equivalent to a constant) satisfying 
the SLP.\ \par 
{\,} {\,} It is also important to know whether there exists an 
estimate for $ \Vert T_{\lambda }^{-1}\Vert $ in terms of the lower 
spectral parameter
$$
 \delta _{T}={\,} \inf_{j}\Big \vert \lambda_{j}(T)\Big \vert .
$$ 
The following quantity is responsible for such a property, 
$$ c_{1}(\delta )={\,} c_{1}(\delta ,{\rm Mult}(L^{p}(\mu )))={\,} 
\sup\Big \{\Big \Vert T^{-1}\Big \Vert \colon  T\in {\rm Mult}(L^{p}(\mu )),{\,} 
\Big \Vert T\Big \Vert \leq 1,{\,} \delta _{T}\geq \delta {\,} \Big \} ,
$$
 where $ 0<\delta \leq 1$.
 
 It is known that for some  function systems 
(say, for complex exponentials $ e^{i\lambda x}$, $ \lambda \in \sigma \subset 
\mathbb{C}$ in certain Banach spaces), even if the multiplier algebra 
is inverse closed, it does not imply that we automatically have a norm 
estimate for inverses (i.e., it may happen that $ c_{1}(\delta )= \infty 
$ for some $ \delta >0$);  see \cite{Nik2009} and the references therein. 
However, for the multiplier algebras appearing in this paper, the situation 
is better: the inverse closedness yields an ``automatic" norm estimate 
for inverses (i.e., $ c_{1}(\delta )<$ $ \infty $, $ \forall \delta 
>0$); see, for instance, Lemma 2.2 below.
\ \par 
{\,} {\,} {\,} In the present paper, we limit ourselves 
to the Hilbert space case, $ p={\,} 2$, and $ \mu ={\,} wm$ (except for  
a few general remarks). In fact, the (open) problem of the spectral 
localization property was the main motivation for the present study. 

Recall that, in general, if a bounded operator $ T\colon H\longrightarrow H$ 
on a Hilbert space $ H$ has a Riesz (unconditional) basis $ (e_{j})$ 
of eigenvectors, $ Te_{j}={\,} \lambda_{j}e_{j}$ ($ j\in J$), then, 
of course, the spectral localization property holds:  $ \sigma (T)={\,} {\rm clos}\{\lambda 
_{j}\colon  j\in J\}$. One could hope that if we replace ``the Riesz basis" 
by ``the (Schauder) basis", then the SLP would still be true. At least, the SLP holds 
for multipliers defining the basis property: $ (e_{j})_{j\geq 1}$ is 
a (Schauder) basis if and only if every sequence $ (\lambda_{j})$ 
of bounded variation $ \sum _{j}\vert \lambda_{j}-\lambda_{j+1}\vert <{\,} 
\infty $ is a multiplier, and if such a sequence is separated from 
zero, $ \inf_{j}\vert \lambda_{j}\vert >{\,} 0$, then the inverse $ (1/\lambda 
_{j})$ is again of bounded variation, and hence a multiplier. 

However, in general, this is not the case: given a complex number $ \alpha $, $ \vert \alpha 
\vert =1$ (not a root of unity), there exists a Muckenhoupt weight 
$ w\in (A_{2})$ such that $ (\alpha ^{j})_{j\in \mathbb{Z}} \in {\rm Mult}(L^{2}(w))$ 
but $ \sigma (T_{\alpha })={\,} \overline{\mathbb{D}}$ (the closed 
unit disc), in particular $ (1/\alpha ^{j})_{j\in \mathbb{Z}}$ is not a multiplier; see 
\cite{Nik2009}.  We show below (Theorem 5.13)  that  the existence of the hidden spectrum $ \sigma (T_{\lambda 
})\setminus\text{clos} \{\lambda_{j}:j\in \mathbb{Z}\}$ in such examples 
is caused by a kind of ``forced holomorphic extension'' of the function 
$ j\longmapsto \lambda_{j}$.
\ \par 
{\,} {\,} For the main class of weights $ w$ considered in this 
paper, namely, the \it ``L\'evy-Khinchin-Schoenberg weights" \rm ({\rm LKS}, for 
short) described below, we will see that the following alternative holds: 
either such a weight $ w\in L^{1}(\mathbb{T})$ 
satisfies the integrability condition $ 1/w\in L^{1}(\mathbb{T})$, and 
then $\mathcal{S}_0 \subset {\rm Mult}(L^{2}(w))$ (here $\mathcal{S}_0$ stands for finitely supported sequences), 
or $ 1/w\not\in L^{1}(\mathbb{T})$, and then ${\rm dim} \, {\rm Mult}(L^{2}(w))< \infty $ (and, in fact, $ {\rm Mult}(L^{2}(w))=\Big\{{\rm const}\Big\}$ in the generic case); 
in both cases the {\rm SLP} holds for ${\rm Mult}(L^{2}(w))$.

The Muckenhoupt condition $ w\in (A_{2})$ (and consequently the fact 
that the exponentials $ (e^{ijx})_{j\in \mathbb{Z}}$ form a Schauder 
basis in $ L^{2}(w)$) plays no essential role for the {\rm SLP}: we will see 
that {\rm LKS} weights satisfying $ 1/w\in L^{1}(\mathbb{T})$ can obey ($ w\in 
(A_{2})$), or disobey ($ w\not\in (A_{2})$) the Muckenhoupt condition, and still 
have the {\rm SLP}. 

Speaking informally, our main message regarding the 
{\rm SLP} is the following. We say that a point $ \zeta\in \mathbb{T}$ is a 
\textit{singularity} of a weight $ w$ if there is no neighborhood $ V$ 
of $ \zeta$ such that $ 0< \inf_{V}w\leq  \sup_{V}w< \infty 
$; then, our results show that\ \par 
 - the {\rm SLP} holds for weights $ w$ having a finite set of singularities 
and ``behaving well" (monotone, or slightly better, see Comment 3.10) 
at every singular point;\ \par 
 - the {\rm SLP} may fail if $ w$ has infinitely many singularities 
(see Theorem 5.13).
 
{\,} {\,} {\,} In Section \ref{sec2}, we develop a kind of general 
scheme to treat the multipliers for ``difference defined" Besov-Dirichlet 
spaces. Of course, it is largely inspired by the famous Beurling-Deny potential 
theory  \cite{BeD1958}, \cite{Den1970}, but there are some new details 
for the case of the discrete group $ \mathbb{Z}$ that we consider. In particular, these spaces are defined by a symmetric 
matrix $\mathcal{C}= (c_{j,k})$, $ c_{j,k}\geq 0$, in such a 
way that their basic properties hold for an arbitrary such matrix (in 
particular, the SLP), and more specific ones - for $\mathcal{C}$ satisfying the 
so-called ``non-splitting condition." In our principal applications 
(Sections \ref{sec3} and \ref{sec4}), $\mathcal{C}$ is a Toeplitz matrix $ (c_{\vert j-k\vert })$, 
and the most complete information is obtained for power-like sequences 
$ (c_{n})$ (the Riesz potential spaces).\ \par 
{\,} {\,} {\,} In Section \ref{sec3}, following P. L\'evy and A. Khinchin 
(and many others, in particular I. Schoenberg, J. von Neumann, M. G. Krein, et al.) we 
introduce a class of remarkable weights for which we will be able to describe 
all multipliers.\ \par 
{\,} {\,} {\,} In Section \ref{sec4}, we complete the program of Section 
\ref{sec3}, giving a capacitary description of multipliers of $ L^{2}(w)$ with 
a L\'evy-Khinchin-Schoenberg ({\rm LKS}) weight. We also discuss a simpler characterization of multipliers 
which does not involve capacities, for {\rm LKS} weights $w$ with quasi-metric  property. 
  In particular, this non-capacitary characterization is valid for multipliers of 
  Besov-Dirichlet spaces of fractional order  which correspond to 
 weights $w(e^{i \theta})=|e^{i \theta}-1|^{\alpha}$, $0<\alpha<1$. For such weights the SLP 
 is a discrete analogue of its continuous counterpart due to Devinatz and Hirschman \cite{DH1959} for multipliers on the group $\mathbb{R}$. 
\ \par 
{\,} {\,} {\,} Section \ref{sec5} is concerned with certain non-{\rm LKS} weights $w$ which can be represented 
as products, or sums of reciprocals, of 
{\rm LKS} weights. Such weights, with a  
finite set of singularities $\zeta_k=e^{i \theta_k}$ ($k=1,2, \ldots, N$) on $\mathbb{T}$ are no longer associated with spaces of Besov-Dirichlet type. Nevertheless, we will show that
 the class of multipliers $ {\rm Mult}(L^{2}(w))$  
permits a complete description in terms of embedding theorems similar 
to those of Sections \ref{sec2}-\ref{sec4}, and has the {\rm SLP}.  In the special case $$
w(e^{i \theta}) = \sum_{k=1}^{d} \,  a_k |e^{i \theta}- \zeta_k|^{-\alpha},  \quad a_k>0, \quad 0<\alpha<1 ,
$$ 
the characterization of $ {\rm Mult}(L^{2}(w))$ depends on the geometry of the points $\{\zeta_k\}$. (A similar characterization holds by duality for $w = \prod_{k=1}^d  |e^{i \theta}- \zeta_k|^{\alpha}$.) 

In particular, if $d$ is a prime 
number, then either $ {\rm Mult}(L^{2}(w))$ coincides with $ {\rm Mult}(L^{2}(w_\alpha))$, where $w_\alpha=|e^{i \theta}-1|^\alpha$, provided $\{\zeta_k\}$ 
are not the set of vertices of a regular polygon, or otherwise with $ {\rm Mult}(L^{2}(w_\alpha (e^{id \theta})))$ where $w_\alpha (e^{id \theta})=|e^{i d \theta}-1|^\alpha$ is equivalent to  an  {\rm LKS} weight with zeros  
at the roots of unity of order $d$. If $d$ is not a prime number the answer is more complicated;   it depends on the divisors of $d$ and involves ``aliases'' of regular polygons (see Theorem 5.1). 

For weights of this type 
with \textit{infinitely} many singularities,
$$
w(e^{i \theta}) = \sum_{k=1}^\infty a_k |e^{i \theta}- \zeta_k |^{-\alpha}, 
$$
where $a_k > 0$, $\sum_{k=1}^\infty a_k < \infty$, and $0<\alpha<1$, it was shown in \cite{Nik2009} that 
the {\rm SLP} may actually fail. 

In this case $ w=w_{-\alpha }\star \nu $, where $ \nu = \sum _{k=1}^\infty  a_{k}\delta _{\zeta_k}$, and  $ {\rm Mult}(L^{2}(\nu ))\subset  {\rm Mult}(L^{2}(\nu \star w_{-\alpha }))$. In Section 5, we complete 
these results of \cite{Nik2009}, by giving a description of ${\rm Mult}(L^{2}(\nu))$ 
in order to show, as mentioned above, that the nature of the hidden 
spectrum of a multiplier $ \lambda =(\lambda_{j})$ lies in 
a  ``forced holomorphic extension'' of the symbol $ j\longmapsto \lambda_{j}$, $j \in \mathbb{Z}$ 
(Theorem 5.13).

\section{Discrete Besov-Dirichlet spaces}\label{sec2} 

 In this section, we work with sequence spaces on 
$ \mathbb{Z}$ (having in mind $ \mathcal{ F}L^{2}(\mathbb{T},w)$ with $ 1/w\in 
L^{1}(\mathbb{T})$, see Section 3 below). An obvious key observation 
is that multipliers of a ``difference defined space" always obey  
the {\rm SLP} (see Lemma 2.2 below). Given a matrix $ \mathcal{C}={\,} (c_{j,k})$\it , 
$ c_{j,k}=$ $ c_{k,j}\geq $ $ 0$\rm , $ c_{j,j}={\,} 0$ ($ j,k\in \mathbb{Z}$) 
and an exponent $ p$, $ 1\leq p<\infty $, we define a (little discrete) 
\it Besov-Dirichlet space $ \mathcal{ B}_{0}^{p}(c_{j,k})$ \rm on $ \mathbb{Z}$ 
in two steps: first, set
$$ \mathcal{ B}^{p}(c_{j,k})={\,} \Big \{x=(x_{j})_{j\in 
\mathbb{Z}}:{\,} \Big \Vert x\Big \Vert ^{p}={\,} \Big \Vert x\Big \Vert 
^{p}_{\mathcal{ B}^{p}  (\mathcal{ C})}={\,} \displaystyle \sum _{j,k}c_{j,k}\Big \vert 
x_{j}-x_{k}\Big \vert ^{p}<{\,} \infty \Big \}$$
equipped with the corresponding (semi)norm $ \Vert \cdot \Vert $. 
A special case important for applications (see Section \ref{sec3}) corresponds  
to $ p={\,} 2$ and $ c_{j,k}={\,} \vert j-k\vert ^{-(1+\alpha )}$, 
$ 0<\alpha <1$. Explaining the terminology, recall a continuous prototype 
of this space, namely, the homogeneous Besov space $ B^{p,p}_{\alpha }(\mathbb{R})$ 
corresponding to the norm  
$$
 \Big \Vert f\Big \Vert ^{p}={\,} \displaystyle \int _{\mathbb{R}}\displaystyle \int 
_{\mathbb{R}}\left ({\frac{\displaystyle \vert f(x)-f(y)\vert }{\displaystyle \vert 
x-y\vert ^{\alpha }}} \right )^{p}{\frac{\displaystyle dxdy}{\displaystyle \vert 
x-y\vert }} ,
$$
 as well as the Dirichlet space of holomorphic functions on the unit disc 
$ \mathbb{D}={\,} \{z\in \mathbb{C}:{\,} \vert z\vert <1\}$ defined 
by
$$ \Big \Vert f\Big \Vert ^{2}={\,} \displaystyle \sum _{n\geq 
1}\Big \vert \hat f(n)\Big \vert ^{2}n={\,} \displaystyle \int _{\mathbb{D}}\Big \vert 
f'(z)\Big \vert ^{2}{\frac{\displaystyle dxdy}{\displaystyle \pi }}$$
$$ ={\,} \displaystyle \int _{\mathbb{T}}\displaystyle \int _{\mathbb{T}}\Big \vert 
{\frac{\displaystyle f(z)-f(\zeta)}{\displaystyle z-\zeta}} \Big \vert ^{2}dm(z)dm(\zeta) ,
$$
where $ m$ stands for normalized Lebesgue measure on $ \mathbb{T}$. (The 
preceding expression for the Dirichlet norm is known as \textit{Compton's 
formula} which goes back to  the 1930s).

The celebrated Beurling-Deny theorem \cite{BeD1958} shows that a Hilbert 
space seminorm $ \Vert \cdot \Vert $ is of the form $ \Vert \cdot \Vert 
_{\mathcal{ B}^{2}  (\mathcal{ C})}$ for some matrix $ \mathcal{ C}$ if and only 
if it is contractive for all Lipschitz maps $ \Phi :\mathbb{C}\longrightarrow \mathbb{C}$ 
such that $ \vert \Phi (z)-\Phi (\zeta)\vert \leq {\,} \vert z-\zeta\vert 
$ and $ \Phi (0)={\,} 0$: $ \Vert \Phi (x)\Vert \leq {\,} \Vert x\Vert 
$ for every complex sequence $ x=(x_{j})_{j\in \mathbb{Z}}$, $ \Phi (x)=(\Phi 
(x_{j}))_{j\in \mathbb{Z}}$.\ \par 
\ \par 
{\,} In order to avoid unnecessary complications we often assume 
that\ \par 
\ \par 
  \centerline{ \it the matrix $ \mathcal{ C}={\,} (c_{j,k})$ does not 
split (into two or more blocks)\rm :}
\ \par 
\rm \noindent if $ A\subset \mathbb{Z}$ is such that $ c_{j,k}=$ $ 0$ for every $ j\in 
A$ and $ k\in \mathbb{Z}\backslash A$ then either $ A={\,} \emptyset$ 
or $ A={\,} \mathbb{Z}$.\ \par 
{\,} {\,} We denote by $ e_{n}$ the standard $ 0-1$ sequence, $ e_{n}={\,} 
\mathcal{ F}z^{n}={\,} (\delta _{nj})_{j\in \mathbb{Z}}$, and observe 
that (for a non-splitting $ \mathcal{ C}$) $ \Big \Vert e_{n}\Big \Vert ^{p}={\,} 
2\displaystyle \sum _{j}c_{j,n}>{\,} 0$ ($ \forall n \in  \mathbb{Z}$).\ \par 
\ \par\noindent  
\bf 2.1. Lemma. {\it Assume $\mathcal{C}$ is a non-splitting matrix. Then the following statements 
hold.\ \par \ \par\noindent 
\bf (1) \it One has $ e_{n}\in \mathcal{ B}^{p}(c_{j,k})$ if and only if 
$ \displaystyle \sum _{j}c_{j,n}<$ $ \infty $.}\ \par\noindent 
(2) \it Assume $ \displaystyle \sum _{j}c_{j,n}<$ $ \infty $ for every 
$ n\in \mathbb{Z}$, and let 
$$ 
\mathcal{ S}_{0}={\,} {\rm Lin}(e_{n}:n\in \mathbb{Z})
$$
be a vector space of finitely supported sequences. Then $ \Vert \cdot 
\Vert $ is a norm on $ \mathcal{ S}_{0}$.\ \par \ \par\noindent  
\bf (3) \it If one of the coordinate functionals $ \varphi _{n}:(x_{j})_{j\in 
\mathbb{Z}}\longmapsto $ $ x_{n}$ is bounded on $ \mathcal{ S}_{0}$ (respectively, 
on $ \mathcal{ B}^{p}(c_{j,k})$), then all of them are bounded.\ \par 
\ \par 
{\,} {\,} \bf Proof. \rm Statement (1) is obvious. We prove (3) first. Without 
loss of generality, suppose $ \varphi _{0}$ is bounded. The non-splitting 
hypothesis implies that for every $ n\in \mathbb{Z}$ there exists a sequence 
$ n_{0}=0$, $ n_{1}$,..., $ n_{k}={\,} n$ (called the \it chain \rm joining 
$ 0$ and $ n$) such that $ c_{n_{j}, n_{j+1}}>{\,} 0$, for all $ 
j={\,} 0,...,k-1$. (Indeed, if $ A$ is the set of all $ n\in \mathbb{Z}$ 
joinable to $ 0$, then $ c_{j,k}=$ $ 0$ for every $ j\in A$ and $ 
k\in \mathbb{Z}\backslash A$, and so $ A=\mathbb{Z}$. Following \cite{BeD1958}, 
the existence of such a chain can  also be called ``$\mathcal{C}$ -connectedness 
of $ \mathbb{Z}$.") Hence, for every $ x=(x_{j})_{j\in \mathbb{Z}}\in \mathcal{ B}^{p}$, 
$$ 
\Big \vert x_{n}\Big \vert \leq {\,} \Big \vert x_{0}\Big \vert 
+{\,} \displaystyle \sum _{j=0}^{k-1}\Big \vert x_{n_{j}}-x_{n_{j+1}}\Big \vert 
$$ 
$$ \leq {\,} \Big \vert \varphi _{0}(x)\Big \vert +{\,} 
a\Big (\displaystyle \sum _{j=0}^{k-1}c_{n_{j}  ,n_{j+1}}\Big \vert x_{n_{j}}-x_{n_{j+1}}\Big \vert 
^{p}\Big )^{1/p}\leq {\,} A\Big \Vert x\Big \Vert _{\mathcal{ B}^{p}} , 
$$ 
and so $ \Big \Vert \varphi _{n}\Big \Vert \leq {\,} A\leq {\,} 
\Big \Vert \varphi _{0}\Big \Vert +$ $ \displaystyle \sum _{j=0}^{k-1}(c_{n_{j}  ,n_{j+1}})^{-1/p}$.\ \par 
{\,} {\,} To prove (2), notice that if $ x\in \mathcal{ S}_{0}$, $ x_{j}={\,} 0$ 
for $ \vert j\vert >N$, and $ \Vert x\Vert ={\,} 0$, the same reasoning 
as for (3) gives $ x_{k}={\,} 0$ for every $ k\in \mathbb{Z}$ (fix 
a $ j$ with $ \vert j\vert >N$ and join $ k$ to $ j$ by a chain). 
 \qed \ \par \ \par \noindent  
\bf Remark. \rm For a general (non-zero) symmetric matrix $ \mathcal{C}= (c_{j,k})_{j,k\in \mathbb{Z}}$, $ c_{j,k}\geq  0$, instead of the non-splitting 
hypothesis, we can introduce in $\mathbb{Z}$ an equivalence relation 
$ R$ saying that $ nRm$ if there exists a sequence of integers 
$ j_{1} =n, j_{2},..., j_{s}=m$ such that $ c_{j_{k}, j_{k+1}}>  
0$ for $ j=1, \ldots, s-1$ ($nRn$ for every $ n$ by 
definition). Different cosets of this relation form a partition of 
$\mathbb{Z}$ (finite or not), and the non-splitting hypothesis says that 
there is only one such ``$ C$-connected" component $ E_{0}$ (equal to 
$ \mathbb{Z}$).\ \par 
  These \it cosets are denoted by $ E_{0}, \ldots, E_{k}, \ldots$ . \rm  
It is clear that $ \Vert x\Vert ^{p}= \sum _{k}\Vert x\chi _{E_{k}}\Vert 
^{p}$ for every $ x\in \mathcal{S}_{0}$. The reasoning of Lemma 2.1 shows 
that $ \Vert \cdot \Vert $ \it is a norm on $\mathcal{S}_{0}(E_{k})=  
{\rm Lin} (e_{j}:\, \,  j\in E_{k})$ if and only if $ E_{k}$ is infinite, \rm 
and \ \par 
\ \par 
  \centerline{ $ \Vert \cdot \Vert $ \it is a norm on $\mathcal{S}_{0}$ 
if and only if each $ E_{k}$ is infinite\rm .}
\ \par 
\, \,  \rm In this paper, the most important case 
is the \it Toeplitz matrix case $ c_{j,k}= c_{\vert j-k\vert }$, \rm 
where $ c= (c_{n})_{n\geq 0}$ is a given non-negative (non-zero) 
sequence (see Sections \ref{sec2} and \ref{sec3}). In this case, the principal parameter 
is the number $ D$ defined by 
$$
D= D({ \mathcal{C}}) :=\,  GCD\{k\geq 1: \,  c_{k}>  
0\}.$$ 
If $ p_{1}, \ldots, p_{m}$ are such that $ D(\mathcal{C})= GCD\{p_{j}:   j = 
1, \ldots, m\}$, then the $ C$-connected component $ E_{0}$ containing 
$ n= 0$ consists of the numbers $ \sum _{j=1}^{m}p_{j}n_{j}$ ($ n_{j}\in 
\mathbb{Z}$), and so is $ D \mathbb{Z}$ (B\'ezout's theorem); the other cosets 
are $ j+D \mathbb{Z}$, $ j= 1, \ldots, D-1$. Therefore, \it for a non-zero 
Toeplitz matrix $\mathcal{C}$\rm , $ \Vert \cdot \Vert $ \it is always a norm 
on $\mathcal{S}_{0}$.\ \par \ \par  
{\,} {\,} \it In what follows, we always assume that $\mathcal{S}_{0}\subset 
 \mathcal{B}^{p}(\mathcal{C})$, and $ x\longmapsto  \Vert x\Vert _{\mathcal{B}^{p}  ( \mathcal{C})}$ 
is a norm on $\mathcal{S}_{0}$.  \rm We define a ``little $ \mathcal{ B}^{p}$ space" 
by 
$$ \mathcal{ B}_{0}^{p}(c_{j,k})={\,} {\rm span}_{\mathcal{ B}^{p}  (c_{j,k}  )}\Big (e_{n}\colon{\,} 
n\in \mathbb{Z}\Big ) ,
$$ 
where $ span$ means the ``closed linear span" (or, better,  the {\it completion}  
of $ \mathcal{ S}_{0}$\rm ,$ \Vert \cdot \Vert $), and the {\it multipliers}  
of $ \mathcal{ B}_{0}^{p}(c_{j,k})$ \rm by the following (standard) requirement: 
$$ 
{\rm Mult}(\mathcal{ B}_{0}^{p}(c_{j,k}))={\,} 
$$
$$ 
={\,} \{(\lambda_{j})_{j\in \mathbb{Z}}\colon \, \,  T_{\lambda }\colon e_{n}\longmapsto 
\lambda_{n}e_{n}, \, \,  \forall n \in \mathbb{Z}, {\,} \textit{extends  to a bounded 
map on} \, \, \mathcal{ B}_{0}^{p}(c_{j,k})\}  .
$$ 

It is clear that $ {\rm Mult}(\mathcal{ B}_{0}^{p}(c_{j,k}))$ 
is a (commutative) unital Banach algebra.\ \par 
\ \par \noindent  
\bf 2.2. Lemma.\it The following statements hold.\ \par \ \par \noindent  
\bf (1) \it Let $ \lambda ={\,} (\lambda_{j})_{j\in \mathbb{Z}}$ 
be a sequence of complex numbers. Then 
$$ \lambda \in {\rm Mult}(\mathcal{ B}_{0}^{p}(c_{j,k})) \Leftrightarrow {\,} 
\lambda \in l^{\infty }(\mathbb{Z}){\,\,} \textit{and} {\,\,} \displaystyle \sum _{k}\Big \vert 
x_{k}\Big \vert ^{p}\mu _{k}^{p}\leq {\,} C^{p}\Vert x\Vert ^{p}_{\mathcal{ B}^{p}},  \, \,  
 \forall x\in \mathcal{ B}_{0}^{p}(c_{j,k}) ,
 $$ 
 where $ \mu _{k}^{p}={\,} \mu _{k}^{p}(\lambda )=:{\,} \displaystyle \sum _{j}c_{j,k}\Big \vert 
\lambda_{j}-\lambda_{k}\Big \vert ^{p}$, and $ C$ is a positive 
constant.\ \par 
{\,} {\,} In particular, $ {\rm Mult}(\mathcal{ B}_{0}^{p}(c_{j,k}))$ obeys 
the {\rm SLP}, and 
$$
 \Vert T_{1/\lambda }\Vert \leq  \delta ^{-2}(\Vert T_{\lambda 
}\Vert +\Vert \lambda \Vert _{\infty })+  \delta ^{-1}, \, \, \textit{where}  \, \, \delta 
={\,} \inf_{j}\vert \lambda_{j}\vert >{\,} 0 .$$
\ \par \noindent  
\bf (2) \it If $ C=C(\lambda )$ is the best possible constant in (1), 
then 
$$ C(\lambda )-\Big \Vert \lambda \Big \Vert _{l^{\infty }}\leq 
{\,} \Big \Vert T_{\lambda }\Big \Vert \leq  C(\lambda )+\Big \Vert \lambda 
\Big \Vert _{l^{\infty }} . 
$$
\ \par 
{\,} {\,} {\,} \bf Proof. \rm If $ \lambda \in {\rm Mult}(\mathcal{ B}_{0}^{p}(c_{j,k}))$, then clearly 
$\lambda \in l^{\infty }(\mathbb{Z})$,  since $ \lambda_{j}$ 
are eigenvalues of a bounded operator. Hence, we can assume that $ \lambda \in l^{\infty 
}(\mathbb{Z})$. We have 
$$ \lambda \in {\rm Mult}(\mathcal{ B}_{0}^{p}(c_{j,k})) \Leftrightarrow  
 \displaystyle \sum _{j,k}c_{j,k}\Big \vert \lambda_{j}x_{j}-\lambda_{k}x_{k}\Big \vert 
^{p}\leq a^{p}\Big \Vert x\Big \Vert ^{p}, \, \,  \forall x\in \mathcal{ B}_{0}^{p}(c_{j,k}) ,
$$
 where $ a={\,} \Vert T_{\lambda }\Vert $ for ($\Rightarrow$),  
or $ \Vert T_{\lambda }\Vert \leq {\,} a$ for ($\Leftarrow$). 
Let $ x={\,} (x_{j})_{j\in \mathbb{Z}}$ be a finitely supported sequence.  Then 
$$ 
\sum _{j,k}c_{j,k}\Big \vert \lambda_{j}\Big \vert 
^{p}\Big \vert x_{j}-x_{k}\Big \vert ^{p}\leq {\,} \Big \Vert \lambda \Big \Vert 
_{l^{\infty }}^{p}\Big \Vert x\Big \Vert ^{p} ,$$
 and hence 
 $$ 
 \Big \Vert T_{\lambda }x\Big \Vert ={\,} \Big (\sum _{j,k}c_{j,k}\Big \vert \lambda_{j}x_{j}-\lambda_{j}x_{k}+{\,} \lambda_{j}x_{k}-\lambda_{k}x_{k}\Big \vert ^{p}\Big )^{1/p}
$$
 $$ 
 \leq {\,} \Big \Vert \lambda \Big \Vert _{l^{\infty }}\Big \Vert 
x\Big \Vert {\,} +{\,} \Big (\sum _{j,k}c_{j,k}\Big \vert 
\lambda_{j}-\lambda_{k}\Big \vert ^{p}\Big \vert x_{k}\Big \vert ^{p}\Big )^{1/p}$$
$$ ={\,} \Big \Vert \lambda \Big \Vert _{l^{\infty }}\Big \Vert 
x\Big \Vert +{\,} \Big (\sum _{k}\Big \vert x_{k}\Big \vert ^{p}\mu_{k}^{p})^{1/p} . 
$$
If the right-hand side inequality in (1) holds, we obtain
$$
 \Big \Vert T_{\lambda }x\Big \Vert \leq {\,} \Big (\Big \Vert 
\lambda \Big \Vert _{l^{\infty }}+C\Big )\Big \Vert x\Big \Vert ,
$$ 
which shows that $ \lambda \in {\rm Mult}(\mathcal{ B}_{0}^{p}(c_{j,k}))$ 
and $ \Big \Vert T_{\lambda }\Big \Vert \leq $ $ \Big \Vert \lambda \Big \Vert 
_{l^{\infty }}+C$. Conversely, if $ \lambda $ is a multiplier, we have 
as before, 
$$ \Big (\sum _{k}\Big \vert x_{k}\Big \vert ^{p}\mu_{k}^{p}\Big )^{1/p}\leq {\,} \Big \Vert \lambda \Big \Vert _{l^{\infty }}
\Big \Vert x\Big \Vert +{\,} \Big \Vert T_{\lambda }x\Big \Vert ,
$$
 and so the right hand side inequality follows with $ C\leq {\,} \Big \Vert \lambda 
\Big \Vert _{l^{\infty }}+$ $ \Big \Vert T_{\lambda }\Big \Vert \leq {\,} 
2\Big \Vert T_{\lambda }\Big \Vert $.\ \par 
{\,} {\,} {\,} The {\rm SLP} follows from this 
description of multipliers by means of the  embedding theorem: if 
$ \lambda =(\lambda_{j})_{j\in \mathbb{Z}}\in 
{\rm Mult}(\mathcal{ B}_{0}^{p}(c_{j,k}))$ and $ \delta ={\,} \inf_{j}\vert \lambda 
_{j}\vert >{\,} 0$, then  
$$ \mu _{k}^{p}(1/\lambda )={\,} {\,}  \sum _{j}c_{j,k}{\frac{\Big \vert \lambda_{j}-\lambda_{k}\Big \vert ^{p}}{\Big \vert 
\lambda_{j}\lambda_{k}\Big \vert ^{p}}} \leq {\,} {\frac{1}{\delta^{2p}}} \mu _{k}^{p}(\lambda ) ,
$$ 
and hence $ C(1/\lambda )\leq {\,} C(\lambda )/\delta ^{2}$, 
$ \Vert T_{1/\lambda }\Vert \leq {\,} C(1/\lambda )+\Vert 1/\lambda \Vert 
_{\infty }\leq {\,} \delta ^{-2}(\Vert T_{\lambda }\Vert +\Vert \lambda 
\Vert _{\infty })+{\,} \delta ^{-1}$.\ \par 
{\,} {\,} It is clear that (2) is proved as well.  \qed \ \par 
\ \par \noindent 
\bf Remark. \rm A slightly different estimate of $ \Vert T_{1/\lambda }\Vert 
$ follows by a direct computation: 
$$ \Big \Vert T_{1/\lambda }x\Big \Vert = \Big (\Big \vert x_{0}/\lambda 
_{0}\Big \vert ^{p} + \displaystyle \sum _{j,k}c_{j,k}\Big \vert (x_{j}/\lambda 
_{j})-(x_{k}/\lambda_{k})\Big \vert ^{p}\Big )^{1/p}
$$
$$ =\Big (\Big \vert x_{0}/\lambda_{0}\Big \vert ^{p}+ \displaystyle \sum _{j,k}c_{j,k}\Big \vert {\frac{\displaystyle \lambda 
_{k}(x_{j}-x_{k})-\lambda_{j}(x_{k}-x_{j})+\lambda_{k}x_{k}-\lambda_{j}x_{j}}{\displaystyle 
\lambda_{j}\lambda_{k}}} \Big \vert ^{p}\Big )^{1/p} 
$$
$$ \leq {\,} 2\delta ^{-1}\Big \Vert \lambda \Big \Vert _{\infty 
}\Big \Vert x\Big \Vert +{\,} \delta ^{-2}\Big \Vert T_{\lambda }x\Big \Vert , 
$$
 so $ \Big \Vert T_{1/\lambda }\Big \Vert \leq $ $ 2\delta ^{-1}\Big \Vert 
\lambda \Big \Vert _{\infty }+$ $ \delta ^{-2}\Big \Vert T_{\lambda }\Big \Vert 
$.\ \par 
\ \par 
{\,} {\,} Lemma 2.2 allows us to decide when the multiplier algebra 
is nontrivial, that is, $ {\rm Mult}(\mathcal{ B}_{0}^{p}(c_{j,k}))\not= {\,} \Big \{{\rm const}\Big \}$.\ \par 
\ \par \noindent  
\bf 2.3. Lemma. \it Under the hypothesis of non-splitting, the following statements 
are equivalent.\ \par \ \par \noindent  
\bf (1) $ \mathcal{ S}_{0}\subset {\,} {\rm Mult}(\mathcal{ B}_{0}^{p}(c_{j,k}))$\it .\ \par \ \par \noindent  
\bf (2) $ {\rm Mult}(\mathcal{ B}_{0}^{p}(c_{j,k}))\not= $ $ \Big \{{\rm const}\Big \}$\it .\ \par \ \par \noindent  
\bf (3) \it All $ \varphi_{n}$ (see Lemma 2.1(3)) are bounded 
on $ \mathcal{ S}_{0}$.\ \par \ \par \noindent  
\bf (4) $ (e_{n})_{n\in \mathbb{Z}}$ \it is a minimal sequence in $ \mathcal{ B}_{0}^{p}(c_{j,k})$.\ \par 
\ \par 
{\,} {\,} \bf Proof. \rm Clearly, (1) $\Rightarrow $ (2) and 
(3) $\Rightarrow $ (4) ($ \varphi _{n}(e_{k})={\,} \delta _{nk}$); moreover,  (4) $\Rightarrow $ (3) since a sequence biorthogonal  to $(e_n)$ coincides with $ \varphi_{n}$ on $\mathcal{S}_0$. Hence, (3) $\Leftrightarrow $ (4), and consequently 
(4) $\Rightarrow $ (1) is also easy: if $ \sum _{n}\vert \lambda_{n}\vert 
\cdot \Vert \varphi _{n}\Vert \cdot \Vert e_{n}\Vert <{\,} \infty $, 
a multiplier $ T_{\lambda }x={\,} \sum _{n}\lambda_{n}\varphi _{n}(x)e_{n}$, 
$ x\in \mathcal{ S}_{0}$, is bounded.\ \par 
{\,} {\,} {\,} Let us show (2) $\Rightarrow $ (3): assume if $ \varphi _{0}$ 
is not bounded, so are all $ \varphi _{k}$ (see Lemma 2.1(3)), and 
let $ \lambda \in {\rm Mult}(\mathcal{ B}_{0}^{p}(c_{j,k}))$. By Lemma 2.2, there 
exists a constant $ C$ such that $ \displaystyle \sum _{n}\Big \vert x_{n}\Big \vert 
^{p}\mu _{n}^{p}\leq $ $ C^{p}\Vert x\Vert ^{p}_{\mathcal{ B}^{p}}$ for 
every $ x\in \mathcal{ S}_{0}$. Since $ x_{n}={\,} \varphi _{n}(x)$ 
this implies that $ \mu _{n}^{p}(\lambda ) :=$ $ \displaystyle \sum _{j}c_{j,n}\Big \vert 
\lambda_{j}-\lambda_{n}\Big \vert ^{p}={\,} 0$ for all $ n\in \mathbb{Z}$. 
Now, given $ n\in \mathbb{Z}$, the block splitting hypothesis implies 
the existence of a chain $ n_{0}=0$, $ n_{1}$,..., $ n_{k}=$ $ n$ such 
that $ c_{n_{j}, n_{j+1}}>$ $ 0$ for all $ j=$ $ 0,...,k-1$ (see arguments 
of Lemma 2.1), and hence $ \Big \vert \lambda_{n_{j}}-\lambda_{n_{j+1}}\Big \vert 
^{p}=$ $ 0$ for all $ 0\leq {\,} j\leq {\,} k-1$, which implies 
$ \lambda_{n}={\,} \lambda_{0}$, $ \forall n\in \mathbb{Z}$.  \qed 
\ \par 
\ \par \noindent  
\bf Remark. \rm In general (i.e., without the non-splitting hypothesis), 
properties (1), (3), and (4) are still equivalent. As to (2), using 
the mentioned above partition $\mathbb{Z}= \bigcup_{k\geq 0}E_{k}$ 
into a union of $C$-connected components, it is easy to see that if  $ \varphi _{n}$ 
is bounded for some $ n\in E_{k}$, then all $ \varphi _{j}\in E_{k}$ are bounded, 
and this property is equivalent to the fact that ${\rm Mult}(\mathcal{B}^{p}_{0}(c_{i,j})_{i, j\in 
E_{k}})\not= \Big\{{\rm const}\Big\}$. Conversely, if  $ \varphi _{n}$ is unbounded 
for some $ n\in E_{k}$ (and hence all $ \varphi _{j}\in E_{k}$ 
are unbounded), then every multiplier $ \lambda \in {\rm Mult}(\mathcal{B}^{p}_{0}(c_{i,j}))$ 
is constant on $ E_{k}$ (with the same proof as above). In particular, \it if 
all $ \varphi _{j}$, $ j\in \mathbb{Z}$, are unbounded, then 
${\rm Mult}(\mathcal{B}^{p}_{0}(c_{i,j}))$ 
consists of all $ l^{\infty }(\mathbb{Z})$ sequences constant on every 
$ E_{k}$, $ k=0, 1, \ldots$ .\ \par 
\ \par \noindent  
\bf 2.4. Example. \rm Let $ c_{k,k+1}={\,} c_{k,k-1}={\,} 1$ 
($ \forall k \in \mathbb{Z}$) and $ c_{k,j}={\,} 0$ for all other indices $ (k,j)\in 
\mathbb{Z}^{2}$, so that 
$$ \Big \Vert x\Big \Vert ^{p}={\,} 2\displaystyle \sum _{k\in 
\mathbb{Z}}\Big \vert x_{k}-x_{k+1}\Big \vert ^{p}, \, \,  x\in \mathcal{ S}_{0} .
$$
 Clearly, if $ p>1$, the functional $ \varphi _{0}$ is not bounded 
(and it is, if $ p={\,} 1$), and by Lemma 2.3, $ {\rm Mult}(\mathcal{ B}_{0}^{p}(c_{j,k}))=$ 
$ \Big \{{\rm const}\Big \}$. For $ p=1$, $ \mathcal{ B}_{0}^{1}(c_{j,k})$ is 
the space of sequences tending to zero and of bounded variation; it 
is easy to see that $ {\rm Mult}(\mathcal{ B}_{0}^{1}(c_{j,k}))={\,} \mathcal{ B}^{1}(c_{j,k}))$,  
that is,  the space of all sequences of bounded variation.\ \par 
\ \par   \noindent 
\bf 2.5. Comments. (1) \rm On the definition of $ \mathcal{ B}_{0}^{p}(c_{j,k})$: \it what is the completion of $ \mathcal{ S}_{0}$? \rm It is clear that 
if the coordinate functionals $ \varphi _{n}$ are bounded on $ \mathcal{ S}_{0}$, then 
we can realize the completion of $ \mathcal{ S}_{0}$ (i.e., the space $ \mathcal{ B}_{0}^{p}(c_{j,k})$) 
as a sequence space. Much less clear is what happens beyond this condition. 
For the norm from Example 2.4, we will give a description of $ \mathcal{ B}_{0}^{2}(c_{j,k})$ 
in Section 3.\ \par 
\ \par   \noindent 
\bf (2) \rm In general, beyond the scope of Besov-Dirichlet spaces, 
the lack of minimality of a sequence $ \mathcal{ E}={\,} (e_{n})$ in 
a Banach space $ X$ does not prevent the multiplier algebra 
$$
 {\rm Mult}(\mathcal{ E}) ={\,} 
\{(\lambda_{n}):{\,} e_{n}\longmapsto \lambda_{n}e_{n}, \, \,  \forall 
n, \, {\,}
$$
$$
\textit{extends to a  bounded linear map on}\,\, {\rm span}_{X}(\mathcal{ E})\}
$$ 
to be nontrivial. A standard example 
is given by reproducing kernel sequences in a holomorphic space, say 
the Hardy space $ H^{2}$ on the unit disc $\mathbb{D}$: whatever is a sequence $ e_{n}={\,} 
(1-\overline{w}_{n}z)^{-1}$, $ \vert w_{n}\vert <1$, the sequence $ \lambda 
_{n}={\,} \varphi (\overline{w}_{n})$, where $ \varphi \in H^{\infty }$, 
is in $ {\rm Mult}(\mathcal{ E})$. (If $ (w_{n})$ is not a Blaschke sequence, 
all multipliers are of that form.) Similar examples involve exponentials 
$ e^{i\lambda x}$, ${\rm Im} (\lambda )>0$ in the space $ L^{2}(0,\infty )$.\ \par 
\ \par   \noindent 
\bf (3) \rm It is not clear how to express the minimality property 
of Lemma 2.3 in terms of $ c_{j,k}$. For an important partial case 
(the principal spaces of this paper), where $ c_{j,k}={\,} c_{\vert j-k\vert 
}$ and $ c_{n}\geq 0$, $ 0<{\,} \sum c_{n}<{\,} \infty $, we 
will give certain criteria in Section \ref{sec3}, in different forms. For example, the 
condition 
$$\sum _{n\geq 1}{\frac{1}{ \sum_{k=1}^{n}c_{k}k^{2}+n^{2} \sum_{k>n}c_{k}}} < \infty $$
is necessary, and 
$$\sum _{n\geq 1}{\frac{1}{\sum _{k=1}^{n}c_{k}k^{2}}} <\infty $$ 
is sufficient for the minimality of $ (e_{n})$ in 
the space $ \mathcal{B}_{0}^{2}(c_{\vert j-k\vert })$. In particular, 
for $ c_{k}={\,} (1+\vert k\vert )^{-(1+\alpha )}$ ($ \alpha>0$), the 
minimality holds if and only if $ \alpha <1$.\ \par 
\ \par   \noindent 
\bf (4) \rm One more property of spaces $ \mathcal{ B}_{0}^{2}(c_{j,k})$ that  
we will need in Section \ref{sec3} is the uniqueness of the matrix $ \mathcal{ C}={\,} 
(c_{j,k})$ in the definition of a Besov-Dirichlet norm, as  stated 
in the following lemma.\ \par 
\ \par  \noindent 
\bf 2.6. Lemma. \it Let $ \mathcal{ C}=$ $ (c_{j,k})$ and $ \mathcal{ C}'=$ 
$ (c'_{j,k})$ be two matrices satisfying the above conditions (i.e., non-splitting, and 
$e_n \in \mathcal{ B}_{0}^{2}$, $\forall \, n \in \mathbb{Z}$; see 
Sec. 2.1), which define the same 2-norm: $ \Big \Vert x\Big \Vert _{\mathcal{ B}^{2}  (\mathcal{ C})}={\,} 
\Big \Vert x\Big \Vert _{\mathcal{ B}^{2}  (\mathcal{ C}')}$ for every $ x\in 
\mathcal{ S}_{0}$. Then $ \mathcal{ C}={\,} \mathcal{ C}'$.\ \par 
\ \par 
{\,} {\,} \bf Proof. \rm Let $ x\in \mathcal{ S}_{0}$ and $ x_{\theta }={\,} 
(x_{k}e^{ik\theta })$, where $ \theta \in (-\pi ,\pi )$. We have 
$$
 \Big \Vert x_{\theta }\Big \Vert ^{2}_{\mathcal{ B}^{2}  (\mathcal{ C})}= 
 \displaystyle \sum _{j,k}c_{j,k}\Big \vert e^{ij\theta }x_{j}-e^{ik\theta 
}x_{k}\Big \vert ^{2}
$$
$$ 
={\,} \displaystyle \sum _{j,k}c'_{j,k}\Big \vert e^{ij\theta 
}x_{j}-e^{ik\theta }x_{k}\Big \vert ^{2}={\,} \Big \Vert x_{\theta }\Big \Vert 
^{2}_{\mathcal{ B}^{2}  (\mathcal{ C}')} , 
$$ 
and integrating in $ \theta \in (-\pi ,\pi )$, we get 
$$
 \displaystyle \sum _{j,k}c_{j,k}\Big (\Big \vert x_{j}\Big \vert ^{2}+{\,} 
\Big \vert x_{k}\Big \vert ^{2}\Big)={\,} \displaystyle \sum _{j,k}c'_{j,k}\Big (\Big \vert 
x_{j}\Big \vert ^{2}+ \Big \vert x_{k}\Big \vert ^{2}\Big ) .  
$$ 
Using the preceding equality, we deduce 
$$ 
\displaystyle \sum _{j,k}c_{j,k}{\rm Re} (e^{i(j-k)\theta }x_{j}\overline{x}_{k})={\,} 
\displaystyle \sum _{j,k}c'_{j,k}{\rm Re} (e^{i(j-k)\theta }x_{j}\overline{x}_{k}) .
$$ 
Letting $ \theta =0$, we obtain 
$$ 
\displaystyle \sum _{j,k}(c_{j,k}-c'_{j,k})x_{j}\overline{x}_{k}={\,} 
{\rm Re} \Big (\displaystyle \sum _{j,k}(c_{j,k}-c'_{j,k})x_{j}\overline{x}_{k}\Big )=  
 0 ,
 $$ 
 for every $ x\in \mathcal{ S}_{0}$, which implies $ c_{j,k}-c'_{j,k}=$ 
$ 0$ for all $ j$, $ k$. $ \qed $\ \par 
\ \par  \noindent 
\bf 2.7. Corollary. \it Let $ \mathcal{ C}=$ $ (c_{j,k})$ (satisfying the 
above conditions). Then, the shift operator $ S(x_{j})_{j\in \mathbb{Z}}={\,} 
(x_{j-1})_{j\in \mathbb{Z{\,} }}$is an isometry on $ \mathcal{ B}_{0}^{2}(c_{j,k})$ 
if and only if $ \mathcal{ C}=$ $ (c_{j,k})$ is a Toeplitz matrix $ c_{j,k}={\,} 
c_{\vert j-k\vert }$, where $ (c_{k})$ is a sequence satisfying $ c_{k}\geq 
0$, $ \forall k\geq 1;$ $ c_{0}={\,} 0$, $ 0<{\,} \sum _{k}c_{k}<{\,} 
\infty $.\ \par 
\ \par 
{\,} {\,} \rm Indeed, 
$$ \Big \Vert Sx\Big \Vert ^{2}_{\mathcal{ B}^{2}  (\mathcal{ C})}={\,} 
\displaystyle \sum _{j,k}c_{j,k}\Big \vert x_{j-1}-x_{k-1}\Big \vert ^{2}={\,} 
\displaystyle \sum _{j,k}c_{j+1,k+1}\Big \vert x_{j}-x_{k}\Big \vert ^{2},$$ 
and if $ \Big \Vert Sx\Big \Vert ^{2}_{\mathcal{ B}^{2}  (\mathcal{ C})}={\,} \Big \Vert 
x\Big \Vert ^{2}_{\mathcal{ B}^{2}  (\mathcal{ C})}$ for every $ x\in \mathcal{ S}_{0}$, 
we obtain by Lemma 2.6 $ c_{j+1,k+1}={\,} c_{j,k}$ for all $ j$, 
$ k$. Setting $ c_{j}={\,} c_{j,0}$ we get $ c_{j,k}=$ $ c_{\vert j-k\vert 
}$.\ \par 
{\,} {\,} Clearly, the converse is true as well.  \qed \ \par 
\ \par   \noindent

\section{L\'evy-Khinchin-Schoenberg weights}\label{sec3}

 The following lemma describes the spaces $ L^{2}(\mathbb{T}, \mu)$ for which the norm 
$$
 \Big \Vert \mathcal{ F}f\Big \Vert\!:={\,} \Big \Vert f\Big \Vert 
_{L^{2}  (\mu )}, \, \,  f\in \mathcal{ P} , 
$$
 (here $ \mathcal{ P}={\,} {\rm Lin}(z^{n}:{\,} n\in \mathbb{Z}), z \in \mathbb{T},$ is the set 
of all trigonometric polynomials) is a Besov-Dirichlet norm $ \Big \Vert \cdot 
\Big \Vert _{\mathcal{ B}^{2}  (\mathcal{ C})}$ on $ \mathcal{ S}_{0}$ (for a matrix 
$ \mathcal{ C}$). Note that $ \mathcal{ F}z^{n}={\,} e_{n}$, and hence 
$ \mathcal{ FP}={\,} \mathcal{ S}_{0}$. Speaking of the norms $ \Big \Vert \cdot 
\Big \Vert _{\mathcal{ B}^{2}  (\mathcal{ C})}$ we always suppose that the (Hermitian) 
matrix $ \mathcal{ C}$ satisfies conditions of Section \ref{sec2} (infinite $ C$-connected components $E_k$ for all k,  
and $ 0<{\,} \sum _{j}c_{j,k}<{\,} \infty $).\ \par 
\ \par  \noindent 
\bf  3.1. Lemma. \it Let $ \mu $ be a Borel measure on $ \mathbb{T}$. 
The following statements are equivalent.\ \par \ \par  \noindent 
\bf (1) $ x\longmapsto {\,} \Big \Vert \mathcal{ F}^{-1}x\Big \Vert _{L^{2}  (\mu 
)}$ \it is a Besov-Dirichlet norm on $ \mathcal{ S}_{0}$, $ \mathcal{ F}^{-1}x={\,} 
\displaystyle \sum _{k}x_{k}e^{ikt}$.\ \par \ \par  \noindent 
\bf (2) $ \mu ={\,} wm$\it , where  
$$ w(e^{it})={\,} 4 \sum _{k\geq 1} c_{k} \sin^{2}(kt/2) , $$
\it and $ c_{k}\geq 0$, $ 0<{\,} \sum_{k=1}^{\infty} c_{k}<{\,} \infty $; 
in this case, $ \mathcal{ F}L^{2}(\mathbb{T},w)={\,} \mathcal{ B}^{2}_{0}(c_{\vert 
j-k\vert }/2)$.\ \par 
\ \par 
{\,} {\,} \bf Proof. \rm For (2) $\Rightarrow $ (1), we simply 
observe that $ 4 \sin^{2}(kt/2)={\,} \vert 1-e^{ikt}\vert ^{2}$, and 
hence 
$$
 \Big \Vert \mathcal{ F}^{-1}x\Big \Vert ^{2}_{L^{2}  (w)}={\,} 
\displaystyle \sum _{k\geq 1}c_{k}\displaystyle \int _{\mathbb{T}}\Big \vert \mathcal{ F}^{-1}x(e^{it})\Big \vert 
^{2}\Big \vert 1-e^{ikt}\Big \vert ^{2}dm
$$
$$ 
={\,} \displaystyle \sum _{k\geq 1}c_{k}\displaystyle \int _{\mathbb{T}}\Big \vert 
\mathcal{ F}^{-1}x(e^{it})-{\,} e^{ikt}\mathcal{ F}^{-1}x(e^{it})\Big \vert ^{2}dm
$$
 $$ = \displaystyle \sum _{k\geq 1}c_{k}\displaystyle \int _{\mathbb{T}}\Big \vert 
\mathcal{ F}^{-1}x(e^{it})- e^{ikt}\mathcal{ F}^{-1}x(e^{it})\Big \vert ^{2}dm={\,} 
\displaystyle \sum _{k\geq 1}c_{k}\displaystyle \sum _{j\in \mathbb{Z}}\Big \vert 
x_{j}-x_{j-k}\Big \vert ^{2}
$$ 
 $$ ={\,} {\frac{\displaystyle 1}{\displaystyle 2}} \displaystyle \sum 
_{j,l\in \mathbb{Z}}c_{\vert j-l\vert }\Big \vert x_{j}-x_{l}\Big \vert ^{2}={\,} 
\Big \Vert x\Big \Vert^2 _{\mathcal{ B}^{2}  (\mathcal{ C})} , 
$$ 
where $ c_{j,l}=$ $ c_{\vert j-l\vert }/2$.\ \par 
{\,} {\,} For (1) $\Rightarrow $ (2), we apply Corollary 2.7: 
the shift operator $ S$ is an isometry on $ L^{2}(\mathbb{T},\mu )$, 
and hence on $ \mathcal{ F}L^{2}(\mathbb{T},\mu )$, and so if $ x\longmapsto \Big \Vert 
\mathcal{ F}^{-1}x\Big \Vert _{L^{2}  (\mu )}$ is a Besov-Dirichlet norm, 
$ \Big \Vert \mathcal{ F}^{-1}x\Big \Vert^2 _{L^{2}  (\mu )}={\,} \displaystyle \sum _{j,l\in 
\mathbb{Z}}c_{j,k}\Big \vert x_{j}-x_{k}\Big \vert ^{2}$, the matrix 
$ \mathcal{ C}=$ $ (c_{j,k})$ is a Toeplitz one, that is, there exists 
a sequence $ (c_{k})$ such that $ c_{j,k}={\,} c_{\vert j-k\vert }$. 
Now, the same computation as before but read in the opposite way shows 
that 
$$ 
\Big \Vert \mathcal{ F}^{-1}x\Big \Vert^2 _{L^{2}  (\mu )}=  
 \displaystyle \sum _{j,l\in \mathbb{Z}}c_{\vert j-k\vert }\Big \vert x_{j}-x_{k}\Big \vert 
^{2}={\,} \Big \Vert \mathcal{ F}^{-1}x\Big \Vert ^{2}_{L^{2}  (w)} , 
$$
for every polynomial $ p={\,} \mathcal{ F}^{-1}x\in \mathcal{ P}$, where 
$ w(e^{it})=$ $ \sum _{k\geq 1}c_{k}\Big \vert 1-e^{ikt}\Big \vert ^{2}$. 
So 
$$
 \displaystyle \int _{\mathbb{T}}\Big \vert p\Big \vert ^{2}d\mu 
={\,} \displaystyle \int _{\mathbb{T}}\Big \vert p\Big \vert ^{2}wdm \quad  
( \forall p\in \mathcal{ P}) ,
$$ 
which obviously implies $ \mu ={\,} wm$.  \qed 
\ \par 
\ \par  \noindent 
\bf 3.2. Comments. (1) \rm Weights $ w$ of the type 3.1(2) first appeared in 
\cite{Lev1934}, \cite{Khi1934} as characteristic exponents of stationary stochastic 
processes with independent increments and continuous time, nowadays often 
called L\'evy processes. A vast theory and numerous applications of these 
processes are known, including deep connections with potential 
theory. The weights themselves were characterized by I. Schoenberg \cite{Sch1938}  
(see also J. von Neumann and I. Schoenberg \cite{vNS1941}):\ \par 
\it - a non negative function $ w\in C(\mathbb{T})$, $ w(1)={\,} 0$, 
is of the form $ w(e^{it})=$ $ 4\sum _{k\geq 1}c_{k} \sin^{2}(kt/2)$, 
where $ c_{k}\geq 0$, $ \sum c_{k}<{\,} \infty $, if and only if 
$ w$ is ``conditionally negative definite" in the following sense: $ \sum 
_{j,k}w(z_{j}\overline{z}_{k})a_{j}\overline{a}_{k}\leq {\,} 0$ 
for every choice of points $ z_{j}\in \mathbb{T}$ and numbers $ a_{j}\in 
\mathbb{C}$ such that $ \sum _{j}a_{j}={\,} 0$, or equivalently, \par 
- if and only if $ e^{-\epsilon w}$ is positive definite for every 
$ \epsilon >0$.\ \par 
{\,} {\,} \rm Schoenberg and von Neumann obtained these characterizations 
as a step in their  solution of a metric geometry problem, in order to describe 
the so-called ``screw lines" on a Hilbert space. The same class of functions 
appeared in the Beurling-Deny potential theory, see \cite{BeD1958}, \cite{Den1970}.\ \par \ \par  \noindent 
\bf (2) \rm Several properties of weights of this class (we call it the 
L\'evy-Khinchin-Schoenberg class, {\rm LKS}) are known; for instance $ w\in {\rm LKS}{\,} 
\Rightarrow {\,} w^{\epsilon }\in {\rm LKS}$, and $1/(w^{\epsilon})$ is positive 
definite if $0<\epsilon \le 1$ (\cite{Sch1938}, \cite{vNS1941}, \cite{Kre1944}; see also \cite{Lan1972}, Sec. 
VI.3.13). It is also clear that  
$$ 
w(e^{it})= \sum _{k\geq 1}c_{k}\vert 1-e^{ikt}\vert ^{2}={\,} 0 
$$
if and only if $ e^{it}$ is a root of unity of order $ d={\,} {\rm GCD}\{k\colon  
c_{k}>0\}$, and so the zero set of $ w$ is always finite (if $ w\not\equiv 0$). 
A generic $ w\in {\rm LKS}$ has only one zero at $ e^{it}={\,} 1$, but 
then $ w(e^{idt})$ is again an {\rm LKS} weight with zeros at the $ d$-th roots 
of unity $ e^{ikt/d}$, $ k={\,} 0,1, \ldots, d-1$. The infinite 
$C$-connected components property always holds for
 $ \mathcal{ F}L^{2}(\mathbb{T}, w)$, $ w\in {\rm LKS}$ 
(see Remark after Lemma 2.1).\ \par 
\ \par 
{\,} {\,} We now derive first consequences of Lemma 3.1 and 
Section \ref{sec2}, in particular, a preliminary form of a description of the 
algebra $ {\rm Mult}(L^{2}(w))$ for $ w\in {\rm LKS}$. For this, we need 
the following simple lemma.\ \par 
\ \par  \noindent 
\bf 3.3. Lemma. \it Let $ \mu$ be a Borel measure 
on $ \mathbb{T}$, and let $ \mu ={\,} \mu _{s}+wm$ be 
its Lebesgue decomposition ($\mu_s$ is the singular part of $\mu$, and $ w\in L^{1}(\mathbb{T})$). 
The following statements are equivalent.\ \par \ \par  \noindent 
\bf (1) \rm $ \varphi _{n}:f\longmapsto \hat f(n)$ (defined on trigonometric 
polynomials) extends to a bounded functional on $ L^{2}(\mathbb{T},\mu )$.\ \par \ \par  \noindent 
\bf (2) \rm $ (z^{k})_{k\in \mathbb{Z}}$ is a minimal sequence in $ L^{2}(\mathbb{T},\mu 
)$.\ \par \ \par  \noindent 
\bf (3) \rm $ 1/w\in L^{1}(\mathbb{T})$.\ \par \ \par  \noindent 
\bf (4) \rm $ L^{2}(\mathbb{T},w)\subset {\,} L^{1}(\mathbb{T})$.\ \par 
\ \par 
{\,} {\,} \bf Proof. \rm It is clear that  (1) $\Leftrightarrow $ 
(2), (3) $\Rightarrow $ (4) (by Cauchy's inequality) and  (4) $\Rightarrow $ 
(1). Let us show that  (1) $\Rightarrow $ (3). Indeed, if $ f\longmapsto \hat 
f(0)$ is bounded, then there exists $ g\in L^{2}(\mathbb{T},\mu )$ such that 
$ \int _{\mathbb{T}}fdm={\,} \int _{\mathbb{T}}fgd\mu $ for every trigonometric 
polynomial $ f$. Hence, $ m={\,} g\mu ={\,} g(\mu _{s}+wm)$, 
and so $ g={\,} 0$ $ \mu _{s}$-a.e. and $ 1={\,} gw$ $ m$-a.e., 
which gives $ \int _{\mathbb{T}}{\frac{1}{w}} dm={\,} \int _{\mathbb{T}}g^{2}wdm\leq 
{\,} \int _{\mathbb{T}}g^{2}d\mu <{\,} \infty $. \qed \ \par 
\ \par \noindent 
\bf 3.4. Theorem. \it Let $ w(e^{it})= 4\sum _{k\geq 1}c_{k} \sin^{2}(kt/2)$ 
be an {\rm LKS} weight, $ c_{k}\geq 0$, $ 0<{\,} \sum _{k}c_{k}< \infty $. 
Then $(\lambda _{k})\in {\rm Mult} (L^{2}(w))$ if and only if  
$$
 \lambda \in l^{\infty}(\mathbb{Z}) \quad  {\rm and} \quad \sum _{k}\Big \vert 
\hat f(k)\Big \vert ^{2}\mu _{k}^{2}\leq  C^{2}\Vert f\Vert ^{2}_{L^{2}  (w)}, 
\quad  \forall f\in \mathcal{P}, 
$$ 
where $ \mu _{k}^{2}= \mu _{k}^{2}(\lambda ) := \displaystyle \sum _{j}c_{\vert 
k-j\vert }\Big \vert \lambda _{j}-\lambda _{k}\Big \vert ^{2}$, and 
the following alternative holds:\ \par 
\ \par \noindent 
\bf (1) \it either $ {\frac{1}{w}} \in L^{1}(\mathbb{T})$, and then $ \mathcal{S}_{0}\subset {\,} {\rm Mult}(L^{2}(w))$,\ \par 
\ \par \noindent 
\bf (2) \it or $ {\frac{1}{w}} \not\in L^{1}(\mathbb{T})$, and then ${\rm Mult}(L^{2}(w))$ 
consists of all sequences constant on every 
$ C$-connected component $ E_{k} = k + D \mathbb{Z}$, $ k= 0, \ldots , D-1$, of $\mathbb{Z}$ 
(see Remark after Lemma 2.1 for definitions); in particular, $ {\rm dim} \, {\rm Mult}(L^{2}(w))= 
D< \infty $.\ \par 
\ \par 
{\,} {\,} \bf Proof. \rm (1) Lemma 3.1 gives $ \mathcal{ F}L^{2}(\mathbb{T},w)={\,} 
\mathcal{ B}^{2}_{0}(c_{\vert j-k\vert }/2)$, and by Lemma 3.3 $ \varphi _{n}$ 
are bounded on $ \mathcal{ B}^{2}_{0}(c_{\vert j-k\vert }/2)$, so that 
Lemmas 2.3 and 2.2 are applicable and yield the statement.\ \par \ \par
\noindent 
(2) The references to the same lemmas show that  
all functionals $ \varphi _{n}$ are unbounded, and consequently 
the space of multipliers is finitely dimensional (see Remark to Lemma 
2.3 above for details).  \qed \ \par 
\ \par \noindent 
\bf Remark. \rm The inequality 
$$\sum _{k \in \mathbb{Z}}\Big \vert 
\hat f(k)\Big \vert ^{2}\mu_k^{2}\leq C^{2}\Vert f\Vert ^{2}_{L^{2}  (w)} , \quad 
 \forall f\in \mathcal{ P} , $$
is equivalent to the embedding 
$\mathcal{ F}L^{2}(\mathbb{T},w)\subset l^{2}(\nu)$, where $\nu= (\mu_k^{2})_{k \in \mathbb{Z}}$. Here and 
below we use the notation $l^2(\nu)=l^{2}(\mathbb{Z},\nu)$ 
for the weighted $l^2$ space with norm 
$$
||x||_{l^{2}(\nu)} = \Big ( \sum_{k \in \mathbb{Z}}\Big \vert x_k\Big \vert^{2}\nu_{k}\Big)^{\frac 1 2},   
$$
where $\nu= (\nu_k)_{k \in \mathbb{Z}}$ ($\nu_k>0$). For a capacitary characterization of this embedding property  
see Section 4 below.\ \par 
\ \par  \noindent 
\bf 3.5. On the condition $ 1/w\in L^{1}(\mathbb{T})$ for {\rm LKS} weights. 
\rm The integrability condition $ 1/w\in L^{1}(\mathbb{T})$ plays a  
key role in the description of $ {\rm Mult}(L^{2}(\mathbb{T},w))$ in Theorem 
3.4. We discuss it below using the following simple observation.\ \par 
\ \par  \noindent 
\bf 3.6. Lemma. \it Let $ w(e^{it})=$ $ 4\sum _{k\geq 1}c_{k} \sin^{2}(kt/2)$ 
be an {\rm LKS} weight, $ c_{k}\geq 0$, $ 0<$ $ \sum _{k}c_{k}<$ $ \infty $. 
Then, 
$$ 
{\frac{4t^{2}}{\pi ^{2}}} \sum _{k\leq \pi /t}c_{k}k^{2}\leq 
w(e^{it})\leq t^{2}\sum _{k=1}^{N}c_{k}k^{2}+{\,} 4\sum _{k>N}c_{k} ,$$
\it for   $t\in (0,\pi )$ and for every $ N\geq 0$. In particular, 
$$ 
c \displaystyle \sum _{n\geq 1}{\frac{\displaystyle 1}{\displaystyle 
\displaystyle \sum _{k=1}^{n}c_{k}k^{2}+n^{2}\displaystyle \sum _{k>n}c_{k}}} 
\leq {\,} \displaystyle \int _{\mathbb{T}}{\frac{\displaystyle dm}{\displaystyle 
w}} \leq {\,} C\displaystyle \sum _{n\geq 1}{\frac{\displaystyle 1}{\displaystyle 
\displaystyle \sum _{k=1}^{n}c_{k}k^{2}}} , 
$$ 
\it with appropriate absolute constants $ 0<c\leq C<\infty $.\ \par 
\ \par 
{\,} {\,} \bf Proof. \rm For $ 0\leq {\,} kt/2\leq {\,} \pi 
/2$, one has $ (kt/\pi )^{2}\leq {\,} \sin^{2}(kt/2)$, and $ \sin^{2}(kt/2)\leq 
$ $ (kt/2)^{2}$ for every $ k$.\ \par 
{\,} {\,} For the integral $ \int _{-\pi }^{\pi }{\frac{dt}{w(e^{it})}} 
$, we first integrate  around $t=0$: $ \int _{-\pi /N}^{\pi /N}={\,} 2\sum 
_{n\geq N}\int _{\pi /n+1}^{\pi /n}$ and 
$$ 
c \, {\frac{\displaystyle \pi (1/n-1/(n+1))}{\displaystyle (1/n^{2})\displaystyle \sum 
_{k=1}^{n}c_{k}k^{2}+\displaystyle \sum _{k>n}c_{k}}} \leq \displaystyle \int _{\pi 
/n+1}^{\pi /n}{\frac{\displaystyle dm}{\displaystyle w}}
 $$
$$ 
\leq {\,} C \, {\frac{\displaystyle \pi (1/n-1/(n+1))}{\displaystyle 
(1/n^{2})\displaystyle \sum _{k=1}^{n}c_{k}k^{2}}} ,
$$
 which gives the  estimate claimed above if the only zero of $ w(e^{it})$ is at $t=0$. 
 If there are other zeros of $ 
w(e^{it})$ then, using comments 3.2(2), we can write $ w(e^{it})={\,} w_{1}(e^{idt})$, 
where $ w_{1}$ a {\rm LKS} weight with the only zero at $ t=0$, and the inequalities 
follow from $ \displaystyle \int _{\mathbb{T}}{\frac{\displaystyle dm}{\displaystyle 
w}} ={\,} \displaystyle \int _{\mathbb{T}}{\frac{\displaystyle dm}{\displaystyle 
w_{1}}} $.  \qed \ \par 
\ \par  \noindent 
\bf 3.7. Examples. \it As before, let $ w(e^{it}) = 4\sum _{k\geq 1}c_{k}\sin^{2}(kt/2)$.\ \par 
\ \par \noindent 
\bf (1) \rm Let $ c_{1}={\,} 1$ , $ c_{k}={\,} 0$  for 
$ k>1$; then $w =4 \sin^{2}(t/2)$. It follows that 
 $ 1/w\not\in L^{1}(\mathbb{T})$ and $ \mathcal{ F}L^{2}(\mathbb{T}, 4 \sin^2(t/2)dt
)=$ $ \mathcal{ B}^{2}_{0}(c_{\vert j-k\vert }/2)$, which is 
the completion of $ \mathcal{ S}_{0}$ in the norm 
$$ 
\Vert p\Vert _{L^{2}  (w)}={\,} \Vert \mathcal{ F}p\Vert 
_{\mathcal{ B}^{2}}={\,} \Big( \sum _{k\in \mathbb{Z}}\vert x_{k}-x_{k-1}\vert 
^{2}\Big)^{1/2} ,$$
for every polynomial $ p\in \mathcal{ P}$. Note that the completion $ \mathcal{ B}^{2}_{0}(c_{\vert 
j-k\vert }/2)$ is not a sequence space, but it can naturally be identified 
with $ L^{2}(\mathbb{T},4 \sin^2(t/2)dt)$. The functionals $ \varphi _{n}$ 
are not continuous, and hence $ {\rm Mult}(L^{2}(w))=
 \Big\{{\rm const}\Big\}$. 
The same conclusion is still true for any finitely supported sequence 
$ (c_{k})_{k\geq 1}$, or for sequences ``rapidly" tending to zero considered below.\ \par 
\ \par \noindent 
\bf (2) \it Power-like kernels $ c_{k}\approx {\frac{1}{k^{1+\alpha }}} 
$ and the spaces $ L^{2}(\mathbb{T},\vert 1-e^{it}\vert ^{\alpha })$, 
$ 0<{\,} \alpha <{\,} 2$. \rm We use the notation  $ c_{k}\approx b_{k}$ in the following sense: 
$$ 
c_{k}\approx b_{k} \Longleftrightarrow  ab_{k}\leq {\,} c_{k}\leq 
{\,} Ab_{k}\, \,  \text{for large} \, \,   k \, \, ( k\geq K) \, \, \text{and} \, \,   0<a\leq A<\infty . 
$$
Then, with appropriate constants $ C>0$ (which may be different in different 
entries) we have by 3.5 (for $ 0<t<\pi $), 
$$ 
w(e^{it})\leq {\,} Ct^{2}\sum _{k\leq \pi /t}k^{1-\alpha 
}+  C\sum _{k>\pi /t}1/k^{1+\alpha }\leq {\,} Ct^{\alpha }+{\,} Ct^{\alpha 
}={\,} Ct^{\alpha } ,
$$
 $$ 
 w(e^{it})\geq {\,} {\frac{t^{2}}{\pi ^{2}}} \sum _{k\leq 
\pi /t}c_{k}k^{2}\geq {\,} ct^{2}\sum _{k\leq \pi /t}k^{1-\alpha }\geq 
{\,} ct^{\alpha } ,
$$
so that $ w(e^{it})\approx {\,} \vert t\vert ^{\alpha }$ as 
$ t\longrightarrow 0$.\ \par 
{\,} {\,} \textit{Conclusion:}  For $ 0<\alpha <1$, we have $ 1/w\in 
L^{1}(\mathbb{T})$, $ \mathcal{ F}L^{2}(\mathbb{T},\vert 1-e^{it}\vert ^{\alpha })=$ 
$ \mathcal{ B}^{2}_{0}({\frac{\displaystyle 1}{\displaystyle  \vert j-k \vert 
^{1+\alpha }}} )$ (with the equivalence of norms), $ \mathcal{ S}_{0}\subset 
{\,} {\rm Mult}(L^{2}(\mathbb{T},\vert 1-e^{it}\vert ^{\alpha }))$ and 
$$ 
(\lambda_{k})\in {\rm Mult}(L^{2}(\mathbb{T},w)) \Leftrightarrow
 \lambda \in l^{\infty }(\mathbb{Z}) \, \,   \text{and} \, \,  \mathcal{ F}L^{2}(\mathbb{T},\vert 
1-e^{it}\vert ^{\alpha })\subset {\,} l^{2}(\nu) , 
$$  
where $ \nu _{k}=$ $ \mu _{k}(\lambda )^{2}:= \displaystyle \sum _{j}{\frac{\displaystyle 
\vert \lambda_{j}-\lambda_{k}\vert ^{2}}{\displaystyle (\vert j-k\vert 
+1)^{1+\alpha }}} .$ (See Section \ref{sec4} for a characterization of the last  
embedding.) \ \par 
{\,} {\,} {\,} \rm For $ 1\leq \alpha <2$, we have $ 1/w\not\in L^{1}(\mathbb{T})$,  
and hence $ {\rm Mult}(L^{2}(w))= \Big\{{\rm const}\Big\}$.  
Clearly, for larger $ \alpha $ ($ \alpha \geq 2$) 
the preceding equality holds as well.  \qed \ \par 
\ \par 
{\,} {\,} The following elementary lemma explains the condition 
$ 1/w\in L^{1}(\mathbb{T})$ for an {\rm LKS} weight $ w(e^{it})=$ $ 4\sum _{k\geq 
1}c_{k} \sin^{2}(kt/2)$ in the ``critical band" between $ c_{k}={\,} 1/k$ 
(decreasing too slowly  since $ \sum_{k}c_{k}={\,} \infty $) and $ c_{k}={\,} 
1/k^{2}$ (decreasing too fast since 
$$ 
{\frac{2}{\pi ^{2}}} \sum _{k\geq 1}{\frac{1}{k^{2}}} 
\sin^{2}(kt/2)={\,} {\frac{t}{2\pi }} (1-{\frac{t}{2\pi }} ), \quad  0<t<2\pi , 
$$
and so $ 1/w\not\in L^{1}(\mathbb{T})$). By the way, the last observation 
shows that, for this integration question, without loss of generality  
we can assume that $ \sum _{k}c_{k}k^{2}={\,} \infty $ (if not, then 
surely $ 1/w\in L^{1}(\mathbb{T})$).\ \par 
\ \par  \noindent 
{\bf 3.8. Lemma.} {\it Let $ x\longmapsto c(x)$ ($ x\in [0,\infty )$) 
be a positive piecewise differentiable function such that $ c_{k}= c(k)$, and }
$$
 x\longmapsto x^{\gamma } c(x) \, \,  \textit{eventually decreases  
for some} \,\,  \gamma, \,  1<\gamma <3 .
$$
{\it Then} 
$$
 w(e^{it}) \approx {\,} W(t)=:{\,} t^{2}\displaystyle \int _{0}^{\pi 
/t}c(x)x^{2}dx \quad  {\rm as} \, \,    t\longrightarrow 0 , 
$$
{\it and consequently}  
$$ 
1/w\in L^{1}(\mathbb{T}) \Leftrightarrow  \displaystyle \sum _{n\geq 
1}{\frac{\displaystyle 1}{\displaystyle \displaystyle \sum _{k=1}^{n}c_{k}k^{2}}} 
<{\,} \infty .
$$
\ \par 
{\,} {\,} {\bf Proof.} First note that $ \lim_{x\longrightarrow \infty 
}xc(x)={\,} 0$ and $ \int _{0}^{\infty }c(x)dx<{\,} \infty $\it . 
\rm By Lemma 3.6, it suffices to prove that there is a constant $ C>0$ 
such that 
$$ 
Ct^{2}\sum _{k\leq \pi /t}c_{k}k^{2}\geq {\,} \sum 
_{k>\pi /t}c_{k}, \quad \text{or} \quad Ct^{2}\int _{0}^{1/t}c(x)x^{2}dx\geq  \int _{1/t}c(x)dx .
$$
\ \par 
{\,} {\,} \rm By the hypothesis $ \gamma x^{\gamma -1}c(x)+x^{\gamma 
}c'(x)\leq {\,} 0$ (for $ x\geq a>0$), and hence $ \gamma c(x)+xc'(x)\leq 
$ $ 0$. Integrating over $ [y,b]$ and letting $ b\longrightarrow \infty $, 
we obtain 
$$ 
(\gamma -1)\int _{y}^{\infty }c(x)dx-{\,} yc(y)\leq 0 \quad  (\text{for} \,  y\geq a) .
$$ 
Multiplying by $ y$ and integrating over $ [a,s]$ we get 
$$ 
0\geq {\,} (\gamma -1)\int _{a}^{s}\Big (\int _{y}^{\infty }c(x)dx \Big) d(y^{2}/2)-  
\int _{a}^{s}y^{2}c(y)dy
$$ 
$$
 ={\,} (\gamma -1)2^{-1}\Big (s^{2}\int _{s}^{\infty }c(x)dx-{\,} 
a^{2}\int _{a}^{\infty }c(x)dx+{\,} \int _{a}^{s}y^{2}c(y)dy \Big) - 
 \int _{a}^{s}y^{2} c(y) dy , 
$$ 
$$
{\frac{3-\gamma }{2}} \int _{a}^{s}y^{2}c(y)dy \geq {\,} 
(\gamma -1)2^{-1}s^{2}\int _{s}^{\infty }c(x)dx- {\rm const} , 
$$ 
which is equivalent to the inequality claimed above (with any constant 
$ C>{\,} {\frac{3-\gamma }{\gamma -1}} $).  \qed \ \par 
\ \par 
\noindent
\bf Remark. \rm Yet another combination of ``regularity conditions" on the 
behaviour of $ c_{k}$ as $ k\to \infty $ leads to the following 
criterion:\ \par 
\ \par 
 \it Assume $ c_{k}= c(\vert k\vert )$ ($ k\in \mathbb{Z}\setminus\{0\}$) where $ c: [1,\infty )\longrightarrow \mathbb{R}_{+}$ is a function 
satisfying $ c(t)t^{2}\uparrow \infty $ (eventually) and $ c(xy)\leq 
A c(x) c(y)$ ($ x,y\geq 1$) (in particular, all $ c(t)= t^{-\gamma }$, $ \gamma \leq 
2$ satisfy these conditions). Then $ {\frac{1}{w}} \in L^{1}({\Bbb T})$ if and only if 
$ \sum _{k\geq 1}{\frac{1}{k^{3}c_{k}}} < \infty $ (or $ \int _{1}^{\infty 
}{\frac{dx}{x^{3}c(x)}} < \infty $).

\rm Indeed, as in Lemma 3.8, we compare functions 
$ B(x)= \int _{0}^{x}c(t)t^{2}dt$ and $ C(x)= x^{2}\int _{x}^{\infty 
}c(t)dt$. We have 
$$ a\cdot c({x}/{2}){\frac{x^{3}}{8}} \leq  \int _{x/2}^{x}c(t)t^{2}dt\leq 
 B(x)\leq  A\cdot c(x)x^{2}\int _{1}^{x}dt= Ac(x)x^{3} ,$$
where $ a>0, A>0$ are constants, and
$$ C(x)= x^{2}\int _{x}^{\infty }c(t)dt= x^{3}\int 
_{1}^{\infty }c(xy)dy\leq  a\cdot x^{3}c(x)\int _{1}^{\infty }c(y)dy=  
A\cdot x^{3}c(x) .$$ 
So, by Lemma 3.6, if $ \int _{1}^{\infty }{\frac{dx}{x^{3}c(x)}} 
<$ $ \infty $, we obtain $ \int _{1}^{\infty }{\frac{dx}{B(x)}} < \infty $, and hence $ {\frac{1}{w}} \in L^{1}({\Bbb T})$. Conversely,  
 if $ {\frac{1}{w}} \in L^{1}({\Bbb T})$, then $ \int _{1}^{\infty }{\frac{dx}{B+C}} 
< \infty $ and $ {\frac{1}{B+C}} \geq  {\frac{1}{Ax^{3}c(x)+Ax^{3}c(x)}} 
$, whence $ \int _{1}^{\infty }{\frac{dx}{x^{3}c(x)}} < \infty $. \qed 

\ \par 
\ \par \noindent 
\bf 3.9. On the Muckenhoupt condition $ w\in (A_{2})$ for {\rm LKS} weights. \rm 
\ \par \ \par\noindent 
\bf (1) \rm Condition $ w\in (A_{2})$ is not so transparent as $ 1/w\in L^{1}(\mathbb{T})$ 
even for {\rm LKS} weights. Recall that by definition   
$$ 
w\in (A_{2}) \Leftrightarrow  \Big ({\frac{\displaystyle 1}{\displaystyle 
\vert I\vert }} \displaystyle \int _{I}wdm\Big )\Big ({\frac{\displaystyle 1}{\displaystyle 
\vert I\vert }} \displaystyle \int _{I}{\frac{\displaystyle 1}{\displaystyle w}} 
dm\Big ) \le C , \quad  \forall I\subset \mathbb{T} ,
$$ 
where $ I$ is an arc (interval), and $C$ is a constant which 
does not depend on $I$. Using P. Jones's $ (A_{p})$-factorization 
theorem ($ w\in (A_{p}){\,} \Leftrightarrow {\,} w={\,} v_{0}v_{1}^{1-p}$, 
where $ v_{0},v_{1}\in (A_{1})$, see \cite{Duo2001}, p.150), we get $ w\in 
(A_{2})$ $ \Leftrightarrow $ $ w=$ $ v_{0}/v_{1}$, where $ v_{0},v_{1}\in 
(A_{1})$. For a weight $ w\in {\rm LKS}$, a good sufficient condition for 
$ w\in (A_{2})$ is simply $ v={\,} 1/w\in (A_{1})$, which means 
that there exists $C>0$ such that 
$$ 
C v(x)\geq {\,} {\frac{1}{\vert I \vert }} \int _{I}v \, dm, \quad \text{for a.e.} {\,\, } x\in I ,
$$ 
for every arc (interval) $ I \subset \mathbb{T}$.\ \par 
{\,} {\,} Identifying $ \mathbb{T}={\,} (-\pi ,\pi)$, it is 
easy to see that for $ 1/w\in (A_{1})$ it suffices to check 
$$
{\frac{\displaystyle C}{\displaystyle w(y)}} \geq {\,} 
{\frac{\displaystyle 1}{\displaystyle y-x}} \displaystyle \int _{x}^{y}{\frac{\displaystyle 
dt}{\displaystyle w(t)}}  \quad \text{for all} \quad 0<x<y<\pi .
$$ 
If the generating function $ c(k)={\,} c_{k}$ satisfies the 
condition of Lemma 3.8, one can replace $ w$ by $ W$ from this Lemma. 
The needed inequality $ {\frac{\displaystyle C}{\displaystyle W(y)}} 
\geq $ $ {\frac{\displaystyle 1}{\displaystyle y-x}} \displaystyle \int _{x}^{y}{\frac{\displaystyle 
dt}{\displaystyle W(t)}} $ (for all $ 0<x<y<\pi $), follows from the following 
H\"older type condition (which defines ``power-like" behaviour of $W$):  
$$
 {\frac{W(y)}{W(t)}} \leq {\,} C(y/t)^{\gamma }, \quad 
\text{where} \quad  0<t<y<\pi \, \, \text{and} \, \,  0<\gamma <1 . 
$$ 
Indeed, the preceding condition implies 
$$ 
{\frac{\displaystyle 1}{\displaystyle y-x}} \displaystyle \int _{x}^{y}{\frac{\displaystyle 
W(y)}{\displaystyle W(t)}} dt\leq {\,} {\frac{\displaystyle C}{\displaystyle y-x}} 
\displaystyle \int _{x}^{y}(y/t)^{\gamma }dt$$
$$
 ={\frac{\displaystyle Cy^{\gamma }}{\displaystyle (y-x)(1-\gamma 
)}} (y^{1-\gamma }-x^{1-\gamma })={\,} {\frac{\displaystyle C}{\displaystyle 1-\gamma 
}} \, {\frac{\displaystyle 1-(x/y)^{1-\gamma }}{\displaystyle 1-(x/y)}} 
\leq {\frac{\displaystyle C}{\displaystyle 1-\gamma }} .
$$
\ \par  \noindent 
\bf (2) \it {\rm LKS} weights satisfying $ 1/w\in L^{1}(\mathbb{T})$ but not 
$ w\in (A_{2})$. \rm In the notation of Lemma 3.8, let $ c(x)={\,} x^{-2}(\log(ex))^{\beta 
}$, $ \beta >1$. Then  
$$ 
w(e^{it})\approx  W(t)=: \,  t^{2}\displaystyle \int _{0}^{\pi 
/t}c(x)x^{2}dx ={\,} t^{2}\displaystyle \int _{0}^{\pi /t}(\log(ex))^{\beta 
}dx\approx {\,} t(\log(e/t))^{\beta } , 
$$
which gives $ \int _{0}^{1}{\frac{dt}{w(e^{it})}} <{\,} \infty $. 
On the other hand, 
$$ 
\int _{0}^{1}{\frac{dt}{w(e^{it})^{1+\epsilon }}} \geq 
{\,} \text{const} \int _{0}^{1}{\frac{dt}{(t(\log(e/t))^{\beta })^{1+\epsilon }}} 
={\,} \infty, 
$$
for every $ \epsilon >0$, and so $ w\not\in (A_{2})$.  \qed 

\ \par 
\ \par \noindent 
\bf 3.10. Typical asymptotic behaviour of {\rm LKS} weights at $0$. \rm Consider 
an {\rm LKS} weight $ w(e^{it})=4 \sum_{k\geq 1}c_{k} \sin^{2}{\frac{kt}{2}} 
$ as a function of the argument $ t$, $ -\pi \le t \le \pi $. The following claim 
shows that an arbitrary ``mildly regular" (convexity-like) behaviour 
is permitted for an {\rm LKS} weight as $ t\longrightarrow 0$. Since the coefficients 
$ k\longmapsto c_{k}$ are nearly monotone, the resulting function is 
always power-like: $ \vert t\vert ^{2}\preceq  w(e^{it}) \preceq 
\vert t\vert ^{\epsilon }$ as $ t\longrightarrow 0$ for some $ \epsilon >0$. 
Here $\phi(t)\preceq \psi(t)$  means that $\phi(t) = O(\psi(t))$ as $t\longrightarrow 0$. 
\ \par 
\ \par\noindent  
\bf Claim. \it Let $ u:(0,\infty )\longrightarrow (0,\infty )$ be an 
(eventually) increasing piecewise differentiable function such that for some $ -1< \alpha < 
1$ the function $ s\longmapsto s^{\alpha }u'(s)$ (eventually) decreases. 
Then there are $ c_{k}>0$, $ \sum _{k}c_{k}< \infty $, such that\ \par 
\ \par 
  \centerline{ $ w(e^{it})\approx  t^{2}u(1/\vert t\vert )$ 
\rm as $ t\longrightarrow 0$.}
\ \par \noindent 
\it In particular, $ u(s)= s^{\beta }$, $ 0<\beta <2$, gives 
$ w(e^{it})\approx $ $ t^{2-\beta }$.\ \par 
\ \par 
 \rm Indeed, define the function $ c(\cdot )$ by $ c(\pi y)(\pi 
y)^{2}= u'(y)$, $ y>0$. It follows that, for $ 1< \gamma =:  \alpha +2< 3$, the function 
 $ c(\pi y)(\pi y)^{\gamma }= (\pi 
y)^{\gamma -2}u'(y)$ 
eventually decreases (and hence, $ \int _{0}^{\infty }c(x)dx<\infty $). Then by Lemma 3.8 (for $ t>0$),\ \par 
\ \par 
  \centerline{ $ w(e^{it})\approx W(t)= t^{2}\int _{0}^{\pi 
/t}c(x)x^{2}dx= t^{2}\int _{0}^{\pi /t}u'(x/\pi )dx=$}
\ \par 
  \centerline{ $ = \pi t^{2}(u(1/t)-u(0))\approx  t^{2}u(1/t)$ 
\rm as $ t\longrightarrow 0$.}  \qed 
\ \par  \ \par\noindent 
\bf 3.11. Remarks on trivial multipliers  
for non-LKS weights. \rm It is easy to see that ${\rm Mult}(L^{2}(w))= 
 \Big\{{\rm const}\Big\}$ implies $ 1/w\not\in L^{1}(\mathbb{T})$ for every $ w\in L^{1}(\mathbb{T})$ 
(not only for {\rm LKS} weights). On the other hand, the converse is \textit{generally not true for non-{LKS}
weights} $ w\in L^{1}(\mathbb{T})$: indeed, let\ \par 
\ \par 
  \centerline{ $ 1/w= \sum_{k\in \mathbb{Z}}a_{k}\vert z-\alpha ^{k}\vert 
^{-1}$,}
\ \par \noindent 
\rm where $ \sum _{k}a_{k}< \infty $ ($ a_{k}>0$), $ 0< \inf_{k}(a_{k}/a_{k+1})\leq 
\sup_{k}(a_{k}/a_{k+1})< \infty $, and let $ \alpha \in \mathbb{T}$, 
$ \alpha ^{k}\not= 1$ ($ \forall k\in \mathbb{Z}$). It is clear that 
the series converges a.e. (it is in $ L^{p}(\mathbb{T})$ for every $ 0<p<1$) 
and $ w\in L^{\infty }(\mathbb{T})$, but $ 1/w\not\in L^{1}(\mathbb{T})$. Then it 
can be proved by the same reasoning as in \cite{Nik2009} that the corresponding rotations 
defined by $ T_{\alpha ^{k}}z^{n}=(\alpha ^{k})^{n}z^{n}$ ($ \forall n\in \mathbb{Z}$) 
are (non-trivial) multipliers of  $ L^{2}(\mathbb{T},w)$, and hence 
${\rm dim} \, {\rm Mult}(L^{2}(w))=\infty$.   
\ \par  \ \par

\noindent 
\bf 3.12. Remarks on duality of multipliers   
for general weights. \rm   Suppose $w^{\pm 1} \in L^1(\mathbb{T})$. Then 
\ \par  \ \par\noindent 
\bf (1) \rm $\displaystyle{\lambda \in {\rm Mult}(L^{2}(w))\Leftrightarrow \overline{\lambda} \in {\rm Mult}(L^{2}(1/w))}$,  
where $\overline{\lambda} = (\overline{\lambda}_{j})$. \ \par \ \par Indeed, using the duality 
$\langle f, g \rangle =  \int _{\mathbb{T}}f \overline{g}dm$ 
yields  $ (L^{2}(w))^{*}= L^{2}(1/w)$, and  $ T_{\lambda }^{*}= 
(\overline{\lambda }_{j})_{j\in \mathbb{Z}}$. \ \par  \ \par\noindent 
\bf (2) \rm In general, $\lambda \in {\rm Mult}(L^{2}(w))$ $\not\Rightarrow$ $\overline{\lambda} \in {\rm Mult}(L^{2}(w))$, and consequently 
$$
{\rm Mult}(L^{2}(w)) \not= {\rm Mult}(L^{2}(1/w)), 
$$
even for weights with the {\rm SLP}  (see Example 5.8 below). However, for {\rm LKS} weights, obviously 
${\rm Mult}(L^{2}(w)) = {\rm Mult}(L^{2}(1/w))$.  
 \ \par  \ \par\noindent 
\bf (3) \rm For every weight $W \in L^1(\mathbb{T})$, 
$$
 \overline{\lambda} \in {\rm Mult}(L^{2}(W))\Leftrightarrow\widetilde{\lambda} \in {\rm Mult}(L^{2}
 (W))\Leftrightarrow\lambda \in {\rm Mult}(L^{2}(\widetilde{W})) ,
 $$
 where $\widetilde{\lambda} =(\lambda_{-j})$, and $\widetilde{W}(z) = W(\overline{z})$ ($z \in \mathbb{T}$). In particular, 
 if $w^{\pm 1} \in L^1(\mathbb{T})$, then 
 $$
 {\rm Mult}(L^{2}(w)) = {\rm Mult}(L^{2}(1/\widetilde{w})) . 
 $$
  Indeed, 
  \begin{equation*}
\begin{aligned}
\Vert T_{\overline{\lambda}}f\Vert ^{2}_{L^2(W)} & = \int_{\mathbb{T}}\Big \vert \sum \overline{\lambda}_{j}\hat f(j) z^{j}\Big \vert ^{2}W(z)dm =  
\int_{\mathbb{T}}\Big \vert \sum \lambda_{j}\overline{\hat f(j)} z^{-j}\Big \vert 
^{2}W(z)dm\\ & = \int_{\mathbb{T}} \Big \vert \sum \lambda_{-j}\overline{\hat f(-j)} z^{j}\Big \vert 
^{2}W(z)dm=\Vert T_{\widetilde{\lambda}}\overline{f}\Vert ^{2}_{L^2(W)}, 
\end{aligned}
\end{equation*}
where $f$, and consequently $ \overline{f} = \sum \overline{\hat f(-j)} z^j$, is  an arbitrary trigonometric polynomial, and 
 $ \Vert \overline{f}\Vert_{L^2(W)}=  \Vert \widetilde{f}\Vert_{L^2(\widetilde{W})}=
 \Vert f\Vert_{L^2(W)}$.

\section{Embedding of Besov-Dirichlet sequence spaces into 
weighted \\ $l^2$ spaces}\label{sec4}

In this section we continue to consider  {\rm LKS} weights  $w$ such that 
$1/w\in L^1(\mathbb{T})$, where 
$$ w(e^{it})={\,} 4 \sum _{k\geq 1}c_{k} \sin^{2}(kt/2) , $$
and $ c_{k}\geq 0$, $ 0<{\,} \sum_{k=1}^{\infty} c_{k}<{\,} \infty $.  
To shorten the notation, we will denote by $D$ the corresponding Besov-Dirichlet space:
  $$D= \mathcal{ B}^{2}_{0}(c_{\vert 
j-k\vert }/2)= \mathcal{ F}L^{2}(\mathbb{T},w).$$ Note that $w \in C(\mathbb{T})$, and consequently 
$D\supset l^2=l^2(\mathbb{Z})$. 

 Here we present a characterization  of the 
multiplier algebra ${\rm Mult}(D)={\rm Mult}(L^{2}(w))$  for general {\rm LKS}  weights $ w$ such that 
$1/w\in L^1(\mathbb{T})$ in terms of \textit{capacities} associated with $D$, as well as some \textit{non-capacitary} characterizations for $w$  
satisfying additional ``regularity" conditions (which includes standard 
power-like weights $ \vert e^{it}-e^{i\theta}\vert ^{\alpha }$, $ 0<\alpha <1$). 
To this end we will need some elements of potential theory, whose adaptation 
to our situation is presented below for the reader's convenience.

As was shown in Sec. \ref{sec3},   $\lambda =(\lambda_k)_{k \in \mathbb{Z}}\in {\rm Mult}(D)$ 
if and only if $\lambda \in l^\infty$, and  the embedding $D  \subset l^2(\nu)$ holds; the last 
condition is equivalent to the inequality
$$
||x||_{l^2(\nu)} \le C \,  ||x||_{D}, \quad \forall x \in \mathcal{S}_0, 
$$
where $\nu = \{\nu_j\}_{j \in \mathbb{Z}}$ is a nonnegative weight  given by  
$$ \nu _{j}=  \mu_{j}(\lambda )^{2}= \displaystyle \sum _{k}{\displaystyle c_{|j-k|}
\vert \lambda_{j}-\lambda_{k}\vert ^{2}} .$$ In this section, we consider general embeddings of 
a given Besov-Dirichlet space $ D\subset  l^{2}(\nu )= l^{2}(\mathbb{Z},\nu 
)$, where $ \nu = (\nu _{j})_{j\in \mathbb{Z}}$ is an arbitrary 
nonnegative weight on $ \mathbb{Z}$, not necessarily related to a  
mulltiplier $ \lambda \in {\rm Mult}(D)$. Later on,  we will be applying these results 
to multiplier related weights $\nu$  with $ \nu _{j}= \mu_{j}(\lambda )^{2}$.
\ \par 
\ \par \noindent 
\bf 4.1. Green's kernel. \rm  Let $ w \in {\rm LKS}$ and $1/w \in L^1(\mathbb{T})$. 
As was mentioned in Sec. \ref{sec3}, this yields that $w^{\epsilon} \in {\rm LKS}$, and $1/(w^\epsilon)$ is positive definite, for all $0<\epsilon\le 1$. In particular, both  
$1/w$ and $1/({w^{1/2}})$ are positive definite. Consider the discrete Green kernel $\Big (g_{m-j}\Big)_{j, m \in  \mathbb{Z}}$,
where $g = \mathcal{F} (1/w) \ge 0$, so that  
$$
 1/w(e^{it}) = 
\sum_{j \in \mathbb{Z}}  g_j e^{i j t}.   
$$  
The Green potential 
$G x$ is defined by $g*x$,  
$$
(G x)_{m} = \sum_{j \in \mathbb{Z}} g_{m-j} \,  x_j, \quad m \in \mathbb{Z}.   
$$
We will also need   the corresponding potential operator $K x$ defined by $\kappa*x$, 
$$
(Kx)_m = \sum_{j \in \mathbb{Z}} \kappa_{m-j} \,  x_j, \quad m \in \mathbb{Z},    
$$
where $\kappa = \mathcal{F}  (1/{w^{1/2}}) \ge 0$, so that $\kappa*\kappa=g$. 

It follows that both $G$ and $K$ have  nonnegative symmetric kernels. Moreover, 
by Parceval's theorem $K$ is an isometry from $l^2$ 
onto $D$. Consequently, for any nonnegative weight $\nu = \{\nu_j\}_{j \in \mathbb{Z}}$, 
the embedding $D \subset l^2(\nu)$ is equivalent to 
the weighted norm inequality 
$$
|| K x||_{l^2(\nu)} \le C \, ||x||_{l^2}, \quad \forall x \in \mathcal{S}_0.  
$$
Since $G=K^2$ ($G x=(\kappa*\kappa)*x= \kappa*(\kappa*x)$), the preceding inequality is equivalent to the corresponding 
weighted norm inequality for Green's potentials: 
$$
||G (x \nu)||_{l^2(\nu)} \le C \, ||x||_{l^2(\nu)},  \quad \forall x \in \mathcal{S}_0.  
$$ 
Here we use the notation $x y=(x_k y_k)_{k \in \mathbb{Z}}$. \ \par 
\ \par \noindent 
\bf 4.2. Capacities and equilibrium potentials. \rm   The capacity of a nonempty set $J \subset \mathbb{Z}$ associated with the Besov-Dirichlet space 
$D$ is given by (\cite{FOT2011}, Sec. 2.1): 
$$
\text{Cap} (J) := \inf \left \{ ||x||^2_{D}: \, \, x \in D, \, \, x_{j} \ge 1  \quad \text{if} \, \, j \in J \right \}.    
$$
In this definition we can restrict ourselves to $x$ satisfying $0\le x\le 1$, since 
by the contraction property $||x||_D \ge ||\bar x||_D$, where $\bar x = \min [ \max (x, 0), 1]$.

If the set of $x \in D$ such that $x \ge 1$ on $J$ is empty then we set $\text{Cap} (J) = \infty$. 
This capacity can also be defined by means of the operator $K$ with nonnegative kernel 
$\kappa=\mathcal{F}(1/w^{1/2})$ introduced above: 
$$
\text{Cap} (J) = \inf \left \{ ||y||^2_{l^2}: \, \, y \in l^2, \, \,  y \ge 0, \, \, (Ky)_{j} \ge 1  \quad \text{if} \, \, j \in J   \right \}.    
$$

A general theory of capacities associated with nonnegative kernels is presented in  \cite{AH1996}, 
 Sec. 2. Note that in our case the capacity of a single point set $J_0=\{j_0\}$ is always positive: 
$$
\text{Cap} (J_0) = 1/||1/w||_{L^1(\mathbb{T})}>0.   
$$
It follows from this and the contraction property that there exists a unique 
extremal element $x^J= K y^J\in D$ such that 
$y^J \ge 0$, $0\le x^J \le 1$, $x^J=1$ on $J$, and 
$$
||x^J||^2_{D} = ||y^J||^2_{l^2} = \text{Cap} (J),  
$$
provided $\text{Cap} (J)<\infty$. Moreover (see \cite{AH1996}, 
 Sec. 2.5), at least for finite sets $J$,  the capacity $\text{Cap}(\cdot)$ coincides with the dual Green capacity 
$$
\text{Cap} (J) = \sup \Big \{ \sum_{j \in J} z_j : \, \, {\rm supp} \, z \subset J, \, \, 
 (G z)_{j} \le 1  \, \, \text{if} \, \, j \in J, \quad z \ge 0\Big \},  
$$
where the supremum is taken over all nonnegative sequences $z=(z_j)$ supported on $J$ such that 
$G z \le 1$ on $J$. There exists a unique extremal element $z^J\ge 0$ supported on $J$ such 
$x^J=K y^J=G z^J \in D$,  $G z^J=1$ on $J$, and $G z^J \le 1$ on $\mathbb{Z}$.

In summary (see \cite{AH1996},  \cite{FOT2011}), to each finite nonempty set $J\subset \mathbb{Z}$  one can associate a unique extremal element 
(\textit{equilibrium potential}) 
$x^J \in D$ such that 
$$
x^J= K y^J = G z^J, \quad y^J= K z^J, 
$$
$$x^J=1\, \,  {\rm on} 
 \, \,  J, \quad ||x^J||_{l^\infty} = 1, \quad z^J\ge 0, \quad {\rm supp} \, z^J \subset J, 
$$ 
$$
 \quad ||x^J||^2_{D} = ||y^J||^2_{l^2} =||z^J||_{l^1}= \text{Cap} (J).
 $$
\ \par \noindent 
\bf 4.3. Theorem. \it Let $\nu=(\nu_k)_{k \in \mathbb{Z}}$ be a nonnegative sequence. Then   
the inequality 
$$
\sum_{j\in \mathbb{Z}} |x_k|^2 \nu_k  \le C \, ||x||^2_{D}
$$
holds for all $x \in D$ if and only if 
$$
\sum_{k \in J} \nu_k \le C_1 \, {\rm Cap} (J),
$$
for every finite set $J \subset \mathbb{Z}$, where $C_1\le C \le 4 C_1$. \rm \ \par \ \par

Theorem 4.3 is an immediate consequence of a discrete analogue  
of Maz'ya's strong capacitary inequality 
stated in the following lemma. Its proof 
given below is based on an argument due to 
K. Hansson \cite{Han1979} (see also \cite{Maz2011}, Sec. 11.2.2; \cite{FOT2011}, Sec. 2.4).   
Its main idea is a clever use of  equilibrium potentials 
whose properties were discussed above.  \ \par \ \par\noindent 
\bf 4.4. Lemma. \it Let $x =(x_j)_{j \in \mathbb{Z}}\in \mathcal{S}_0$. For $t>0$,  let $N_t = \{j\in \mathbb{Z}: \, |x_j|\ge t\}$. 
 Then 
$$
\int_0^\infty {\rm Cap} (N_t) \, t \, dt \le 2 \, ||x||^2_D. 
$$ 

\bf Proof. \rm Clearly, the left-hand side of the preceding inequality is finite. 
 Let $x= K y$, where 
$y = \mathcal{F} \big ((1/w)^{1/2} \mathcal{F}^{-1} x\big ) \in l^2$. Notice that $|x| \le K (|y|)$, and 
$y \in l^2$. 
Let  $u^{N_t} = G z^{N_t}$ be the equilibrium potential associated with 
the finite set $N_t$. Here $0\le u^{N_t}\le 1$,  $u^{N_t}=1$ on $N_t$, and   $z^{N_t}\ge 0$,  
  ${\rm supp} \, z^{N_t} \subset N_t$. Since ${\rm Cap} (N_t) = \sum_j z^{N_t}_j$, and $ (K |y|)_j \ge |x_j| \ge t$ on $N_t$, 
we have  
$$
\int_0^\infty {\rm Cap} (N_t) \, t \, dt = \int_0^\infty \sum_j  z^{N_t}_j \,  t \, dt 
\le \int_0^\infty \sum_j  (K |y|)_j   z^{N_t}_j \,   dt$$
$$
 = \int_0^\infty \sum_j  |y|_j   (Kz^{N_t})_j \,   dt=
\sum_j   |y_j|  \int_0^\infty    (Kz^{N_t})_j \,   dt 
$$
$$
 \le ||y||_{l^2} \, \left (\sum_j \left ( \int_0^\infty    (Kz^{N_t})_j \,   dt \right)^{2}\right)^{1/2}. 
$$
We deduce  
$$
\sum_j \left ( \int_0^\infty    (Kz^{N_t})_j \,   dt \right )^{2} 
= 2 \sum_j \int_0^\infty    (Kz^{N_t})_j \,    \int_0^t    (Kz^{N_s})_j \,      ds  dt 
$$
$$
= 2 \int_0^\infty \int_0^t    \sum_j (Kz^{N_t})_j \, (Kz^{N_s})_j \,      ds  dt. 
$$
Since $G$ and $K$ have symmetric kernels, and $G=K^{2}$, we have 
$$
 \sum_j (Kz^{N_t})_j \, (Kz^{N_s})_j  =  \sum_j (G z^{N_s})_j  z^{N_t}_j.   
$$
Here $G z^{N_s}=u^{N_s}$ is the equilibrium potential associated with $N_s$. 
 Consequently, $0 \le (G z^{N_s})_j\le 1$ for 
all $j \in \mathbb{Z}$, and, since $ z^{N_t}$ is supported in $N_t$,  
$$
\sum_j  (G z^{N_s})_j  z^{N_t}_j 
 \le \sum_{j \in N_t} z^{N_t}_j =  {\rm Cap} (N_t). 
$$
Hence, 
$$
\sum_j \left ( \int_0^\infty    (Kz^{N_t})_j \,   dt \right )^{2} 
\le 2 \int_0^\infty   \int_0^t   {\rm Cap} (N_t) \,  ds \, dt  
=
2 \int_0^\infty   
{\rm Cap} (N_t) \, t \, dt. 
$$ 
Combining the preceding inequalities, 
we deduce 
$$
\left ( \int_0^\infty {\rm Cap} (N_t) \, t \, dt \right)^{1/2}  \le 2 ||y||_{l^2} =2 ||x||_D. 
$$
\qed \ \par \ \par

\bf Proof of theorem 4.3.\rm To prove the ``if'' part of Theorem 4.3, 
we assume without loss of generality that $x\in \mathcal{S}_0$, and estimate 
$$
\sum_{j\in \mathbb{Z}} |x_k|^2 \nu_k  = 2 \int_0^\infty \left ( \sum_{k: \, |x_k|\ge t} \nu_k \right) t \, dt 
$$
$$
\le 2 C_1 \int_0^\infty {\rm Cap} (N_t) \, t \, dt  
\le 4C_1 \, ||x||^2_{D}.
$$
The ``only if'' part is obvious. \qed

\ \par \ \par \noindent 
\bf 4.6. A non-capacitary characterization of the embedding $D \subset l^2(\nu)$. \rm The capacitary condition for 
the embedding $D \subset l^2(\nu)$
in Theorem 4.3 can be restated in the equivalent ``energy'' form which does not 
use capacities: 
$$
\sum_{j \in J} \sum_{m\in J} g_{m-j} \nu_j \nu_m \le C \sum_{j \in J} \nu_j,  
$$
for every finite set $J \subset \mathbb{Z}$. (In the continuous case 
this was first noticed by D. R. Adams.) Indeed, denote by $\nu^J = \chi_J \nu$ the 
sequence $\nu$ restricted to $J$. Then, 
if the preceding condition holds, it follows that, 
for every $y \in l^2$ ($y \ge 0$)  such that $K y\ge 1$ on $J$,  
$$
\sum_{j \in J} \nu_j \le \sum_{j\in J} (K y)_j \nu_j = \sum_j y_j (K \nu^J)_j \le ||y||_{l^2} ||K \nu^J||_{l^2} 
$$
$$
=||y||_{l^2} \left (\sum_{j \in J} \sum_{m\in J} g_{m-j} \nu_j \nu_m \right)^{1/2} \le  ||y||_{l^2} 
C^{1/2}   \left ( \sum_{j \in J} \nu_j \right)^{1/2} . 
$$
Consequently, 
$$
 \sum_{j \in J} \nu_j \le C ||y||_{l^2}^2 . 
$$
Minimizing over all such $y$, we obtain 
$$
 \sum_{j \in J} \nu_j \le C \, {\rm Cap}(J) . 
$$ 
Conversely, suppose that the preceding condition holds. Obviously,  
$$
\sum_{j \in J} \sum_{m\in J} g_{m-j} \nu_j \nu_m = || K \nu^J||^2_{l^2} . 
$$
By duality,
$$
|| K \nu^J||_{l^2}= \sup_{y:  \, ||y||_{l^2} \le 1} \left \vert \sum (K \nu^J)_j y_j \right \vert . 
$$
Since $||y||_{l^2} \le 1$, envoking Theorem 4.3 we estimate
$$
\left \vert  \sum_j (K \nu^J)_j y_j \right \vert = \left \vert  \sum_{j \in J} (K y)_j \nu_j \right \vert \le 
||K y||_{l^2(\nu)} 
\left (\sum_{j \in J} \nu_j \right)^{1/2} \le 2 C^{1/2} \left (\sum_{j \in J} \nu_j \right)^{1/2} . 
$$
Thus,
$$
\sum_{j \in J} \sum_{m\in J} g_{m-j} \nu_j \nu_m \le 4 C \sum_{j \in J} \nu_j. 
$$ 
\qed \ \par \ \par\noindent 
\bf 4.7. Quasi-metric Green kernels. \rm Suppose that  the discrete Green's kernel $(g_{j-m})$, where 
$g = \mathcal{F} (w^{-1})$ ($g=(g_j)>0$) has the 
following quasi-metric property: 
$$
1/g_{j+m}\le \varkappa \left ( 1/g_{j} + 1/g_{m}\right ), 
\quad j, m \in \mathbb{Z},  
$$ 
for some constant $\varkappa>0$. 
Then in the energy condition 
$$
 \sum_{j \in J}  \sum_{m \in J} g_{j-m} \nu_j \, \nu_m \le C \, 
 \sum_{j \in J} \nu_j,  
 $$
which characterizes the embedding $D \subset l^2(\nu)$, it suffices to assume that $J$ is a 
quasi-metric ball:
$$J = \{j \in   \mathbb{Z}: \, \, g_{j-m}  >1/r \}, \quad m \in \mathbb{Z},  \, \, r>0. $$ 
Obviously, if $(g_j)$ is nonincreasing for $j\ge 0$ then  
$J$ is an interval:
$ \{ j \in \mathbb{Z}: \, \, n_1 \le j \le n_2\}. 
 $
 
 This is a special case of a general result on quasi-metric kernels. 
  Let $(\Omega, \nu)$ be a measure space.  A symmetric, measurable kernel $G:\Omega \times \Omega \rightarrow (0, + \infty] $ 
  is called quasi-metric if $d = 1/G$ satisfies the quasi-triangle inequality 
  $$
  d(x,y) \leq \varkappa
\left (d(x,z) + d(z,y) \right)
$$
for some $\varkappa > 0$ independent of $x,y,z \in \Omega$. By $B(x,r)=\{ \, y \in \Omega: \, d(x,y)<r \}$ denote the quasi-metric ball of radius 
$r>0$ centered at $x \in \Omega$. Consider the integral operator 
$$
T f (x) = \int_\Omega G(x,y) f(y) d \nu(y), \quad x \in \Omega. 
$$ 
The following theorem is due to F. Nazarov. 
  \ \par \ \par\noindent 
\bf 4.8. Theorem. \it Let $(\Omega, \nu)$ be a measure space with $\sigma$-finite measure $\nu$.   
Let $G$ be a quasi-metric kernel on $\Omega$, 
and let $d = 1/G$ be the corresponding quasi-metric such that $\nu(B)< \infty$ for every 
quasi-metric ball $B=B(x,r)$. Then 
$$
||T f||_{L^2(\Omega, \nu)} \le C \, ||f||_{L^2(\Omega, \nu)}, 
$$
for all $f \in L^2(\Omega, \nu)$,  if and only if there exists a constant $c=c(\kappa)>0$ such that, for every $B$, 
$$\int_B\!\int_B G(x,y) \, d \nu (x) \, d \nu (y) \le C_1 \, \nu(B).$$
 Moreover, 
there exists a constant $c=c(\varkappa)>0$ such that $C/c \le C_1 \le c \, C$.  
\rm  \ \par \ \par

In particular, this theorem is applicable to the weighted norm inequality for the discrete Green's operator 
with kernel $(g_{j-m})$, where $g = \mathcal{F} (w^{-1})$: 
$$
||G (x \nu)||_{l^2(\nu)} \le C \, ||x||_{l^2(\nu)},  \quad \forall x \in \mathcal{S}_0,   
$$ 
provided $d = 1/g$ is a quasi-metric on $\mathbb{Z}$. As was demonstrated above, the preceding  inequality 
is equivalent to the embedding $D \subset l^2(\nu)$. 

The quasi-metric property holds in many important cases, in particular, for Green's kernel 
associated with the weight $w_\alpha(e^{it})= |e^{it}-1|^{\alpha}$ ($-1<\alpha<1$) discussed in the next 
subsection. For such weights, both capacitary and non-capacitary characterizations 
of multipliers ${\rm Mult}(D)={\rm Mult}(L^2(w_\alpha))$ are available. 
\ \par \ \par \noindent 
\bf 4.9. Example: Besov-Dirichlet spaces of fractional order. \rm Let $f \in L^2 (\mathbb T, w_\alpha)$ 
where $w_\alpha(e^{it})= |e^{it}-1|^{\alpha}$ ($0<\alpha<1$). Then  
$$
w_\alpha(e^{it}) \approx 4 \sum_{j \ge 1} c_j \sin^2 (jt/2), \quad \text{where} \, \, c_j = \frac {1} {(j+1)^{1+\alpha}} .
$$
Let $x = \{x_j\}$, where $x_j= \hat f(j)$, $j \in \mathbb Z$. It follows that 
$$ ||f||^2_{L^2 (\mathbb T, w_\alpha)} \approx ||x||^2_{D_{\alpha/2}} = \sum_{j \in \mathbb Z} \sum_{m \in \mathbb Z} 
\frac{|x_j-x_{m}|^2}{ (|j-m|+1)^{1+\alpha}}. $$ 
Here $D_{\alpha/2}= D=\mathcal{ B}^{2}_{0}(c_{\vert 
j-m\vert }/2)$ is a Besov-Dirichlet space on $\mathbb Z$ of fractional order $\alpha$.

Next, a sequence $\lambda=(\lambda_j)_{j \in \mathbb Z}$ is a multiplier of $D_{\alpha/2}$:
$$
||\lambda  x||_{D_{\alpha/2}} \le C ||x||_{D_{\alpha/2}}, \quad \forall x \in D_{\alpha/2}, 
$$ 
if and only $\lambda \in {\rm Mult} \left(L^2 (w_\alpha)\right)$, 
or by duality  $\lambda \in {\rm Mult} \left(L^2 (w_{-\alpha})\right)$. 

One can rewrite this condition using the discrete Riesz potential 
$\mathcal{R}_{\alpha /2}$ of order $ \alpha /2$, where 
$$
\mathcal{R}_{\alpha /2}y:= \Big (\sum _{j\in \mathbb{Z}}{\frac{y_{j}}{(\vert 
m-j\vert +1)^{1-\alpha /2}}} \Big )_{m\in \mathbb{Z}} ,
$$
in the following way (letting $ x= \mathcal{R}_{\alpha /2}y$, 
$ y\in \mathcal{S}_{0}$, and taking into account that $\mathcal{R}_{\alpha /2}\mathcal{S}_{0}$ 
is dense in $ D_{\alpha /2}$): 
$$
||\lambda \, (\mathcal{R}_{\alpha/2} \, y)||_{D_{\alpha/2}} \le C ||y||_{l^2(\mathbb Z)} ,   \quad  
\forall y \in \mathcal{S}_0 .  
$$

The corresponding Green kernel $g=(g_{j-m})$ is equivalent to the Riesz kernel of order $\alpha$ since 
$g_j \approx 1/(|j|+1)^{1-\alpha}$, $j \in \mathbb{Z}$.

By Lemma 2.2 (1), 
$\lambda \in {\rm Mult} (D_{\alpha/2}) $ if and only if 
$\lambda \in l^\infty(\mathbb Z)$, and the sequence $\mu=(\mu_j )_{j \in \mathbb Z}$ defined by:  
$$
\mu_j =\mu_j^\alpha(\lambda):= \left (\sum_{m \in \mathbb Z}   \frac{| \lambda_j - \lambda_m|^2} 
{(|j-m|+1)^{1+ \alpha}}\right)^{1/2}, \quad j \in \mathbb Z,
$$ 
is a multiplier from $D_{\alpha/2}$ to $l^2$:
$$
\sum_{j \in \mathbb Z} \mu_j^{\alpha}(\lambda)^2 \, |x_j|^2 \le C^2 ||x||^2_{D_{\alpha/2}} , \quad \forall x \in D_{\alpha/2} . 
$$
Equivalently, $D_{\alpha/2} \subset l^2(\nu)$ where $\nu = \mu^2$. 
The preceding inequality holds if and only if $\mu$ obeys the capacitary 
condition of Theorem 4.3: 
$$
\sum_{j \in J} \mu_j^{\alpha}(\lambda)^2 \le C \, {\rm Cap}_\alpha (J) , 
$$ 
for every finite set $J \subset \mathbb{Z}$. Here ${\rm Cap}_\alpha (\cdot)$ is the capacity 
 associated with the Besov-Dirichlet space $D=D_{\alpha/2}$ (see Sec. 4.2).

Since the corresponding Green kernel $(g_{j-m})$ has the  quasi-metric property,  
and $g_j$ is decreasing for $j\ge 0$, it follows from Theorem 4.8 that the embedding 
$D_{\alpha/2} \subset l^2(\nu)$ ($\nu = \mu^2$) is equivalent to the energy 
condition 
$$
\sum_{j \in J} \sum_{m \in J} \frac{ \mu_j^{\alpha}(\lambda)^2 \, \mu_m^{\alpha}(\lambda)^2} {(|j-m|+1)^{1-\alpha}} \le C \sum_{j\in J} \mu_j^{\alpha}(\lambda)^2,
$$ 
for every \textit{interval} $J$ in $\mathbb{Z}$. This characterization of multipliers $\mu\colon D_{\alpha/2} \longrightarrow l^2$ is due to Kalton and Tzafriri \cite{KT1998}. 
\ \par \ \par\noindent 
\bf 4.10. Multipliers in pairs of Besov-Dirichlet spaces.\rm 
To treat weights with several power-like singularities considered below, we will need classes  of multipliers acting from $D_{\beta/2}$ to $D_{\alpha/2}$. They will be  
denoted by ${\rm Mult} (D_{\beta/2}\longrightarrow D_{\alpha/2})$; for $\alpha=\beta$ we will continue 
to use the notation ${\rm Mult} (D_{\alpha/2})$. We remark that 
${\rm Mult} (D_{\beta/2}\longrightarrow D_{\alpha/2})$ coincides with the class of Fourier multipliers 
${\rm Mult}(L^2(w_\beta) \longrightarrow  L^2(w_\alpha))$ defined in a similar way. The following characterization 
of multipliers is similar to the continuous case (see \cite{MSh2009}), but there are certain differences which we 
need to take into account (see Remark 4.12 below). \ \par \ \par\noindent 
\bf 4.11. Theorem. (1) \it  Let $0< \beta \le \alpha <1$. Then $\lambda \in {\rm Mult} (D_{\beta/2}\longrightarrow D_{\alpha/2})$ 
if and only if $\lambda \in l^\infty$, and $\mu \in  {\rm Mult} (D_{\beta/2} \longrightarrow   l^2)$, where 
$\mu=(\mu_j)$ is defined by 
$$
\mu_j=\mu_j^\alpha(\lambda):= \left (\sum _{m} \frac{\vert \lambda_{j}-\lambda_{m}\vert ^{2}}{(\vert j-m\vert +1)
^{1+\alpha }}\right)^{1/2}, \quad j \in \mathbb{Z} . 
$$ 
Equivalently, 
$$
\lambda \in l^\infty, \quad {\rm and}\quad 
\sum_{j \in J} \mu_j^{\alpha}(\lambda)^2 \le C \, {\rm Cap}_\beta (J) , 
$$ 
for every finite set $J \subset \mathbb{Z}$, where $C$ does not depend on $J$. \ \par \ \par  \noindent 
\bf (2) \it Let $0< \alpha < \beta <1$. Then $\lambda \in {\rm Mult} (D_{\beta/2} \longrightarrow   D_{\alpha/2})$ if and only if $\lambda \in {\rm Mult} (D_{(\beta-\alpha)/2} \longrightarrow   l^2)$, and 
$\mu \in  {\rm Mult} (D_{\beta/2} \longrightarrow   l^2)$. Equivalently, 
$$
\sum_{j \in J} \, |\lambda_j|^2 \le C \, {\rm Cap}_{\beta-\alpha} (J) , \quad {\rm and}\quad  
\sum_{j \in J} \mu_j^{\alpha}(\lambda)^2 \le C \, {\rm Cap}_\beta (J) , 
$$ 
for every finite set $J \subset \mathbb{Z}$, where $C$ does not depend on $J$.  
\rm \ \par \ \par

\bf Proof. \rm If $\lambda \in {\rm Mult}(L^2(w_\beta) \longrightarrow  L^2(w_\alpha))$, 
and $\Vert 
T_{\lambda }\Vert=\Vert 
T_{\lambda }\Vert_{L^2(w_{\beta})  \longrightarrow L^2(w_{\alpha})}$ is the multiplier norm, 
then 
for all $n \in \mathbb{Z}$,  
$$
\vert \lambda _{n}\vert \cdot \Vert w_{\alpha }\Vert_{L^1(\mathbb{T})}^{1/2}=
\Vert T_{\lambda }z^{n}\Vert_{L^{2}  (w_{\alpha })}\leq  \Vert 
T_{\lambda }\Vert  \Vert z^{n}\Vert_{L^{2}  (w_{\beta})}= \Vert T_{\lambda }\Vert \cdot \Vert w_{\beta 
}\Vert_{L^1(\mathbb{T})}^{1/2}. 
$$
Consequently, $\lambda \in l^\infty$.

\rm (1) Suppose $0< \beta \le \alpha <1$, and $\lambda \in {\rm Mult}(L^2(w_\beta) \longrightarrow  L^2(w_\alpha))$. Since $\lambda \in l^\infty$, we have  
$$
\sum_{j \in \mathbb{Z}} \,  |\lambda_j|^2 \mu_j^\alpha(x)^2=\sum_{j \in \mathbb{Z}} |\lambda_j|^2  \sum_{m\in \mathbb{Z}} \frac{\vert x_{j}- x_{m}\vert ^{2}}{(\vert j-m\vert +1)
^{1+\alpha }}\le C \, ||x||^2_{D_{\alpha/2}} \le C \, ||x||^2_{D_{\beta/2}}. 
$$
Hence, as in the case $\alpha=\beta$ considered above, we see that $\lambda \in {\rm Mult}(L^2(w_\beta) \longrightarrow  L^2(w_\alpha))$ if and only if $\lambda \in l^\infty$, and 
$$
\sum_{j \in \mathbb{Z}} \,  |x_j|^2 \mu_j^\alpha(\lambda)^2=\sum_{j \in \mathbb{Z}} |x_j|^2 \sum_{m\in \mathbb{Z}} \frac{\vert \lambda_{j}- \lambda_{m}\vert ^{2}}{(\vert j-m\vert +1)
^{1+\alpha }}\le C \, ||x||^2_{D_{\beta/2}} , 
$$
for all $x \in D_{\beta/2}$, or equivalently $\mu \in  {\rm Mult} (D_{\beta/2} \longrightarrow   l^2)$. 
By Theorem 4.3 the preceding condition is equivalent to 
$$
\sum_{j \in J} \mu_j^{\alpha}(\lambda)^2 \le C \, {\rm Cap}_\beta (J) , 
$$ 
for every finite set $J \subset \mathbb{Z}$. \ \par \ \par

(2) Suppose $0< \alpha < \beta <1$, and $\lambda \in {\rm Mult}(L^2(w_\beta) \longrightarrow  L^2(w_\alpha))$. Then by duality, $\lambda \in {\rm Mult}(L^2(w_{-\alpha}) \longrightarrow  L^2(w_{-\beta}))$. Hence, using  interpolation 
for operators acting in $L^2$ spaces with weights, we obtain $\lambda \in {\rm Mult}(L^2(w_{\beta-\alpha}) \longrightarrow  L^2).$ By Theorem 4.3 this implies 
$$
\sum_{j \in J} \, |\lambda_j|^2 \le C \, {\rm Cap}_{\beta-\alpha} (J) , 
$$ 
for all  $J \subset \mathbb{Z}$.  On the other hand, as in the case $\alpha\ge \beta$, we have that $\lambda \in {\rm Mult}(L^2(w_\beta) \longrightarrow  L^2(w_\alpha))$ if, for all $x \in D_{\beta/2}$,  the following pair of inequalities hold:
$$
\sum_{j \in \mathbb{Z}} \,  |x_j|^2 \mu_j^\alpha(\lambda)^2 \le C ||x||^2_{D_{\beta/2}} ,  
$$ 
$$
\sum_{j \in \mathbb{Z}} \,  |\lambda_j|^2 \mu_j^\alpha(x)^2 \le C ||x||^2_{D_{\beta/2}} .  
$$ 
Moreover, if the second inequality holds, then the first one is necessary in order that 
$\lambda \in {\rm Mult}(L^2(w_\beta) \longrightarrow  L^2(w_\alpha))$. 

It remains to show that if $\lambda \in {\rm Mult}(L^2(w_{\beta-\alpha}) \longrightarrow  L^2)$, then 
the last inequality holds. 
  Clearly, if $\lambda \in {\rm Mult}(L^2(w_{\beta-\alpha}) \longrightarrow  L^2)$ 
is a bounded Fourier multiplier, or equivalently, $\lambda \in {\rm Mult} (D_{(\beta-\alpha)/2} \longrightarrow   l^2)$, 
it follows that  
$$
|| \lambda \, \mu^{\alpha} (x) ||_{l^2} \le C \, ||\mu^{\alpha} (x)||_{D_{(\beta-\alpha)/2}} , 
$$
where $\mu^{\alpha} (x) =(\mu_j^\alpha(x))$. 
We note that by definition 
$$
 ||\mu^{\alpha} (x)||_{D_{(\beta-\alpha)/2}} = || \mu^{\beta-\alpha} [ \mu^{\alpha} (x)] ||_{l^2} .  
$$
Let us show that 
$$ ||\mu^{\beta-\alpha} [ \mu^{\alpha} (x)] ||_{l^2}  \le C  ||\mu^{\beta} (x)||_{l^2} , $$
 where $C$ 
depends only on $\alpha, \beta$. This is a discrete analogue of Lemma 4.2.1 \cite{MSh2009}: by the triangle 
inequality, 
$$
 ||\mu^{\beta-\alpha} [ \mu^{\alpha} (x)] ||^2_{l^2}
  \le \sum_{k \in \mathbb{Z}}\sum_{j \in \mathbb{Z}} \sum_{m \in \mathbb{Z}} \, 
 \frac{| x_{j+k} -x_k +x_{m+k} -x_{j+k+m }|^2} {(|j|+1)^{1+\alpha} (|m|+1)^{1+\beta-\alpha} } .
 $$ 
 The triple sum on the right-hand side is symmetric with respect to $j$ and $m$, and so 
 it is enough to consider the  case  $|m|\ge |j|$. 
 We estimate 
 $$
  \sum_{k \in \mathbb{Z}}\sum_{j \in \mathbb{Z}} 
\frac{| x_{j+k} -x_k|^2}{(|j|+1)^{1+\alpha}}    \sum_{|m|\ge |j|} \frac{1}{ (|m|+1)^{1+\beta-\alpha}} 
\le C \,  \sum_{k \in \mathbb{Z}}\sum_{j \in \mathbb{Z}} 
\frac{| x_{j+k} -x_k|^2}{(|j|+1)^{1+\beta}}  = C ||\mu^{\beta} (x)||^2_{l^2} . 
$$
Similarly, interchanging the order of summation in the remaing term and replacing $k$ with $n=k+m$, 
we estimate  
$$
\sum_{k \in \mathbb{Z}}\sum_{j \in \mathbb{Z}} \frac{1}{(|j|+1)^{1+\alpha}} \sum_{|m|\ge |j|}  \frac{|x_{m+k} -x_{j+k+m }|^2} {(|m|+1)^{1+\beta-\alpha} } \le C ||\mu^{\beta} (x)||^2_{l^2} . 
 $$ 
 
 Combining the above estimates and taking into account that $||\mu^{\beta} (x)||_{l^2} = 
 ||x||^2_{D_{\beta/2}}$, 
 we conclude that 
 $\lambda \in {\rm Mult}(D_{(\beta-\alpha)/2}\longrightarrow L^2)$ implies 
$$
|| \lambda \, \mu^{\alpha} (x) ||_{l^2} \le C \, ||x||^2_{D_{\beta/2}} , 
$$
for all $x \in D_{\beta/2}$. Thus, $\lambda \in {\rm Mult}(D_{\beta/2}\longrightarrow D_{\alpha/2})$ if and only if 
$\mu=\mu^\alpha(\lambda) \in  {\rm Mult} (D_{\beta/2} \longrightarrow   l^2)$ and $\lambda \in {\rm Mult}(D_{(\beta-\alpha)/2} \longrightarrow  l^2)$.  By Theorem 4.3 both conditions have the corresponding capacitary characterizations (see  Sec. 4.9).  
\qed
\ \par \ \par\noindent 
\bf 4.12. Remark. \rm In the continuous case of multipliers in pairs of Sobolev spaces $W^{\alpha, p}(\R^n)$,  it is known that 
${\rm Mult}(W^{\beta, p}(\R^n)\longrightarrow W^{\alpha, p}(\R^n)) = \{0\}$ if $\alpha> \beta>0$  (see \cite{MSh2009}, Sec. 2.1), contrary to the discrete case ${\rm Mult} (D_{\beta/2} \longrightarrow   D_{\alpha/2})$. 

%
%
%
%

\section{Weights with several {\rm LKS} singularities}\label{sec5}

In this section we are concerned with weights which are equivalent to products, or sums of the reciprocals of {\rm LKS} weights. 
Such weights may have finitely many singularities, and generally are no longer {\rm LKS}-weights. Nevertheless, we will be able to characterize multipliers ${\rm Mult} \, (L^2(w))$, and show that   they obey the Spectrum Localization Property ({\rm SLP}). It is known that the  {\rm SLP} fails for weights with infinitely many singularities of this type (see 5.12-5.13 below). 

In this section it will be convenient to use the following notation for weights on $\mathbb{T}$: $w_\alpha=|e^{it}-1|^\alpha$, and 
$w_\alpha^{\theta}= |e^{it}-e^{i \theta}|^\alpha$ ($\alpha \in \mathbb{R}$). We will consider weights of the type 
$$ 
w= \displaystyle \prod _{j=0}^{d-1}w_{\alpha _{j}}^{\theta_{j}} ,
$$ 
where $\theta_{j}\in  \mathbb{R}$ are pairwise distinct $({\rm mod} (2\pi))$ points.
\ \par \ \par \noindent 
\bf 5.1. Weights with several singularities of the same order.

\ \par 
\noindent \bf Theorem. \it  Let $ w=$ $ \displaystyle \prod _{j=0}^{d-1}w_{\alpha }^{\theta_{j}}$,  
where $\theta_{j}$ are pairwise distinct $({\rm mod} (2\pi))$ points.
Let 
$ d=d_{1}>d_{2}>...>d_{n}= 1$ be all distinct divisors of $d$. We denote by $d_s$ 
the largest among the 
divisors  such that the set of singularities $ \sigma =$ $ \{e^{i \theta_j}:$   
$ 0\leq j<d-1\}$ is a union of vertices of $n_s=d/d_s$ distinct regular $d_s$-sided polygons. Then ${\rm Mult} (L^{2}(w))={\rm Mult} (L^{2}(w_{\alpha }(e^{it d_{s}}))$, 
and $ \lambda \in {\rm Mult} (L^{2}(w))$ if and only if $\lambda \in l^\infty$, and 
$$
\sum_{j \in J} \underset{d_s \vert \, j -m} {\sum_{m \in \mathbb{Z}}} \, 
  \frac{|\lambda_j - \lambda_m|^2}{(|j-m|+1)^{1+\alpha}} \le  C \, {\rm Cap}_\alpha (J) , 
$$ 
for every finite set $J \subset \mathbb{Z}$.
\ \par 
\ \par 
\noindent \bf Remarks.  (1) \rm If $d_s=1$, i.e.,  $\sigma$ has no rotational symmetries 
(the generic case), then $ {\rm Mult} (L^{2}(w))= {\rm Mult} (L^{2}(w_{\alpha 
}))$. \ \par \ \par 
\noindent \bf (2) \rm If $d$ is a prime number then either $d_s=d$ or $d=1$. However, generally  
$d_s$  is not necessarily a prime number. For instance, 
if $d=12$, then $\sigma$ may consist of the vertices of either a single regular dodecagon,  or two regular hexagons, 
or three squares, 
or four regular triangles, or six $2$-gons, i.e., pairs of opposite points. Then $d_s=12, 6, 4, 3, 2$, respectively, and there are $n_s=d/d_s$ distinct regular polygons, so that each of them is a rotation 
of the set of roots of unity of order $d_s$. If $d_s=1$, then $\sigma$ consists of $12$ points on $\mathbb{T}$ 
with no rotational symmetry. 
 \rm \ \par \ \par 
\bf Proof. \rm 1. For the sake of simplicity, we  first consider the case of two singularities. Let $w=w_{-\alpha} w_{-\alpha}^{\theta}$, where $0<\alpha<1$ and $\theta\not=0$ ($\text{mod} (2\pi)$). Note that 
$$w\approx w_{-\alpha} + w_{-\alpha}^{\theta}= |e^{it}-1|^{-\alpha} + |e^{it}-e^{i \theta}|^{-\alpha} .$$ 
By duality ${\rm Mult} \, (L^2(w))$= ${\rm Mult} \, (L^2(1/{\widetilde w}))$, where $\widetilde w(e^{i t}) = w(e^{-i t})$.  
Hence, at the same time 
we  obtain a characterization of multipliers for weights of the type 
$$1/w = 
w_{\alpha} w_{\alpha}^{\theta} = |e^{it}-1|^{\alpha} |e^{it}-e^{i\theta}|^{\alpha} .$$

 It will be more convenient to work with 
convolution operators $T_\lambda f=k\star f$ on $L^2(w)$, where $\lambda = (\lambda_j)_{j \in \mathbb{Z}} = \mathcal{F} k \in l^\infty(\mathbb{Z})$, and consequently  $k$ is a pseudo-measure on $\mathbb{T}$. 
 For a pseudo-measure $k$ on $\mathbb{T}$,  it follows that $\lambda= \mathcal{F} k  \in {\rm Mult} \, (L^2(w))$  
if and only if 
$$||k\star f||_{L^2(w)} \le C ||f||_{L^2(w)},$$
for all trigonometric polynomials $f$, which is equivalent to a pair of inequalities:
$$||k\star f||_{L^2(w_{-\alpha})} \le C ||f||_{L^2(w_{-\alpha} w_{-\alpha}^{\theta})}, $$
$$||k\star f||_{L^2(w_{-\alpha}^{\theta})} \le C ||f||_{L^2(w_{-\alpha} w_{-\alpha}^{\theta})}.$$ 
Using the rotation operator $R_{\theta} f (e^{it}) =  f (e^{i (t+\theta)} )$, and letting $g= R_{\theta} f$, we see that the second inequality is  
equivalent to  
$$||k\star g||_{L^2(w_{-\alpha})} \le C ||g||_{L^2(w_{-\alpha} w_{-\alpha}^{-\theta})} ,$$ 
for all $g \in L^2(w_{-\alpha} w_{-\alpha}^{-\theta})$.  

By duality, we deduce that $\lambda \in {\rm Mult} \, (L^2(w))$ if and only if the following pair of 
inequalities hold: 
$$|| \widetilde k\star f||_{L^2(w_{\alpha} w_{\alpha}^{\theta})} \le C ||f||_{L^2(w_{\alpha})},$$ 
$$|| \widetilde k \star f||_{L^2(w_{\alpha} w_{\alpha}^{-\theta})} \le C ||f||_{L^2(w_{\alpha})},$$
for all $f \in L^2(w_\alpha)$, where $\widetilde k (e^{ix}) = k(e^{-ix})$. 

Let $W=  w_{\alpha} (w_{\alpha}^{\theta}+ w_{\alpha}^{-\theta})$.   Adding up the preceding displayed inequalities, we obtain that they are equivalent to:
$$|| \widetilde k\star f||_{L^2(W)} \le C ||f||_{L^2(w_\alpha)}, \quad \forall f \in  L^2(w_\alpha) .$$
Clearly, if $\theta \not = \pi  + 2\pi n$ ($n \in \mathbb{Z}$), i.e. in the generic case,  
we have 
$$ W(e^{it})=  |e^{it}-1|^{\alpha} | ( |e^{it}-e^{i\theta}|^{\alpha} +  |e^{it}-e^{-i\theta}|^{\alpha} ) 
\approx  |e^{it}-1|^{\alpha}   = w_\alpha .$$
Hence in this case $ T_\lambda = k \star (\cdot)$ is a bounded operator in $L^2(w)$ if and only if  
$\widetilde k \star (\cdot)$ is  a bounded operator in $L^2(w_\alpha)$, or equivalently  $ \lambda \in 
{\rm Mult} \, (L^2(w_\alpha))$, since $\widetilde{w}_\alpha=w_\alpha$. Thus,  multipliers for weights 
with two generic singularities are the same as for weights  with one singularity (characterized 
in Section \ref{sec4}). 

In the non-generic case $\theta  = \pi  + 2\pi n$ ($n \in \mathbb{Z}$), we have 
$$
W(e^{it})  \approx  |e^{it}-1|^{\alpha} |e^{it}+1|^{\alpha} =|e^{2it}-1|^\alpha,
$$
which is equivalent to an {\rm LKS} weight. It follows that $\lambda \in {\rm Mult} \, (L^2(w))$ if and only if 
$\lambda \in l^\infty(\mathbb{Z})$, and 
$$
\underset{j -m\, \rm{even}} {\sum_{j \in \mathbb{Z}} \sum_{m \in \mathbb{Z}}} \, 
  \frac{|\lambda_j x_j - \lambda_m x_m|^2}{(|j-m|+1)^{1+\alpha}} \le C  \, ||x||^2_{D_{\alpha/2}} , \quad \forall  x \in D_{\alpha/2}.
$$ 
Since $\lambda  \in l^\infty(\mathbb{Z})$, it is easy to see, using the same argument 
as in the case of  the weight $|e^{it}-1|^\alpha$, that the preceding inequality holds if and only if 
 $$
\underset{j -m\, \rm{even}} {\sum_{j \in \mathbb{Z}} \sum_{m \in \mathbb{Z}}} \, 
 |x_j|^2  \frac{|\lambda_j- \lambda_m|^2}{(|j-m|+1)^{1+\alpha}} \le C  \, ||x||^2_{D_{\alpha/2}} , \quad \forall  x \in D_{\alpha/2} .
$$ 
Letting 
$$
\mu_j^{\alpha}(\lambda)^2 = \underset{j -m\, \rm{even}} {\sum_{m \in \mathbb{Z}}} \, 
  \frac{|\lambda_j - \lambda_m|^2}{(|j-m|+1)^{1+\alpha}} , \quad j \in   \mathbb{Z} , 
$$
we see that the multiplier problem is reduced 
to the inequality 
$$
\sum_{j \in \mathbb{Z}} \mu_j^{\alpha}(\lambda)^2 \,  |x_j|^2 \le  C  \, ||x||^2_{D_{\alpha/2}}, \quad \forall  x \in D_{\alpha/2} .
$$ 
 Inequalities of this type have been characterized 
in terms of Besov-Dirichlet capacities, or energies associated with $D_{\alpha/2}$ (Sec. 4.9). 
\ \par \ \par 
 \noindent 2. A similar argument works for any number of singularities $d$. Let $J_d=\{0, 1, \ldots, d-1\}$, and let 
 $\sigma=\{e^{i \theta_j}\}_{j \in J_d}$ be the set of singularities 
 of the weight $w=\prod_{j \in J_d} w_{-\alpha}^{\theta_j}$, where $\theta_j \in [0, 2 \pi)$ are 
 pairwise distinct points, and 
 $0 < \alpha <1$. Note that $w \approx \sum_{j \in J_d} w_{-\alpha}^{\theta_j}$. (Analogous results for weights of the type $w=\prod_{j \in J_d} w_{\alpha}^{\theta_j}$ follow by duality.)

 Using the same argument as 
 in the case $d=2$, it is easy to see that 
$\lambda  \in {\rm Mult} \, (L^2(w))$ if and only if the following  
inequality holds: 
$$||\widetilde k\star f||_{L^2(W)} \le C ||f||_{L^2(w_\alpha)}, \quad \forall f \in L^2(w), $$ 
where 
$$
W = \sum_{j\in J_d} \prod_{m\in J_d} w_{\alpha}^{\theta_m-\theta_j} .  
$$

To complete the proof of the Theorem we will need the following lemma which describes the set of singularities of $W$ on $\mathbb{T}$. 
As we will see, this question is related 
to actions of the group of rotations on the finite set $\sigma\subset \mathbb{T}$. 
\ \par \ \par 
\noindent \bf Lemma. \it Let $ W=\sum _{j\in J_{d}}\prod 
_{m\in J_{d}}w_{\alpha }^{\theta _{m}  -\theta _{j}}$ be the weight defined 
above.  The following alternative holds.\ \par  
\noindent (i) Either there exists a divisor $D$ of $ d$, $ 1<D\leq d$, such 
that $ \sigma = \{e^{i\theta _{j}}\}_{j\in J_{d}}$ is the union 
of $ d/D$ (different) regular $ D$-sided polygons, and then, denoting 
by $ d_{s}$ the maximal possible such $D$, we have $ W\approx \vert 
1-e^{id_{s}  t}\vert ^{\alpha }$, \ \par 

\noindent (ii) or $ W\approx \vert 1-e^{it}\vert ^{\alpha }$.\ \par \ \par 
\noindent
\bf Remark. \rm We wish to thank Stephen Montgomery-Smith for pointing out 
that this Lemma and its proof given below are 
related to the orbit-stabilizer theorem and Burnside lemma (see \cite{Ja1985}, Sec. 1.12), and 
can be generalized 
to arbitrary abelian groups.\ \par \ \par

Let us complete the proof of the Theorem assuming the Lemma. 
In case (i) of the Lemma,   $\sigma$ is the union 
of $n_s=d/d_s$ (different) regular $d_s$-sided polygons, and $W$ is equivalent to the  {\rm LKS} weight 
$w_\alpha (e^{i \theta d_s})=|e^{i \theta d_s}-1|^\alpha$. Letting  
$$
\mu_j^{\alpha}(\lambda)^2 = \underset{d_s \, \vert \, j -m } {\sum_{m \in \mathbb{Z}}} \, 
  \frac{|\lambda_j - \lambda_m|^2}{(|j-m|+1)^{1+\alpha}} , \quad j \in   \mathbb{Z} , 
$$
and using the same argument as above we see that 
in this case $\lambda  \in {\rm Mult} \, (L^2(w))$ 
if and only if $\lambda  \in l^\infty(\mathbb{Z})$, and 
$$
\sum_{j \in \mathbb{Z}} \mu_j^{\alpha}(\lambda)^2  \, |x_j|^2 \le  C  \, ||x||^2_{D_{\alpha/2}}, \quad \forall  \, x \in D_{\alpha/2}.
$$ 
In case (ii),  $\sigma$ has no rotational symmetry, and $W \approx w_\alpha$, so that ${\rm Mult} \, (L^2(w))= {\rm Mult} \, (L^2(w_\alpha))$, and the preceding characterization holds with 
$d_s=1$. \qed \ \par 
\ \par

\bf Proof of the Lemma. \rm  1. Without loss of generality we may assume that $ \theta_{0}=
0$. Clearly, $t=0$ is a zero of $W(e^{i t})$. Suppose there exists $t \in (0, 2 \pi)$ such that 
$W(e^{i t})=0$. Then for every $j \in J_d$ there exists $m \in J_d$ such that  
$w_{\alpha}^{\theta_m-\theta_j}(e^{i t})=0$. In other words, there is  a permutation 
$j \mapsto m(j)$ of $J_d$ such that, for all $j \in J_d$, we have  $\theta_{m(j)}= \theta_j + t$ (${\rm mod} \, (2 \pi)$), where $m(j)$ is unique, 
and  
$m(j_1) \not=m(j_2)$ if $j_1 \not= j_2$. Obviously,  $m(j)\not= j$ since $t\not=0$. 
 Adding together these equations for all $j \in J_d$, we see that $t \, d = 0$ (${\rm mod} \, (2 \pi)$), 
i.e., $t = 2 \pi n/ d$ for some $n \in J_d$. It follows that
 $$
 \theta_{m(j)} = \theta_j + \frac{2 \pi n}{d} \, \,  \, ({\rm mod} \, (2 \pi)) ,     \quad  \forall \, j \in J_d . 
 $$ 
 Moreover, for the consecutive iterations of the map $j \mapsto m(j)$ defined by  $m^{(0)}(j) = j$, and $m^{(k)} (j) =m(m^{(k-1)}(j))$ ($k=1, 2, \ldots$),  we see that all  $e^{i \theta_{m^{(k)}(j)}} \in \sigma$, where 
 $$
 \theta_{m^{(k)}(j)} = \theta_j + \frac{2 \pi k n}{d} \, \,  \, ({\rm mod} \, (2 \pi)) ,  
 \quad  \forall \, j \in J_d, \quad  k=0, 1, 2, \ldots .  
 $$

 \ \par \noindent 2. Suppose first that $d$ is a prime number. Then $z_k = m^{(k)} (0) 
 = e^{2 \pi i k n/d}\in \sigma$ ($k=0, 1, 2, \ldots$) are obviously 
 roots of unity of order $d$. It is easy to see that $z_k$ are distinct for $k \in J_d$. Indeed, if $k_1, 
 k_2 \in J_d$, and 
 $z_{k_1} = z_{k_2}$ for $k_2 > k_1$, we have 
 $$
 2 \pi (k_2-k_1) \frac{n}{d}  = 0\quad ({\rm mod} \, (2 \pi)) . 
 $$
 Since $d$ is a prime number, and $0 < k_2-k_1 \le d-1$, it follows that  $n=0$ 
 and consequently $t=0$, which is a contradiction. 
Hence $\sigma=\{z_k\}_{k \in J_d}$ consists of all the roots of unity of order $d$.  In this case 
obviously $W(z)=0$ if and only if $z$ is a root of unity of order $d$, and $W$ has  no repeated zeros. 

Thus, in the non-generic case,  $\sigma$ is  the set of  vertices 
of a regular $d$-sided polygon, $d_s=d$, and $W$ is equivalent to  
$w_\alpha (e^{i \theta d})=|e^{i \theta d}-1|^\alpha$. 

If  $\sigma$ is not the set of  vertices 
of a regular $d$-sided polygon (the generic case), then $d_s=1$, and $W$ equals zero only at $e^{it}=1$, so that $W \approx w_\alpha$. 
 
 \ \par  \noindent 3. 
If $d$ is not a prime number, denote by $d=d_1 > d_2 > \ldots >d_N=1$ all the divisors 
of $d$, and let $n_s=d/d_s$ ($s=1, \ldots, N$).  As was shown above, if 
 $W(e^{it})=0$ for $t \in (0, 2 \pi)$, then  $t= 2 \pi n/d$ for some $n \in J_d$, and there exists a permutation 
$j \mapsto m(j)$ of $J_d$ such that,  
$$
 \theta_{m(j)} = \theta_j + \frac{2 \pi n}{d} \quad ({\rm mod} \, (2 \pi)) ,   \quad  \forall \, j \in J_d . 
 $$ 
 Since every permutation can be decomposed into a union of disjoint cycles, and 
 $ \theta_j  \mapsto \theta_{m(j)}$ is a rotation with the fixed angle $t$, it follows that 
 all cycles in this decomposition must be of the same length $l$. Consequently, the length of 
 the cycle must be  a divisor of $d$, i.e., $l=d_k$ for some $k=1, \ldots, N-1$, and there 
 are $n_k=d/d_k$ disjoint cycles in the decomposition. Geometrically this means that 
 $\sigma=\{e^{i \theta_j}\}_{j\in J_d}$ is the set of vertices of a union of $n_k=d/d_k$ distinct  regular $d_k$-sided polygons. In this case the admissible  values of $n \in J_d$ are $n=0, n_k, 2 n_k, \ldots, (d_k-1)n_k$, and the corresponding admissible values of $\, t = 2\pi j/d_k$ ($j \in J_{d_k}$); i.e.,  $\{e^{it}\}$ are 
 the roots of unity of order $d_k$.

 Let us denote by $d_s$ the length of the largest cycle $l$ (for all 
 possible values of $n \in J_d$), and set $n_s=d/d_s$. Then $\sigma$ is a union of the set of vertices of $n_s$ regular 
 $d_s$-sided polygons, and $d_s$ is the largest such number. In this case,  $k \ge s$, so that 
 $d_s \ge d_k$, 
 and $t= 2 \pi/d_k$ must be a root of unity of order $d_s$. From this it follows that $d_k$ 
 is a divisor of $d_s$, and consequently all admissible  values of $n \in J_d$ are $n=0, n_s, 2 n_s, \ldots, (d_s-1)n_s$.  In other words, the zero set of $W$ coincides  with the roots of unity of order $d_s$, i.e.,   $W$ is equivalent to $w_\alpha(e^{it d_s})=|e^{itd_s}-1|^\alpha$ 
 ($s= 1, 2, \ldots, N-1$). In the generic case $d_s=1$, the points in $\sigma$ \textit{cannot} be represented as  the set of vertices of a union 
of $n_k$ regular $d_k$-sided polygons for any $k=1, 2, \ldots, N-1$. Then $W$ has the only zero 
at $e^{it}=1$, and consequently, $W \approx w_\alpha$. 
 This completes the proof of the Lemma. \qed 
 \ \par \ \par \noindent 
\bf 5.2. Weights with singularities of different orders. \rm For the sake of simplicity let us consider a weight with two generic zeros on $\mathbb{T}$: 
$$w = w_{\beta} w_\alpha^{\theta} =  |e^{it}-1|^\beta |e^{it}-e^{i\theta}|^\alpha, \quad \theta \not= 2 \pi n, \, \, \forall n \in \mathbb{Z} .$$  Without loss of generality we assume 
$0<\alpha \le \beta<1$. In the following theorem we characterize 
bounded convolution operators $T_\lambda= k \star (\cdot)\!: L^2(w)\longrightarrow L^2(w)$, or equivalently multipliers $\lambda  \in {\rm Mult}(L^2(w))$, in terms of multipliers involving the weights $w_\alpha$ and $w_\beta$,  separately. 
Note that all multiplier algebras discussed below are embedded into $l^\infty(\mathbb{Z})$. The norms in the intersection and the sum of the multiplier spaces 
are introduced as usual for a Banach couple $(X^{(1)},X^{(2)})$ ($X^{(1)}, X^{(2)} \subset X$):
$$||x||_{X^{(1)}\cap X^{(2)}} = \max (||x||_{X^{(1)}}, ||x||_{X^{(2)}}), $$
$$ 
||x||_{X^{(1)}+X^{(2)}} = \inf \left \{ ||x^{(1)}||_{X^{(1)}}+||x^{(2)}||_{X^{(2)}}\!: \quad x=x^{(1)}+x^{(2)} \right \} , 
$$
where $x^{(1)}\in X^{(1)}, \, x^{(2)}\in X^{(2)}$. 

We will denote the class of bounded Fourier multipliers acting from $L^2(w_1)$ to $L^2(w_2)$  
by  ${\rm Mult}(L^2(w_1) \longrightarrow  L^2(w_2))$, and  in the case $w_1=w_2=w$ continue to use the notation ${\rm Mult}(L^2(w))$. 

It turns out that 
${\rm Mult}(L^2(w_{\beta }w^{\theta }_{\alpha }))$ can be characterized 
 as the intersection of  ${\rm Mult}(L^2(w_\alpha))$ and the sum of $ {\rm Mult}(L^2(w_\beta))$ and the ``rotated'' 
 multiplier  class  $\hat R_{\theta} ( {\rm Mult}(L^2(w_\beta)\longrightarrow L^2(w_\alpha)))$: 
 $$ {\rm Mult}(L^2(w_{\beta} w_\alpha^{\theta} ))= {\rm Mult}(L^2(w_\alpha)) \, \bigcap \left({\rm Mult}(L^2(w_\beta))+ \hat R_{\theta} ( {\rm Mult}(L^2(w_\beta)\longrightarrow L^2(w_\alpha)))\right).$$ 
 Here $(\hat R_{\theta} \lambda)_j = e^{-ij\theta} \lambda_j$ ($j \in \mathbb{Z}$) is a rotation operator on $\mathbb{Z}$. 
Note that the sum of the multiplier spaces above is not a direct sum since 
$${\rm Mult}(L^2(w_\beta))\bigcap  \hat R_{\theta} ( {\rm Mult}(L^2(w_\beta)\longrightarrow L^2(w_\alpha))) = {\rm Mult}(L^2(w_\beta) \longrightarrow  L^2).$$ 
 \ \par \noindent 
\bf 5.3. Theorem.\label{theorem1}  \it Let $\lambda =\mathcal{F} k \in l^\infty(\mathbb{Z})$.  Let 
$$w = w_{\beta} w_\alpha^{\theta} , \quad \theta \not= \pi n, \,\, \forall n \in \mathbb{Z},$$
where $0<\alpha < \beta<1$.  Then the following statements hold. \ \par \ \par
\noindent \bf (1)\it\, The inequality 
\begin{equation}\label{cond0} 
||k \star f||_{L^2(w)} \le C ||f||_{L^2(w)}
\end{equation}
holds for all $f \in L^2(w)$ if and only if  $k$ can be represented 
in the form 
\begin{equation}\label{rep} 
k=k^{(1)} + k^{(2)},
\end{equation}
where 
$k$, $k^{(1)}$, and $k^{(2)}$ are pseudo-measures satisfying the following conditions: 
\begin{equation}\label{cond1}
 || k \star f||_{L^2(w_\alpha)} \le C_0 ||f||_{L^2(w_\alpha)}, 
\end{equation}
\begin{equation}\label{cond2}
 || k^{(1)} \star f||_{L^2(w_\beta)} \le C_1 ||f||_{L^2(w_\beta)}, 
\end{equation}
\begin{equation}\label{cond3}
 ||(R_{-\theta} k^{(2)} )\star f||_{L^2(w_\alpha)} \le C_2 ||f||_{L^2(w_\beta)}, 
\end{equation}
where $(R_{-\theta} k^{(2)})(e^{it}) =k^{(2)}(e^{i(t+\theta)})$. \ \par \ \par
\noindent \bf (2)\it\, Condition (\ref{cond1})  in statement {\rm (1)} can be replaced with the following condition on $k^{(2)}$:
\begin{equation}\label{cond3a}
 || k^{(2)} \star f||_{L^2} \le C_3 ||f||_{L^2(w_\alpha)}. 
\end{equation}
 \ \par \ \par
\noindent \bf(3)\,\it Decomposition (\ref{rep})  
can be obtained explicitly as follows: 
\begin{equation}\label{rep1} 
k^{(1)} = \eta k, \quad k^{(2)}=(1-\eta) k,
\end{equation}
where $\eta$ is a cut-off function such that $\eta (e^{it}) =1$ 
if $|t|<a$,  and $\eta (e^{it}) =0$ outside $|t|<2a$, for some $0<a<\pi/4$ so that $4a<|\theta|<\pi$, under the additional assumption that  $\eta$ is in the Wiener 
algebra on $\mathbb{T}$, i.e., $\sum_{n \in \mathbb{Z}} |\hat\eta(n)|< +\infty$.  \ \par \ \par
\noindent \bf (4)\it\, In the case $\theta=\pi n$ $(n \in \mathbb{Z})$ inequality \eqref{cond0} holds 
if and only if $k=k^{(1)} + k^{(2)}$, so that \eqref{cond2} and \eqref{cond3} hold.
\rm 

 \ \par \noindent 
\bf 5.4. Remarks. \rm  \bf 1. \rm In the case $\theta=\pi n$ condition  \eqref{cond1} in statement (1) is  replaced with $ || k \star f||_{L^2(w_\alpha w_\alpha^\theta)} \le C_0 ||f||_{L^2(w_\alpha)}$, which is a consequence of \eqref{cond2} and \eqref{cond3}.

\bf 2. \rm Conditions (5.4) and (5.5) automatically imply that $ k^{(1)}$, $k^{(2)}$ 
are pseudo-measures. A direct characterization of the conditions on $k^{(1)}$, $k^{(2)}$ in representation (\ref{rep}) 
in terms of their Fourier coefficients is given below (see Corollary 5.7).  \ \par 
\bf 3. \rm It is easy to see that every function $ f\in L^{2}(w_{\beta 
}w^{\theta}_{\alpha })$ allows a decomposition $ f= f_{1}+f_{2}$ 
where $ f_{1}\in L^{2}(w_{\beta })$ and $ f_{2}\in L^{2}(w^{\theta }_{\alpha 
})$, unique up to a summand from the ``flat" space $ L^{2}(\mathbb{T})= 
L^{2}(w_{\beta })\bigcap L^{2}(w^{\theta}_{\alpha })$. We can write 
this decomposition as follows: 
$$ L^{2}(w_{\beta }w^{\theta}_{\alpha })= L^{2}(w_{\beta 
})+ L^{2}(w^{\theta}_{\alpha })+L^{2}(\mathbb{T}) .$$
Using this decomposition we can restate the principal claim of 
Theorem 5.3 as follows: the arrow on the left-hand side of the following diagram (consequently, 
a multiplier $ T_{\lambda }$) represents a bounded operator if and only if $ \lambda =  
\lambda^{(1)}+ \lambda^{(2)}$, $ T_{\lambda }=T_{\lambda^{(1)}}+ 
T_{\lambda^{(2)}}$, and all the arrows  on the right-hand side represent bounded operators (in the corresponding 
spaces):  \ \par
$$
\renewcommand{\arraystretch}{1.5}
\begin{array}{ccccc}
L^2(w_\beta w_\alpha^{\theta}) &=& \qquad L^2(w_\beta) \!  + \! L^2(w_\alpha^{\theta})\! + \! L^2(\mathbb{T})\\
\downarrow\rlap{$\scriptstyle T_\lambda$}&&\downarrow\rlap{$\scriptstyle T_{\lambda^{(1)}}$} \, \, \quad \searrow\rlap{$\scriptstyle T_{\lambda^{(2)}}$} \qquad \, \, \, \searrow\rlap{$\scriptstyle T_{\lambda^{(2)}}$}\\
L^2(w_\beta w_\alpha^{\theta}) &=& \qquad L^2(w_\beta) \!  + \! L^2(w_\alpha^{\theta}) \! + \! L^2(\mathbb{T})
\end{array}
$$
 \ \par \ \par
\bf Proof of Theorem 5.3. \rm We start with the following lemma which characterizes  the class ${\rm Mult}(L^2(w))$ 
in simpler terms.  \ \par \ \par \noindent 
\bf 5.5. Lemma. \it 
 (i) Under the assumptions of Theorem 5.3, inequality 
  (\ref {cond0}) holds, or equivalently $\lambda =\mathcal{F} k \in {\rm Mult}(L^2(w))$,   if and only if the following pair of inequalities hold: 
 \begin{equation}\label{ineq5}
||k \star f||_{L^2(w_{\beta} w_\alpha^{\theta}  )} \le C \, ||f||_{L^2(w_\beta)},
\end{equation}
\begin{equation}\label{ineq6}
||k \star f||_{L^2(w_\alpha)} \le C \, ||f||_{L^2(w_\alpha)}, 
\end{equation}
provided $\theta \not=\pi n$, $n \in \mathbb{Z}$. (ii) For $\theta=\pi n$, $n \in \mathbb{Z}$, the above statement 
holds if inequality \eqref{ineq6} is replaced with 
\begin{equation}\label{ineq6a}
||k \star f||_{L^2(w_\alpha w_\alpha^\theta)} \le C \, ||f||_{L^2(w_\alpha)}. 
\end{equation}
\ \par \ \par \noindent 
\bf 5.6. Remark \rm 
If $\alpha=\beta\ge 0$ then obviously (\ref{ineq5}) follows from (\ref{ineq6}). \ \par \ \par

\bf Proof of Lemma 5.5. \rm  
Let $\widetilde k(e^{it})= k(e^{-it})$. Notice that 
\begin{equation}\label{weights}
\frac 1 w = \frac{1}{w_{\beta} w_\alpha^{\theta} } \approx  w_{-\alpha}^{\theta} + w_{-\beta}.
\end{equation}
Then by duality, (\ref{cond0}) holds if and only if 
$$
||\widetilde k \star f||^2_{L^2(w_{-\beta})} + ||\widetilde k \star f||^2_{L^2(w_{-\alpha}^{\theta})}
$$
$$ \le C \, \left (||f||^2_{L^2(w_{-\beta})} + ||f||^2_{L^2(w_{-\alpha}^{\theta})} \right) , 
$$
for all trigonometric polynomials $f$, which are dense in  $L^2(w_{-\beta}) \bigcap L^2(w_{-\alpha}^{\theta})$. 
Consequently, (\ref{cond0}) is equivalent to the following pair of inequalities:
 $$
||\widetilde k \star f||^2_{L^2(w_{-\beta})} \le C \, \left (||f||^2_{L^2(w_{-\beta})} + ||f||^2_{L^2(w_{-\alpha}^{\theta})} 
\right) ,
$$
$$
||\widetilde k \star f||^2_{L^2(w_{-\alpha}^{\theta})} \le C \, \left ( ||f||^2_{L^2(w_{-\beta})} + 
||f||^2_{L^2(w_{-\alpha}^{\theta})} 
\right) . 
$$
Using duality and (\ref{weights}) again, we rewrite the preceding inequalities in the  equivalent  form: 
\begin{equation}\label{ineq1}
||k \star f||_{L^2(w_{\beta} w_{\alpha}^{\theta})} \le C \, ||f||_{L^2(w_{\beta})} ,
\end{equation}
 \begin{equation}\label{ineq2}
||k \star f||_{L^2(w_{\beta} w_{\alpha}^{\theta})} \le C \, ||f||_{L^2(w_{\alpha}^{\theta})} .  
\end{equation}
Notice that (\ref{ineq1}) coincides with (\ref{ineq5}). 
Applying the rotation operator $R_{\theta} f (e^{it}) = f(e^{i(t-\theta)})$, we see that (\ref{ineq2})  is equivalent to: 
\begin{equation}\label{ineq3}
||k \star f||_{L^2(w_{\alpha} w_{\beta}^{-\theta})} \le C \, ||f||_{L^2(w_{\alpha})} . 
\end{equation}
 Clearly, $w_{\beta} w_{\alpha}^{\theta}  + w_{\alpha} w_{\beta}^{-\theta} 
\approx w_\alpha$  in the generic case ($\theta \not=\pi n$) since $0<\alpha\le \beta$. Adding together  (\ref{ineq1}) and (\ref{ineq3}) we arrive at the inequality 
$$
||k \star f||_{L^2(w_{\alpha})} \le C \, ||f||_{L^2(w_{\alpha})},   
$$
which coincides with (\ref{ineq6}). Moreover, the preceding inequality is stronger than (\ref{ineq3}) 
since  $w_{\alpha} w_{\beta}^{-\theta} \le 2^\beta w_{\alpha}$. 
Thus, (\ref {cond0}) holds 
if and only if 
both (\ref{ineq5}) and (\ref{ineq6}) hold. 

If $\theta =\pi n$ ($n \in \mathbb{Z}$), then $w_{\beta} w_{\alpha}^{\theta}  + w_{\alpha} w_{\beta}^{-\theta} 
\approx w_\alpha w_\alpha^\theta$, and  $w_{\alpha} w_{\beta}^{-\theta} \le w_{\alpha}w_{\alpha}^{\theta}$. Hence \eqref{ineq3}, and consequently \eqref{cond0}, holds 
if and only if 
both (\ref{ineq5}) and (\ref{ineq6a}) hold. 
\qed \ \par \ \par

We now prove the sufficiency part of Theorem 5.3. 
Suppose that 
$k= k^{(1)} + k^{(2)},$
where  $k^{(1)}$ and $k^{(2)}$ satisfy  (\ref{cond2}) and  (\ref{cond3}) respectively. (The corresponding 
decomposition of the multiplier $\lambda$ will be written in the form $\lambda = \lambda^{(1)} + \lambda^{(2)}$ where $\lambda^{(i)}=\mathcal{F} k^{(i)}$.)  
Notice that   (\ref{cond3}) is equivalent to: 
\begin{equation}\label{ineq3a}
||k^{(2)} \star f||_{L^2(w_{\alpha}^{\theta})} 
\le C_2 ||f||_{L^2(w_{\beta})}.
\end{equation}
 Then clearly, 
$$
||k \star f||_{L^2(w_{\beta} w_{\alpha}^{\theta})} \le 2^\alpha  ||k^{(1)} \star f||_{L^2(w_{\beta})} 
+ 2^\beta || k^{(2)} \star f||_{L^2(w_{\alpha}^{\theta})}
$$
$$
\le (2^\alpha C_1+2^\beta C_2) 
||f||_{L^2(w_{\beta})}.  
$$
This proves  (\ref{ineq5}). If $\theta \not=\pi n$ then 
 combining  (\ref{ineq5}) with (\ref{ineq6}) yields  
(\ref{cond0}) by Lemma 5.5. Similarly, for $\theta= \pi n$, we use  (\ref{ineq5}) together with (\ref{ineq6a})  to see that (\ref{cond0}) holds.

Suppose now that $\theta \not=\pi n$, and condition (\ref{cond3a}) is used in place of (\ref{ineq6}), 
that is, we assume 
$$
\lambda^{(2)} \in  {\rm Mult}(L^2(w_\alpha) \longrightarrow  L^2).
$$
The preceding condition obviously implies 
\begin{equation}\label{cond3aa}
\lambda^{(2)} \in  {\rm Mult}(L^2(w_\alpha))).
\end{equation}
Since  $\lambda^{(1)} \in {\rm Mult}(L^2(w_\beta))$, it follows  that 
$\lambda^{(1)} \in   {\rm Mult}(L^2)=l^\infty(\mathbb{Z})$. Hence 
by interpolation, 
\begin{equation}\label{cond3ab}
\lambda^{(1)} \in {\rm Mult}(L^2(w_\alpha))). 
\end{equation}
Combining (\ref{cond3aa}) and (\ref{cond3ab}) we obtain (\ref{ineq6}). 
This completes the proof of the sufficiency part of Theorem 5.3. \ \par \ \par 

To prove the necessity part,  suppose  (\ref{cond0}) holds, or equivalently, 
$\lambda  \in {\rm Mult}(L^2(w))$.  Then (\ref{ineq5}) holds as well by Lemma 3.5. 
Let $ \eta $ be a cut-off function defined in Theorem 
5.3 (3). Consider decomposition 
(\ref{rep1}), where $k^{(1)}=\eta k$ and $k^{(2)}= (1-\eta)k$. Then 
 clearly, $\lambda^{(1)} \in  {\rm Mult}(L^2(w))$, since 
 $$
 k^{(1)}\star f(e^{it}) = 
  \sum_{j \in \mathbb{Z}}  \hat \eta(j)  e^{ijt} (k\star f_j)(e^{it}),
  $$ 
where $f_j(e^{it})=  e^{-ijt}  f(e^{it})$. Hence,
\begin{equation}\label{ineq7}
 ||k^{(1)}\star f||_{L^2(w)} \le    \sum_{j \in \mathbb{Z}}  |\hat \eta(j)| \, 
 ||k\star f_j||_{L^2(w)} \le C \,  ||f||_{L^2(w)} \,  \sum_{j \in \mathbb{Z}}  
 |\hat \eta(j)|. 
\end{equation}
Since $\lambda  \in {\rm Mult}(L^2(w))$, it follows from the 
preceding inequality that 
\begin{equation}\label{ineq7a}
 ||k^{(2)}\star f||_{L^2(w)}  \le C \,  ||f||_{L^2(w)}
\end{equation}
as well, or equivalently $\lambda^{(2)} \in {\rm Mult}(L^2(w))$.  

Next, we estimate 
$$
 ||k^{(1)}\star f||^2_{L^2(w_\beta)} 
 = \int_{|t|<3a} | k^{(1)}\star f|^2 w_\beta  dt + 
 \int_{|t|\ge 3a} | k^{(1)}\star f|^2 w_\beta dt = I + II. 
 $$
Since $w(e^{it}) \approx w_\beta(e^{it})$ for $|t|<3a$, we deduce from 
(\ref{ineq5}) and (\ref{ineq7}):
$$ 
I \le || k^{(1)}\star f||^2_{L^2(w)} \le C ||k_1\star f||^2_{L^2(w)} 
\le C_1 ||f|| ^2_{L^2(w_\beta)}.
$$ 
To estimate the second term, notice that $k^{(1)}(e^{i(t-\tau)})$ is supported in 
$|t - \tau|\le 2a$. Hence, for $|t|\ge 3a$, we have $|\tau|\ge a$, so that we may 
assume that $f(e^{i\tau})$ is supported in 
$|\tau|\ge a$. In other words, we may replace $f$ in $II$ with 
$f \chi_{|\tau|\ge a}$. 
Moreover, for $|t|\ge 3a$, $w_\beta(e^{it}) \approx 1$. 
It follows, 
$$
II \le || k^{(1)}\star ( f \chi_{|\tau|\ge a})||^2_{L^2} \le 
||\lambda^{(1)}||^2_{l^\infty} ||f \chi_{|\tau|\ge a}||^2_{L^2} 
\le C ||f||^2_{L^2(w_\beta)}.
$$
Combining the preceding estimates, we see that 
$\lambda^{(1)} \in  {\rm Mult}(L^2(w_\beta))$. 

	Similarly, for $k^{(2)}= (1-\eta)k$ and 
	$(R_{-\theta} k^{(2)}) (e^{i \tau})=(1-\eta(e^{i (\tau+\theta)})
k (e^{i (\tau+\theta)})$, we obtain:  
$$ 
|| (R_{-\theta} k^{(2)}) \star f||^2_{L^2(w_{\alpha})} 
= || k^{(2)}\star f||^2_{L^2(w_{\alpha}^{\theta})}
$$
$$ 
 = \int_{|t|<\frac{a}{2}} |k^{(2)}\star f|^2    w_{\alpha}^{\theta} dt + 
 \int_{|t|\ge\frac{a}{2}} |k^{(2)}\star f|^2  w_{\alpha}^{\theta} dt
 $$
 $$ = III+ IV. 
$$
We first estimate $III$. Since $k^{(2)}(e^{i(t-\tau)})$ is supported in $|t-\tau|\ge a$ and $|t|<\frac{a}{2}$, it follows that $f(e^{i \tau})$ can be replaced  with 
$f \chi_{|\tau|>\frac{a}{2}}$. Hence, 
$$
III \le || k^{(2)} \star f \chi_{|\tau|>\frac{a}{2}}||^2_{L^2(w_{\alpha}^{\theta})} 
\le || k^{(2)} \star f \chi_{|\tau|>\frac{a}{2}}||^2_{L^2} 
$$
$$\le ||k^{(2)}||_{l^\infty} || f \chi_{|\tau|>\frac{a}{2}}||^2_{L^2} 
\le C  ||k^{(2)}||_{l^\infty} || f ||^2_{L^2(w_\beta)}.
$$
To estimate $IV$, we notice that, for $|t|\ge\frac{a}{2}$, 
$$
w_{\alpha}^{\theta}(e^{it})  \approx w_{\alpha}^{\theta}(e^{it})   w_{\beta} (e^{it})  
= w(e^{it}) .
$$ 
Consequently, using (\ref{ineq7a}), we estimate: 
$$
IV\le C_1 || k^{(2)} \star f||^2_{L^2(w)} \le C_2 || k \star f||^2_{L^2(w)} \le 
C_3 ||f|| |^2_{L^2(w)}. 
$$
Combining these estimates, we deduce 
$$\lambda^{(2)} \in {\rm Mult}(L^2(w_\beta) \longrightarrow  L^2(w_{\alpha}^{\theta})),
$$ 
which is  equivalent to 
 (\ref{cond3}): 
 \begin{equation}\label{ineq2a}
  || k^{(2)} \star f||_{L^2(w_{\alpha}^{\theta})}\le C 
  ||f||_{L^2(w_\beta)}.
 \end{equation}
 This completes the proof in the case $\theta= \pi n$. 
 
 The above estimates remain true if $\theta\not= \pi n$. Additionally,  
 it follows from Lemma 5.2 that (\ref{ineq6}) holds. 
Using (\ref{ineq7}), (\ref{ineq7a}) with $w_\alpha$ in place of $w$  
 we deduce that (\ref{ineq6}) holds 
 with $k^{(2)}$ in place of $k$: 
 \begin{equation}\label{ineq2b}
  || k^{(2)} \star f||_{L^2(w_\alpha)}\le C 
  ||f||_{L^2(w_\alpha)}.
 \end{equation} 
 Adding together (\ref{ineq2a}) and (\ref{ineq2b}) we obtain:
 $$
 || k^{(2)} \star f||_{L^2}\le C 
  ||f||_{L^2(w_\alpha)}.
 $$
 This proves (\ref{cond3a}), and completes the proof of the necessity 
 part of Theorem 5.3 in the generic case $\theta\not=\pi n$.  \qed \ \par \ \par
 Combining Theorem 5.3 with Theorem 4.11 we obtain the following characterization of multipliers.  
 \ \par \ \par \noindent 
\bf 5.7. Corollary. \it Under the assumptions of Theorem 5.3, 
 inequality (\ref{cond0}) holds if and only if $k=k^{(1)}+k^{(2)}$, where 
 $\lambda^{(i)} = \mathcal{F} k^{(i)} \in l^\infty(\mathbb{Z})$ $(i=1, 2)$,  and for every 
finite $J \subset \mathbb{Z}$, the following conditions hold: 
\begin{equation}\label{cond4}
 \sum_{j \in J} \sum_{m \in \mathbb{Z}} \frac{ \vert  \lambda_j^{(1)} - \lambda_m^{(1)} \vert^2}
 {(|j-m|+1)^{1+ \beta}}   \le C \, \rm{Cap}_\beta (J), 
\end{equation}
\begin{equation}\label{cond5}
 \sum_{j \in J} \sum_{m\in \mathbb{Z}} \frac{ \vert  e^{i j \theta}  \lambda_j ^{(2)} - 
 e^{i m \theta} \lambda_m^{(2)}  \vert^2}
 {(|j-m|+1)^{1+ \alpha}}   \le C \, \rm{Cap}_\beta (J), 
\end{equation}
\begin{equation}\label{cond6}
 \sum_{j \in J}  \vert \lambda_j^{(2)}  \vert^2 \le C \, \rm{Cap}_\gamma (J), 
\end{equation}
where $\gamma = \max (\alpha, \beta-\alpha)$ and $C$ does not depend on $J$, provided 
$\theta \not=\pi n$, $n \in \mathbb{Z}$. 

In the case $\theta=\pi n$,  (\ref{cond0}) holds if and only if $k=k^{(1)}+k^{(2)}$ so that 
 $\lambda^{(i)}\in l^\infty(\mathbb{Z})$ $(i=1, 2)$, and \eqref{cond4}, \eqref{cond5} hold. 
\rm \ \par \ \par 

\noindent 
\bf Remark. \rm In condition \eqref{cond6} it suffices to let 
$\gamma=\alpha$. However, the proof of this assertion is complicated, and we do not present it here. 
(It requires discrete analogues of multiplier estimates obtained earlier by the second author in 
the continuous case; see Sec. 3.2.10 in \cite{MSh2009}.) 
\  \par \ \par

The following example demonstrates that conditions (\ref{cond4})--(\ref{cond6}) 
are essential, and in a sense cannot be relaxed. Moreover, the inequality 
\begin{equation}\label{cond7}
 || k \star f||_{L^2(w_\beta)} \le C ||f||_{L^2(w_\beta)}, 
\end{equation}
is only sufficient, but not necessary for (\ref{cond0}) 
in the case $\alpha< \beta$, contrary to the case $\alpha=\beta$. In other words, we cannot 
let $k^{(2)}=0$ in decomposition (\ref{rep}). On the other hand, condition 
(\ref{cond1}) is only necessary, but not sufficient for (\ref{cond0}). 
\ \par \ \par \noindent 
\bf 5.8. Example. \rm Suppose $w=w_{\beta} w_\alpha^{\theta}$ where $0<\alpha<\beta<1$ 
and $\theta \not=\pi n$ ($n \in \mathbb{Z}$). For $\gamma = \max (\alpha, \beta-\alpha)$, we pick $\delta>0$ so that $\frac{\gamma}{2}< \delta < \frac {\beta}{2}$. 
Let 
$\lambda=(\lambda_j)$, where $\lambda_j = \frac{e^{-ij\theta}}{(|j|+1)^\delta}$ ($j \in \mathbb{Z}$). Then conditions (\ref{cond5})--(\ref{cond6}) 
hold for $k=k^{(2)}$, $k^{(1)}=0$, but (\ref{cond4}) fails for 
$k=k^{(1)}$, $k^{(2)}=0$. In other words,  $\lambda  \in {\rm Mult}(L^2(w))$, but $\lambda  \not\in {\rm Mult}(L^2(w_\beta))$. Moreover, $\overline{\lambda}=(\overline{\lambda}_j) \not\in {\rm Mult}(L^2(w))$. 
If 
$\delta <\frac{\gamma}{2}$, then 
$\lambda  \not\in {\rm Mult}(L^2(w))$. 
\ \par \ \par

The claims in Example 5.8 follow from the well-known fact 
that if $\Lambda=(\Lambda_j)$, where 
$$\Lambda_j= \frac{1}{(|j|+1)^\delta}, \quad j 
\in \mathbb{Z}, \quad 0<\delta<1,$$  
then $\Lambda \in {\rm Mult} (L^2(w_\beta) \longrightarrow  L^2(w_\alpha))$ for $0\le \alpha<\beta<1$ if and only if 
$0<\delta \le \frac{\beta-\alpha}{2}$. Note that in this example 
$\lambda  \in {\rm Mult}(L^2(w))$, but $\overline{\lambda}  \not\in {\rm Mult}(L^2(w))$, since otherwise 
by interpolation the Fourier multiplier $T_\Lambda\!: \, L^2(w_\beta) \longrightarrow  L^2(w_\alpha)$ would be bounded, which fails for 
$\frac{\beta-\alpha}{2}< \delta$. 
\ \par \ \par

 The next theorem shows that, for weights with a finite number of power-like singularities, 
 ${\rm Mult}(L^2(w))$ 
  has the Spectral Localization Property ({\rm SLP}). As above, for the sake of simplicity, 
 we consider weights with two generic zeros, 
 $$w = w_{\alpha}^{\theta} w_\beta, \quad \theta \not= \pi n, \,\,  \forall n \in \mathbb{Z},$$
 where $0<\alpha\le \beta<1$. (It is easy to see that in the case   $\theta = \pi$  the {\rm SLP} holds as well.) 
 \ \par \ \par \noindent 
\bf 5.9. Theorem. \it Under the hypotheses of Theorem 5.3, suppose $\lambda  \in {\rm Mult}(L^2(w))$, so that $\lambda =\mathcal{F} k$, 
  where $k$ is a pseudo-measure on $\mathbb{T}$ 
 such that the inequality 
 \begin{equation}\label{cond1a}
|| k \star f||_{L^2(w)} \le C ||f||_{L^2(w)}
\end{equation}
holds for all $f \in L^2(w)$. If $\, \inf_{j \in \mathbb{Z}} |\lambda_j| >0$, 
then $1/\lambda  \in {\rm Mult}(L^2(w))$. \rm 
\ \par \ \par 

\bf Proof. \rm  Suppose $\lambda  \in {\rm Mult}(L^2(w))$ and 
 $\inf_{j \in \mathbb{Z}} |\lambda_j| =\delta >0$. By Theorem 5.3, 
 $\lambda =\lambda^{(1)}+\lambda^{(2)}$, where the following three conditions hold: 
\begin{equation}\label{cond0b}
\lambda \in  {\rm Mult}(L^2(w_\alpha)), 
\end{equation}
 \begin{equation}\label{cond1b}
\lambda^{(1)}\in  {\rm Mult}(L^2(w_\beta)), 
\end{equation}
 \begin{equation}\label{cond2b}
\lambda^{(2)}  \in   {\rm Mult}(L^2(w_\beta)\longrightarrow  L^2(w_{\alpha}^{\theta}))\bigcap {\rm Mult}(L^2(w_\alpha)\longrightarrow  L^2).
\end{equation}
We will also need the following relations which follow from (\ref{cond0b}) and (\ref{cond2b}) respectively by applying the rotation 
operator $\hat R_{\theta}$:
\begin{equation}\label{cond2c}
\lambda   \in   {\rm Mult}(L^2(w_{\alpha}^{\theta})), 
\end{equation}
\begin{equation}\label{cond2d}
 \lambda^{(2)} \in  {\rm Mult}(L^2(w_{\alpha}^{\theta})\longrightarrow  L^2).
\end{equation}

We first prove Theorem 5.9 under the additional assumption 
  \begin{equation}\label{assump}
 \, \inf_{j \in \mathbb{Z}} |\lambda_j^{(1)}| >0.
 \end{equation}
 Then 
 \begin{equation}\label{rep2}
\frac 1 {\lambda } = \frac {1} {\lambda^{(1)}} + \frac{- \lambda^{(2)}}{\lambda  \, 
 \lambda^{(1)}}. 
 \end{equation}
  Using 
   assumption (\ref{assump}) and the {\rm SLP} for the weight $w_\beta$,  
we deduce from (\ref {cond1b}): 
\begin{equation}\label{cond2e}
\frac 1 {\lambda^{(1)}} \in {\rm Mult}(L^2(w_\beta)). 
\end{equation} 
 Since $1/\lambda^{(1)}\in {\rm Mult}(L^2)=l^\infty$, 
and $\alpha\le \beta$,  using interpolation we see that 
$1/ {\lambda^{(1)}} \in {\rm Mult}(L^2(w_\alpha))$. Applying 
 the rotation operator $R_{\theta}$, 
we obtain  
\begin{equation}\label{cond2f}
\frac 1 {\lambda^{(1)}} \in {\rm Mult}(L^2(w_{\alpha}^{\theta})). 
\end{equation}
By (\ref{cond2c}) and the {\rm SLP} for the weight $w_{\alpha}^{\theta}$, 
\begin{equation}\label{cond2g}
\frac 1 {\lambda} \in  {\rm Mult}(L^2(w_{\alpha}^{\theta})).
\end{equation}
By (\ref{cond2b}), $\lambda^{(2)} \in  {\rm Mult}(L^2(w_\beta)\longrightarrow  L^2(w_{\alpha}^{\theta}))$. 
 Consequently, as a product of three multipliers, 
$$
 \frac{\lambda^{(2)}}{\lambda  \, 
 \lambda^{(1)}} \in  {\rm Mult}(L^2(w_\beta)\longrightarrow  L^2(w_{\alpha}^{\theta})) .$$
 From (\ref{cond2d}) and (\ref{cond2f}), (\ref{cond2g}) it follows that 
$$
 \frac{\lambda^{(2)}}{\lambda  \, 
 \lambda^{(1)}} \in    {\rm Mult}(L^2(w_{\alpha}^{\theta})\longrightarrow  L^2)$$ 
 as well. This proves that decomposition (\ref{rep2}) for $1/\lambda $ is of the same type as for $\lambda $ 
 in Theorem 5.3. Thus, by the sufficiency part of  Theorem 5.9, we conclude that    
 $1/\lambda   \in {\rm Mult}(L^2(w))$. 
 
 We now demonstrate how to remove the additional assumption (\ref{assump}) used above. Suppose 
   \begin{equation}\label{lower}
 \delta =  \inf_{j \in \mathbb{Z}} |\lambda_j| >0.
 \end{equation}
 Let $k=k^{(1)}+k^{(2)}$ be decomposition (\ref{rep}) obtained in   Theorem 5.3. 
Let 
$$Z_1 = \left \{ j \in \mathbb{Z}: \, \, | \lambda_j^{(1)}| \ge \frac \delta 2 \right \}, \quad 
Z_2 = \mathbb{Z}\setminus Z_1.
$$
Then, obviously,  
  \begin{equation}\label{lower1}
\inf_{j \in Z_2} |\lambda_j^{(2)}| \ge \frac \delta 2. 
\end{equation}
We claim that the following inequality holds for every $f \in L^2(w_\beta)$: 
\begin{equation}\label{claim}
\sum_{j \in Z_2} |\hat f(j)|^2 \le C \, ||f||^2_{L^2(w_\beta)}. 
\end{equation}
In other words, the set $Z_2$ is quite meager. 
Indeed, since 
$$\lambda^{(2)} \in  {\rm Mult}(L^2(w_\beta)\longrightarrow  L^2(w_\alpha))\bigcap {\rm Mult}(L^2(w_{\alpha}^{\theta})\longrightarrow  L^2) ,$$
we obtain for every  $f \in L^2(w_\beta)$: 
$$
\sum_{j \in \mathbb{Z}} |\lambda_j^{(2)}|^4 |\hat f(j)|^2 = 
||k^{(2)}\star(k^{(2)}\star f)||^2_{L^2}
$$
$$
\le C \, 
||k^{(2)}\star f||^2_{L^2(w_{\alpha}^{\theta})} \le C  ||f||^2_{L^2(w_\beta)}. 
$$
On the other hand, by (\ref{lower1}), 
$$
\sum_{j \in Z_2} |\hat f(j)|^2 \le \frac{16}{\delta^4} 
\sum_{j \in Z_2} |\lambda_j^{(2)}|^4 |\hat f(j)|^2
$$
$$ 
\le \frac{16}{\delta^4}   \sum_{j \in \mathbb{Z}} |\lambda_j^{(2)}|^4 |\hat f(j)|^2 . 
$$
Combining the preceding estimates, we prove (\ref{claim}). 

From  (\ref{claim}) we deduce:
$$
\sum_{j \in Z_2} |\lambda_j^{(2)}|^2 |\hat f(j)|^2 \le || \lambda^{(2)} ||^2_{l^\infty} 
\sum_{j \in Z_2} |\hat f(j)|^2 \le C   ||f||^2_{L^2(w_\beta)} . 
$$
This yields:
  \begin{equation}\label{intersection}
 \lambda^{(2)}  \chi_{Z_2}  
 \in {\rm Mult}(L^2(w_\beta)\longrightarrow  L^2).
 \end{equation}
 We can now adjust decomposition  (\ref{rep}): 
 $k= k^{(3)}  + k^{(4)}$, where 
 $$\lambda^{(3)} = \lambda^{(1)}  + \lambda^{(2)}  \chi_{Z_2}, \quad 
 \lambda^{(4)} = \lambda^{(2)} -    \lambda^{(2)}  \chi_{Z_2} .$$
 Clearly, by (\ref{lower}) and  (\ref{lower1}), 
   \begin{equation}\label{lower2} 
   \inf_{j \in \mathbb{Z}} |\lambda_j^{(3)}| \ge \frac \delta 2 .
 \end{equation}
 Moreover, (\ref{intersection}) yields that the required conditions 
 still hold for the new components $\lambda^{(3)}$ and $\lambda^{(4)}$: 
 $$\lambda^{(3)} \in {\rm Mult}(L^2(w_\beta)) ,$$
 $$ \lambda^{(4)}  \in {\rm Mult}(L^2(w_\beta)\longrightarrow  L^2(w_{\alpha}^{\theta})), 
 $$
 $$
 \lambda^{(4)} \in {\rm Mult}(L^2(w_\alpha)\longrightarrow  L^2).
 $$  
 Thus,  assumption (\ref{assump}) is redundant in the general case. This proves that $1/\lambda   \in {\rm Mult}(L^2(w))$. 
\qed 

\ \par \ \par \noindent 
\bf 5.10. Remark. \rm The same argument as in Theorem 5.3 works in the non-generic case 
$\theta=\pi$ as well. For instance, for $\alpha=\beta$, i.e., $w=w_\beta w_\beta^{\pi}$,  
or $w=w_{-\beta} w_{-\beta}^{\pi}$ ($0<\beta<1$), 
it follows that $\lambda \in  {\rm Mult}(L^2(w))$ if and only if 
$$\lambda  = \lambda^{(1)} + \hat R_\pi \lambda^{(2)}, $$
where $ \lambda^{(1)}, \lambda^{(2)} \in  {\rm Mult}(L^2(w_\beta)),$ and $
\hat R_\pi \lambda = \big ((-1)^j \lambda_j\big)$ 
($j \in \mathbb{Z}$) is the corresponding rotation operator. An equivalent  direct characterization 
of ${\rm Mult}(L^2(w))$ is given 
in Sec. 5.1. 
\ \par \ \par \noindent 
\bf 5.11. Remark. \rm  More generally, suppose $w = w_{-\beta} \star \nu$ ($0<\beta<1$), where 
$\nu = \sum_{j=0}^{d-1} c_j \, \delta_{\zeta^j}$ is a finite discrete measure 
($c_j > 0$), and $\{\zeta^j\}$ are the roots of unity of order $d$. An explicit  description 
of multipliers ${\rm Mult}(L^2(w))$ is given in Sec. 5.1. There is an alternative 
``sliced'' decomposition:  $\lambda  \in   {\rm Mult}(L^2(w))$ if and only if 
$$\lambda  = \sum_{j= 0}^{d-1}   \hat R_{\zeta^j} \lambda^{(j)},$$ 
where $\hat R_\zeta \lambda=(\bar{\zeta^j} \lambda_j)$ is the corresponding rotation operator, and 
each $\lambda^{(j)} \in    {\rm Mult}(L^2(w_{\beta}))$, $j=1,2, \ldots, d-1$. 
This can be proved by means of a decomposition similar to that used in the proof 
of Theorem 5.3 with smooth cut-off functions $\eta_j$, 
$0\le j\le  d-1$. 
\ \par 
\ \par 
The case of infinitely many singularities, which is 
briefly discussed below, is quite different.\ \par 
 \ \par\noindent 
\bf 5.12. Infinite superposition of {\rm LKS} singularities: 
the hidden spectrum. \rm Continuing Remark 5.11, one can consider the 
case of an infinite combination of {\rm LKS} singularities (in a dual form), 
say (for the sake of simplicity) of the same order, as follows: $ w=  
w_{-\beta }\star \nu $, where $ 0<\beta <1$, and 
$$ 
\nu = \sum _{k\geq 0}c_{k}\delta _{\zeta ^{k}}, \quad 
\text{where} \, \,  c_{k}>0, \, \,  \sum_{k\geq 0} c_{k}<\infty, \, \,  
  \sup_{k\geq 0} \, {\frac{c_{k}}{c_{k+1}}} 
< \infty, 
$$
 with $ \zeta \in \mathbb{T}$ such that $ \zeta ^{k}\not= 1$ ($ \forall 
k\in \mathbb{Z}$). In this case, at the moment, we can only conjecture (but not prove) 
that there is an analogue  
of a ``sliced" decomposition from Remark 5.11 for a multiplier $ T_{\lambda }\in 
{\rm Mult}(L^{2}(w_{-\beta }\star\nu ))$. We recall, however, that  the 
situation can be more complicated: in \cite{Nik2009} it is shown that in 
this case there exists a ``hidden spectrum," i.e. the {\rm SLP} does not hold. 

In fact, we believe that the reason why the ``hidden spectrum" appears 
lies in a kind of holomorphic extension of multipliers $ n\longmapsto \lambda_{n}$ ($n \in \mathbb{Z}$) of the space $ L^{2}(\nu )$, followed with a  ``sliced 
decomposition" mentioned above. The latter property is still a conjecture, 
but the former one (namely, the holomorphic nature of ${\rm Mult}(L^{2}(\nu ))$) 
is confirmed by the following claim, for which we need a bit of notation. 
In order to distinguish Fourier multipliers of $ L^{2}(\nu )$ from 
yet another (pointwise) holomorphic multipliers appearing in the next 
theorem, we temporarily change the notation for the former adding a 
subscript $ ``F"$ (for Fourier): 
$$
{\rm Mult}_{F}(L^{2}(\mathbb{T},\nu ))\!:=  {\rm Mult}(L^{2}(\mathbb{T},\nu 
)) .$$
For $c=(c_k)_{k\ge 0}$, let $1/c=(1/c_k)_{k\ge 0}$, and denote by 
$ l_{a}^{2}(1/c)$  the Hilbert space of 
functions $ f$ holomorphic on the unit disc $ \mathbb{D}$ such that 
$$ 
\Big \Vert f\Big \Vert ^{2}=   \sum _{k\geq 
0}\Big \vert \hat f(k)\Big \vert ^{2}{\frac{1}{c_{k}}} 
<  \infty .  
$$ 
Let $ {\rm Mult}(l_{a}^{2}(1/c))$ stand for (standard, pointwise) 
multipliers of $ l_{a}^{2}(1/c)$: 
$$ 
\varphi \in {\rm Mult}(l_{a}^{2}(1/c)) \Leftrightarrow   
\{ f\in l_{a}^{2}(1/c)  \Rightarrow   \varphi f\in l_{a}^{2}(1/c)\} .$$ 
In particular, if $ l_{a}^{2}(1/c)$ is an algebra (for instance, when $ c_{k}=1/(k+1)^{1+\epsilon }$, 
$ \epsilon >0$), then 
$$
 {\rm Mult}(l^{2}_{a}(1/c))= {\rm mult}(l^{2}_{a}(1/c))= 
l^{2}_{a}(1/c), $$
 where ${\rm mult}(l^{2}_{a}(1/c))$ stands for the closure of polynomials 
in the norm $ \Vert \cdot \Vert_{{\rm Mult}(l^{2}_{a}  (1/c))}$. 
Note that we always have $ {\rm Mult}(l_{a}^{2}(1/c))\subset $ $ l_{a}^{2}(1/c)\subset 
  l^{1}_{a}$ (the Wiener algebra).
     \ \par 
\ \par \noindent 
\bf 5.13. Theorem. \it Under the above assumptions on $\nu$, $c=(c_k)$ and $\zeta$,  
$$
 {\rm Mult}_{F}(L^{2}(\mathbb{T},\nu ))=  \{\lambda_{n}=  
\varphi (\zeta ^{n})  (\forall n\in \mathbb{Z})\!: \, \,  \varphi \in 
{\rm Mult}(l_{a}^{2}(1/c))\} . $$
Moreover, the ``visible spectrum" of a multiplier $ \lambda =  (\lambda 
_{n})_{n\in \mathbb{Z}}=$ $ (\varphi (\zeta ^{n}))_{n\in \mathbb{Z}}$ 
is a continuous curve 
$$ \rm{clos} \{\varphi (\zeta ^{n})\!:  \, \, n\in \mathbb{Z}\}=  
\varphi (\mathbb{T}) ,
$$ 
but the entire spectrum is the $ \varphi $-image of the closed disc 
$ \overline{\mathbb{D}}$:  
$$ 
\sigma (T_{\lambda })=  \varphi (\overline{\mathbb{D}}) , 
$$ 
and, at least for $ \varphi \in {\rm mult}(l^{2}_{a}(1/c))$, every point 
$ z\in \varphi (\mathbb{D})\!\setminus\! \varphi (\mathbb{T})$ is a Fredholm point 
of $ T_{\lambda }$ so that  
$$
 \rm{ind}(T_{\lambda }-zI)=  \rm{dim \, Ker}(T_{\lambda }-zI)=  
\rm{wind}(\varphi -z) . $$ 
  
    \bf Proof. \rm Notice that 
    $$ 
    L^{2}(\mathbb{T},\nu )=  \{(f(\zeta ^{k}))_{k\geq 
0}\!:  \int _{\mathbb{T}}\vert f\vert ^{2}d\nu =  \sum _{k\geq 
0}\vert f(\zeta ^{k})\vert ^{2}c_{k}<  \infty \}=  l_{a}^{2}(c) .$$ 
(We use a natural identification, $ (a_{k})=  (f(\zeta ^{k}))\longmapsto   
\sum _{k\geq 0}a_{k}z^{k}$.) By the hypothesis the backward shift 
$ S^{*}(f(\zeta ^{k}))_{k\geq 0}=  (f(\zeta ^{k+1}))_{k\geq 0}$ 
is a bounded operator on $ L^{2}(\mathbb{T},\nu )=  l_{a}^{2}(c)$. 
But $ S^{*}=  T_{\{\zeta ^{n}  \}  }$ is in $ {\rm Mult}_{F}(L^{2}(\mathbb{T},\nu 
))$, since   
$$ 
S^{*}z^{n}=  S^{*}((\zeta ^{k})^{n})_{k\geq 0}=  
((\zeta ^{k+1})^{n})_{k\geq 0}=  \zeta ^{n}((\zeta ^{k})^{n})_{k\geq 
0}=  \zeta ^{n}z^{n}, \quad  \forall n\in \mathbb{Z} . 
$$
Consequently, any multiplier operator $ T_{\lambda }\in {\rm Mult}_{F}(L^{2}(\mathbb{T},\nu 
))$ commutes with $ S^{*}$, and hence $ T_{\lambda }^{*}$ commutes 
with the shift $ S$ on the dual space $ l_{a}^{2}(1/c)$. So, $ T_{\lambda 
}^{*}=  \varphi (S)\in {\rm Mult}(l_{a}^{2}(1/c))$. For every $ 
F\in l_{a}^{2}(1/c)$, we have (under the bilinear duality $ \langle F,G \rangle =  
\sum _{k\geq 0}\hat F(k)\hat G(k)$): 
 $$ \langle T_{\lambda }z^{n},F \rangle = \langle ((\zeta ^{k})^{n})_{k\geq 
0},T_{\lambda }^{*}F \rangle=  \langle ((\zeta ^{k})^{n})_{k\geq 0},\varphi F \rangle =  
\varphi (\zeta ^{n})F(\zeta ^{n}) $$ 
$$ =  \varphi (\zeta ^{n}) \langle ((\zeta ^{k})^{n})_{k\geq 
0},F \rangle =  \varphi (\zeta ^{n}) \langle z^{n}, F \rangle  .  
$$ 
Hence $ T_{\lambda }z^{n}=  \varphi (\zeta ^{n})z^{n}$ 
for every $ n\in \mathbb{Z}$, i.e., $ T_{\lambda }=  
\varphi (S)^{*}$.\ \par 
    Clearly, the converse is also true, i.e.,  $ \varphi (\zeta ^{n})_{n\in 
\mathbb{Z}}$ is in $ {\rm Mult}_{F}(L^{2}(\mathbb{T},\nu )$ for every $ \varphi \in 
{\rm Mult}(l_{a}^{2}(1/c))$.

    The spectral nature of the adjoint operator $ T_{\lambda }=  
\varphi (S)^{*}$ related to a multiplier $ \varphi \in {\rm Mult}(l_{a}^{2}(1/c))$ 
is well known (see \cite{Nik1986}).  \qed \ \par 
\ \par




\begin{thebibliography}{GMcG1979}
  
  \bibitem[AH1996]{AH1996}   D.~R.~Adams and L.~I.~Hedberg, 
  \it Function Spaces and Potential Theory\rm, 
  Grundlehren der math. Wissenschaften, 314,  Berlin--Heidelberg, Springer, 1996. 
  \ \par 
\ \par 
\bibitem[BeD1958]{BeD1958} A.~Beurling and J.~Deny, \it Espaces de Dirichlet. 
I. Le cas \'el\'ementaire\rm, Acta Math., 99 (1958), 203--224.\ \par 
\ \par 
\bibitem [Den1970]{Den1970} J.~Deny, \it M\'ethodes hilbertiennes en th\'eorie 
du potentiel\rm ,  in: ``Potential Theory," Ed. M.~Brelot, Lectures 
at Summer School  Centro Internazionale Matematico Estivo 
(C.I.M.E.)  in Stresa (Varese), Italy, July 2-10, 1969; 1st. ed. 
C.I.M.E., Ed. Cremonese, Roma, 1970; reprint: Series C.I.M.E. Summer 
Schools,  49, Springer, 2011,  123--201.\ \par 
\ \par 
 \bibitem[DH1959]{DH1959} A.~Devinatz and I.~I.~Hirschman~Jr., \it Multiplier transformations on $\ell^{2,\alpha}$\rm, 
 Ann.  Math. (2), 69 (1959), 575--587.\ \par 
\ \par 
\bibitem [Duo2001] {Duo2001} J.~Duoandikoetxea, \it Fourier Analysis\rm, 
Providence, RI, AMS,  2001.\ \par 
\ \par 
\bibitem [FOT2011] {FOT2011} M.~Fukushima, Y.~Oshima, and M.~Takeda, \it Dirichlet Forms 
and Symmetric Markov Processes \rm  (2nd revised, extended ed.), 
Berlin--New York, De Gruyter, 2011.\ \par 
\ \par 
\bibitem [GRS1960]{GRS1960} I.~M.~Gelfand, D.~A.~Raikov, and G.~E.~Shilov, 
\it Commutative Normed Rings\rm, Fizmatgiz, Moscow, 1960 (Russian); 
English transl.: New York, Chelsea 1964.\ \par 
\ \par 
\bibitem[Gra2008]{Gra2008} L.~Grafakos, \it Classical Fourier Analysis \rm  (2nd 
ed.), Graduate Texts in Math., 249,  Heidelberg--New York, Springer, 2008.\ \par 
\ \par 
\bibitem[GMcG1979]{GMcG1979} C.~C.~Graham and O.~C.~McGehee, \it Essays 
in Commutative Harmonic Analysis\rm, Grundlehren der math. Wissenschaften, 238, 
Berlin--Heidelberg--New York, Springer, 
1979.\ \par 
\ \par 
\bibitem[Ja1985]{Ja1985} N.~Jacobson, \it Basic Algebra I \rm (2nd 
ed.),  New York, W.~H.~Freeman and Co., 
1985.\ \par 
\ \par 
\bibitem [Han1979] {Han1979} K.~Hansson, \it Imbedding theorems of Sobolev type 
in potential theory\rm , Math. Scand., 45 (1979), 77--102. \ \par 
\ \par 
\bibitem [KT1998] {KT1998} N.~J.~Kalton and L.~Tzafriri, \it The behaviour of Legendre and ultraspherical polynomials in $L_p$-spaces\rm, Canad. J. Math., 50 (1998), 1236--1252. \ \par 
\ \par 
\bibitem [Khi1934] {Khi1934} A.~Khinchin, \it Korrelationstheorie der station\"aren 
stochastischen Prozesse\rm, Math. Ann., 109:1 (1934), 604--615.\ \par 
\ \par 
\bibitem [Kre1944]{Kre1944} M.~G.~Krein, \it On the logarithm of an infinitely 
divisible Hermitian positive definite function\rm,  Doklady Acad. Nauk 
SSSR, 43:3 (1944), 99--102 (Russian).\ \par 
\ \par 
\bibitem [Lan1972]{Lan1972} N.~S.~Landkof, 
\it Foundations of Modern Potential Theory\rm,  Grundlehren der math. Wissenschaften, 180, 
Berlin--Heidelberg--New York, Springer, 1972.\ \par 
\ \par 
\bibitem [Lev1934] {Lev1934} P.~L\'evy, \it Sur les int\'egrales dont les 
\'el\'ements sont des variables al\'eatoires ind\'ependantes\rm, Annali 
delle Universit\`a Toscane, Pisa, (2), vol. 3 (1934), pp. 337--366. \ \par 
\ \par 
\bibitem [Lev1937] {Lev1937} P.~L\'evy, \it Th\'eorie de l'Addition des Variables 
Al\'eatoires\rm, Paris, 1937. \ \par 
\ \par 
\bibitem[Maz2011]{Maz2011} V.~Maz'ya, \it Sobolev Spaces, with Applications to Elliptic Partial Differential 
Equations\rm,  2nd augmented ed., Grundlehren der math. Wissenschaften, 342,  Berlin--New York, Springer, 2011.
\ \par 
\ \par 
\bibitem[MSh2009]{MSh2009} V.~G.~Maz'ya and T.~O.~Shaposhnikova, \it Theory of Sobolev Multipliers, with Applications to Differential and Integral Operators\rm,  Grundlehren der math. Wissenschaften, 337,  Berlin--New York, Springer, 2009.
\ \par 
\ \par 
\bibitem [vNS1941]{vNS1941}  J.~von~Neumann and I.~J.~Schoenberg, \it Fourier 
integrals and metric geometry\rm , Trans. Amer. Math Soc., 50 (1941), 226--251. 
\ \par 
\ \par 
\bibitem [Nik1986]{Nik1986} N.~K.~Nikolski, \it Treatise on the Shift Operator\rm, 
Grundlehren der math. Wissenschaften, 273,  Berlin--Heidelberg, Springer, 1986.
\ \par 
\ \par 
\bibitem [Nik2009]{Nik2009} N.~K.~Nikolski, \it The spectral localization property 
for diagonal operators and semigroups\rm , St. Petersburg Math. J., 21:6 
(2009), 1--25.\ \par 
\ \par 
\bibitem [Sch1938] {Sch1938} I.~J.~Schoenberg, \it Metric spaces and positive 
definite functions\rm , Trans. Amer. Math. Soc., 44 (1938), 522--536.\ \par 
\ \par 
\bibitem [Tor1986]{Tor1986} A.~Torchinsky,  \it Real-Variable Methods in 
Harmonic Analysis\rm , Academic Press, New York, 1986.

\end{thebibliography}
\end{document}